\documentclass[a4paper,12pt]{article}
\usepackage[T2A]{fontenc}
\usepackage[cp1251]{inputenc}
\usepackage[cp1251]{inputenc}
\usepackage[dvips]{graphicx}
\usepackage{amsfonts,amsmath,amssymb,amscd,array, euscript}
\begin{document}
\title{Affine and wrap algebras over octonions}
\author{S.V. Ludkovsky}
\date{10.12.2008}
\maketitle

\begin{abstract}
The article is devoted to affine and wrap algebras over quaternions
and octonions. Residues of functions of quaternion and octonion
variables are studied. They are used for construction of such
algebras. Their structure is investigated.

\end{abstract}

\section{Introduction}
Mention that Lie algebras are associative while the Cayley-Dickson
algebras ${\cal A}_r$ with $r\ge 3$ are non associative
\cite{baez,kansol,kurosh,ward}. This induces specific features of
algebras over ${\cal A}_r$. In this article the developed earlier
technique of residues of meromorphic functions of Cayley-Dickson
variables from \cite{luoyst,luoyst2} is used. It is necessary to
mention that theory of functions of Cayley-Dickson variables differ
drastically from that of complex and certainly has many specific
features. It is useful not only for mathematics, but also for
theoretical physics, including quantum mechanics, quantum field
theory, partial differential equations, non commutative geometry,
etc. \cite{emch,hamilt,harvey,lawmich,mensk}.
\par Over the complex field loop algebras are known. They are
defined for meromorphic functions in an open domain $U$ with one
singular marked point $z_0\in U\subset \bf C$ \cite{kacb}. But in
the case of Cayley-Dickson algebras it is possible to consider
meromorphic function with singularities in a closed subset $W$ of
codimension not less than 2. This $W$ may already be of dimension
greater than zero and winding around $W$ may be in any plane
containing $\bf R$. Thus there are winding surfaces around $W$, so
that the loop interpretation is lost. Therefore, analogs of loop
algebras over ${\cal A}_r$ are called here wrap algebras.
\par In this article affine and wrap algebras over quaternions and
octonions are introduced and studied. For this residues of functions
of quaternion and octonion variables are defined and their
properties are described. They are used for construction of such
algebras. Their structure is investigated. All main results of the
paper are obtained for the first time.

\section{Algebras over octonions}
\par To avoid misunderstandings we first introduce our notations and
definitions.
\par {\bf 1. Remark.} Let $V$ be a vector space over the Cayley-Dickson algebra
${\cal A}_r$. This means by our definition that $V=\mbox{ }_0V
i_0\oplus ... \oplus \mbox{ }_{2^r-1}V i_{2^r-1}$, where $\mbox{
}_0V,...,\mbox{ }_{2^r-1}V$ are pairwise isomorphic real vector
spaces, while $ \{ i_0,..., i_{2^r-1} \} $ is the set of the
standard generators of ${\cal A}_r$, $i_0=1$, $i_j^2=-1$,
$i_0i_j=i_j=i_ji_0$ and $i_ji_k = - i_ki_j$ for each $1\le j\ne k
\le 2^r-1$, $2\le r$, ${\cal A}_2=\bf H$ is the quaternion skew
field, ${\cal A}_3= \bf O$ is the octonion algebra. If in addition
each $\mbox{ }_jV$ is an associative real algebra and $(ax)(by) =
(ab)(xy)$ for each $x, y\in V_0$ and all $a, b\in {\cal A}_r$, then
we call $V$ the super-algebra over ${\cal A}_r$. For short we call
it also algebra.
\par Let $\sf g$ be a Lie super-algebra over the Cayley-Dickson
algebra ${\cal A}_r$, $r\ge 2$. By our definition this means that
${\sf g} = \mbox{ }_0 {\sf g} i_0\oplus \mbox{ }_1 {\sf g} i_1
\oplus ... \oplus \mbox{ }_{2^r-1}{\sf g} i_{2^r-1}$, where $\mbox{
}_0{\sf g},...,\mbox{ }_{2^r-1} {\sf g}$ are pairwise isomorphic
real Lie algebras, $\{ i_0, i_1,..., i_{2^r-1} \} $ are the standard
generators of ${\cal A}_r$. The multiplication in $\sf g$ is such
that \par $(1)$ $[xi_k, yi_j] = (-1)^{\chi (k,j) +1} [yi_j, xi_k]$
\\ for each pure states $xi_k\in {\sf g}_ki_k$ and $yi_j\in {\sf
g}_j i_j$, where $\chi (k,j)=0$ if either $k=0$ or $j=0$ or $k=j$,
while $\chi (k,j) =1$ for $k\ne j$ with $k\ge 1$ and $j\ge 1$. The
Jacobi identity takes the form:
\par $(2)$ $[xi_k, [yi_j, zi_s]] + (-1)^{\xi (k,j,s)} [yi_j, [zi_s,
xi_k]] + (-1)^{\xi (k,j,s) +\xi (j,s,k)} [zi_s, [xi_k, yi_j]] =0 $ \\
for each pure states $xi_k\in \mbox{ }_k{\sf g} i_k$, $yi_j\in
\mbox{ }_j{\sf g} i_j$ and $zi_s\in \mbox{ }_s{\sf g} i_s$, where
$\xi (k,j,s)\in \{ 0, 1 \} $ is such that $i_k(i_ji_s) = (-1)^{\xi
(k,j,s)} i_j(i_si_k)$. Moreover, the multiplication is real $\bf R$
bilinear:
\par $(3)$ $[ax, by] = ab [x,y]$ for each $x, y\in \sf g$ and
$a, b \in \bf R$;
\par $[x_1+x_2,y] = [x_1,y] + [x_2,y]$ and $[y,x_1+x_2] = [y,x_1] +
[y,x_2]$ for all $x_1, x_2, y\in {\sf g}$.
\par For short instead of Lie super-algebra we shall also write
Lie algebra or (Lie) algebra over ${\cal A}_r$.
\par Consider the family $Mat_n({\cal A}_r)$ of
$n\times n$ matrices with entries in ${\cal A}_r$.
\par If $X$ and $Y$ are vector spaces over ${\cal A}_r$ then we say
that a mapping $A: X\to Y$  is left ${\cal A}_r$ linear if it is
$\bf R$ linear and $A(z\mbox{ }_0x)= z A\mbox{ }_0x$ for each pure
vector $\mbox{ }_0x\in \mbox{ }_0X$ and every Cayley-Dickson number
$z$. The space of all left ${\cal A}_r$ linear operators from $X$
into $Y$ we denote by $L_l^a(X,Y)$. Suppose that $\sf h$ is a vector
space over ${\cal A}_r$.  Then we denote by ${\sf h}_l^*$ the space
$L_l^a({\sf h},{\cal A}_r)$ of all left ${\cal A}_r$ linear
functionals. Evidently ${\sf h}_l^*$ is the ${\cal A}_r$ vector
space.
\par We call vectors $v_1,...,v_n$ in a vector space $X$ over ${\cal
A}_r$ with $2\le r \le 3$ vector independent, if for each non zero
constants $a_1,...,a_n, b_1,...,b_n$ each vector $\{ a_1v_1b_1 \}
_{q_1(3)} +...+ \{ a_nv_nb_n \} _{q_n(3)} $ is non zero for each
associators $ \{ \alpha _1...\alpha _n \}_{q(n)}$ indicating on an
order of the multiplication, when $r=3$. For $r=2$ these associators
can be dropped, since the quaternion skew field is associative.
\par Consider a $n\times m$ matrix $B$ with entries in ${\cal A}_r$.
If $ 2\le r\le 3$, then ${\cal A}_r$ is alternative and each
equation $ax=b$ in it has the solution, which for $a\ne 0$ is
$x=a^{-1}b$. Therefore, the Gauss algorithm of reducing a matrix to
the step form is applicable in this case. It is seen from the step
form of the matrix $B$ that the rank of $B$ by rows and columns is
the same. Henceforth we consider $2\le r \le 3$ for algebras, if
another is not specified.
\par For example, Conditions $(1-3)$ are satisfied, when ${\sf g} =
\mbox{ }_0g \otimes {\cal A}_r$ and $[xa,yb] = [x,y]ab$ and $xa=ax$
for each $a, b\in {\cal A}_r$ and all $x, y\in \mbox{ }_0{\sf g}$,
where $\mbox{ }_0{\sf g}$ is the real Lie algebra.
\par {\bf 2. Definitions.} A matrix $A$ is called a generalized
Cartan matrix if \par $(C1)$ $a_{j,j} =2$ for each $j$,
\par $(C2)$ $a_{j,k}$ are non positive integers for all $j\ne k$,
\par $(C3)$ $a_{j,k}=0$ implies $a_{k,j}=0$.
\par By a realization of a matrix $A\in Mat_n({\cal A}_r)$ with $2\le
r \le 3$ we call the triple $({\sf h}, \Upsilon , \Upsilon ^{\vee
})$, where $\sf h$ is a vector space over ${\cal A}_r$, $\Upsilon =
\{ \beta _1,...,\beta _n \} \subset {\sf h}^*_l$, $\Upsilon ^{\vee }
= \{ \gamma _1,...,\gamma _n \} \subset \sf h$ so that
\begin{itemize}
\item[1] $\Upsilon $ and $\Upsilon
^{\vee }$ are ${\cal A}_r$ vector independent;
\item[2] $<\gamma _k, \beta _j> = a_{k,j}\in {\cal A}_r$, for all $k,
j=1,...,n$;
\item[3] $n-l=dim_{{\cal A}_r} {\sf h} - n$, where $l$ denotes the
rank of $A$ by rows over ${\cal A}_r$;
\item[4] $<a,b> := b(a) \in
{\cal A}_r$ for any $a\in \sf h$ and $b\in {\sf h}^*_l$ and $<\mbox{
}_0y, \beta _j> \in \bf R$ for each $j=1,...,n$ and $\mbox{ }_0y\in
\mbox{ }_0{\sf h}$; $~ <\mbox{ }_syi_s, \mbox{ }_k\beta i_k> =
(-1)^{\chi (s,k)} <\mbox{ }_k\beta i_k, \mbox{ }_s y i_s>$ for all
pure states.
\end{itemize}
\par Two realizations $({\sf h}, \Upsilon , \Upsilon ^{\vee } \} $
and $({\sf h}_1, \Upsilon _1, \Upsilon ^{\vee }_1 \} $ are called
isomorphic,  if there exists an isomorphism of vector spaces $\phi :
{\sf h}\to {\sf h}_1$ so that $\phi (\Upsilon ^{\vee }) = \Upsilon
^{\vee }_1$ and $\phi ^*(\Upsilon ) = \Upsilon _1$, where $<\phi ^*
(\beta ), \gamma > := <\beta , \phi (\gamma )>$ for all $\beta \in
\Upsilon $ and $\gamma \in \Upsilon ^{\vee }$.
\par {\bf 3. Proposition.} {\it  For each $n\times n$ matrix $A$ there exists a unique
up to an isomorphism its realization.  Realizations of matrices $A$
and $B$ are isomorphic if and only if $B$ can be obtained from $A$
by interchanging its rows or columns.}
\par {\bf Proof.} We can enumerate rows and columns of the matrix
and consider that $A= {{A_1}\choose {A_2}} ,$ where $A_1$ is the
$l\times n$ matrix of rank $l$. Compose the following matrix $C = {
{A_1\quad 0} \choose {A_2 I_{n-l}}}$ and put ${\sf h} = {\cal
A}_r^{2n-l}$, where $I_n$ denotes the unit $n\times n$ matrix. Take
$\beta _j(x)=x_j$ for $j=1,...,n$, where $x=(x_1,...,x_n)$, $x_j\in
{\cal A}_r$, $x\in {\cal A}_r^n$. Then as $\gamma _j$ take rows of
the matrix $C$. This gives the realization of $A$.
\par Vise versa, if for a realization $({\sf h}, \Upsilon , \Upsilon
^{\vee } )$ we complete $\Upsilon $ up to a basis with the help of
$\beta _{n+1},...,\beta _{2n-l}\in {\sf h}^*_l$. Then for suitable
$l\times (n-l)$ matrix $B$ and $(n-l)\times (n-l)$ matrix $D$ of
rank $(n-l)$ we get $(<\gamma _j, \beta _k>) = {{A_1 ~ B}\choose
{A_2 ~ D}}$. Adding to $\beta _{n+1}$ suitable ${\cal A}_r$ linear
combinations of $\beta _1,...,\beta _l$ we can get, that $B=0$,
where $2\le r\le 3$, since the octonion algebra is alternative and
each equation either $ax=b$ or $xa=b$ has the solution with non zero
$a$ either $x=a^{-1}b$ or $x=ba^{-1}$ respectively. Then substitute
$\beta _{n+1},...,\beta _{2n-l}$ on their ${\cal A}_r$ linear
combinations that to get $D=I$. This means the uniqueness of the
realization up to interchanging of rows and columns.
\par If $B$ is obtained from the matrix $A$ by interchanging its
rows and columns, then two realizations are evidently isomorphic,
since $x\mapsto (x_{\sigma (1)},...,x_{\sigma (n)})$ is the
isomorphism of ${\cal A}_r^n$, where $\sigma $ is a transposition,
that is a bijective surjective mapping of the set $ \{ 1,2,...,n \}
$.
\par {\bf 3.1. Remark.} Instead of using left ${\cal A}_r$ linear functionals
it is possible to use right ${\cal A}_r$ linear functionals in the
realization  of a matrix. Denote by ${\sf h}^*_r$ the space of right
${\cal A}_r$ linear functionals on $\sf h$, $\Upsilon _r\subset {\sf
h}^*_r$. If $({\sf h},\Upsilon , \Upsilon ^{\vee } )$ is the
realization of the matrix $A$, then $({\sf h}^*_r, \Upsilon ^{\vee
}, \Upsilon _r)$ is the realization of the transposed matrix $A^T$.
\par For two matrices $A_k$, $k=1, 2$, and their realizations $({\sf
h}_k, \Upsilon _k , \Upsilon ^{\vee }_k )$ we can get the
realization of their direct sum ${{A_1 ~ ~ 0}\choose {0 ~ ~ A_2}}$
of these matrices $({\sf h}_1\oplus {\sf h}_2, \Upsilon _1\times \{
0 \} \cup \{ 0 \} \times \Upsilon _2, \Upsilon ^{\vee }_1\times \{ 0
\} \cup \{ 0 \} \times \Upsilon ^{\vee }_2)$ it is called the direct
sum of their realizations.
\par The matrix $A$ is called decomposable, if after enumeration of
its columns and rows $A$ decomposes into a non trivial direct sum.
Clearly that $A$ can be presented as a direct sum of indecomposable
matrices and hence its representation as a direct sum of
indecomposable realizations. We call $\Upsilon $ the root basis,
while $\Upsilon ^{\vee }$ the dual root basis. Elements of $\Upsilon
$ or $\Upsilon ^{\vee }$ are called simple roots or dual simple
roots respectively. We put \par $Q= \sum_{j=1}^n {\bf Z}\beta _j$,
$Q_+ = \sum_{j=1}^n {\bf Z}_+ \beta _j$ \\
and call $Q=Q(A)$ the root lattice, where ${\bf Z}_+$ denotes the
set of all positive integers.
\par For $\beta = \sum_j k_j \beta _j\in Q$ the number $ht ~ \beta
:= \sum_j k_j$ is called the hight of the element $\beta $.
Introduce the partial ordering $\ge $ on ${\sf h}^*_l$ putting $a\ge
b$ if $a-b \in Q_+$.
\par {\bf 4. Definitions.} Take the matrix $A\in Mat_n({\cal A}_r)$ and its realization
$({\sf h}, \Upsilon , \Upsilon ^{\vee })$. We introduce the
auxiliary algebra $\eta (A)$ over ${\cal A}_r$  with generators
$e_j$, $f_j$, $j=1,...,n$ and $h\in \sf h$ and with defining
relations
\par $(1)$ $[e_k, f_j] = \delta _{k,j} \gamma _j$;
\par $(2)$ $[h, h']=0$ for all $h$, $h'\in \sf h$;
\par $(3)$ $[h, e_j] = <\beta _j,h> e_j$;
\par $(4)$ $[h, f_j] = - <\beta _j,h> f_j$ for all $j=1,...,n$ and
$h\in \sf h$; \par $(5)$ $[au,bv]= (ab) [u,v] = [u,v] (ab)$ for any
$a, b\in {\cal A}_r$ and $u, v\in \mbox{ }_0\eta (A)$.
\par Denote by $\eta _+$ and $\eta _-$ the subalgebra in $\eta (A)$
generated by elements $e_1$,...,$e_n$ and $f_1$,...,$f_n$
respectively.

\par A matrix $A$ is called decomposable, if after an enumeration of
its rows and columns it becomes the non-trivial direct sum of two
matrices.
\par {\bf 5. Theorem.} {\it Let $\eta (A)$ be as in Definition 4.
Then
\par $(1)$ $\eta (A) = \eta _- \oplus {\sf h} \oplus \eta _+$
is the direct sum of vector spaces over ${\cal A}_r$;
\par $(2)$ $\eta _+$ and $\eta _-$ are freely generated by
$e_1$,...,$e_n$ and $f_1$,...,$f_n$ correspondingly;
\par $(3)$ the mapping $e_j\mapsto - f_j$, $f_j\mapsto - e_j$ for each
$j=1,...,n$, $h\mapsto - h$ for all $h\in \sf h$ has a unique
extension up to an involution $\omega $ in the algebra $\eta (A)$;
\par $(4)$ there exists the decomposition into root spaces relative to
$\sf h$ of the form $\eta (A) = (\bigoplus_{\beta \in Q_+} \eta
_{-\beta }) \oplus {\sf h} \oplus ( \oplus _{\beta \in Q_+} \eta
_{\beta })$, where $\eta _{\beta } = \{ x\in \eta (A): ~ [h,x] =
\beta (h) x ~ \forall h\in {\sf h} \} $, moreover, $dim_{{\cal A}_r}
\eta (A) <\infty $ and $\eta _{\beta } \subset \eta _{\pm }$ for
$\pm \beta \in Q_+$, $\beta \ne 0$;
\par $(5)$ among ideals of $\eta (A)$ having the void intersection with
$\sf h$ a maximal ideal $\tau $ exists and $\tau = (\tau \cap \eta
_-) \oplus (\tau \cap \eta _+)$ is the direct sum of ideals.}
\par {\bf Proof.} Let $V$ be a vector space over ${\cal A}_r$ with
a basis $v_1,...,v_n$ so that $v_1,...,v_n\in \mbox{ }_0V$ and
$\beta \in \mbox{ }_0{\sf h}_l^*$. Define the tensor algebra $T(V)$
of $V$ as consisting of all elements $ \{ (b_1 x_1) ...(b_kx_k) \}
_{q(k)}$ and their finite sums, where $b_1,...,b_k\in {\cal A}_r$,
$x_1,...,x_k\in \{ v_1,...,v_n \} $, $k\in \bf N$, $q(s)$ is a
vector indicating on an order of the tensor multiplication
analogously to \cite{luoyst,luoyst2}. This means that $T(V) =
T^0(V)\oplus T^1(V)\oplus ...\oplus T^k(V)\oplus ...$, where $T^0(V)
= {\cal A}_r$, $T^1(V)=V$, $T^k(V)= T^{k-1}(V)\otimes V +
T^{k-2}\otimes T^2(V)+ ...+ V\otimes T^{k-1}(V)$, where the natural
equivalence relation in $T(V)$ is caused by the alternativity of the
octonion algebra $\bf O$ and associativity of the quaternion skew
field $\bf H$. For products of elements here the identities arising
from the fact that $\bf R$ is the center $Z({\cal A}_r)$ of the
algebra ${\cal A}_r$ can be taken into account. If $r=2$, then the
quaternion skew field ${\bf H} = {\cal A}_r$ is associative and in
this particular case associators $\{ * \} _{q(*)}$ are not
necessary.
\par Define the action of generators of the algebra $\eta (A)$ on
the tensor algebra $T(V)$ of the space $V$ by
\par $(6)$ $f_j(a) = v_j\otimes a$ for each $a\in T(V)$;
\par $(7)$ $h(1) = <\beta ,h>1$ and by induction
\par $h(v_j\otimes a) = - <\beta _j,h> v_j\otimes a + v_j\otimes
h(a)$ for each $a\in T^{s-1}(V)$, $j=1,...,n$, $h\in \eta_0$, $h\in
\mbox{ }_0{\sf h}$, $v_j \in \mbox{ }_0V$;
\par $(8)$ $e_j(1) = 0$ and by induction \par $e_k(v_j\otimes a) = \delta
_{k,j} \gamma _k(a) + v_j \otimes e_k(a)$ for every $a\in
T^{s-1}(V)$, $j=1,...,n$, $v_j\in \mbox{ }_0V$.
\par Verify now that Conditions 4$(1-5)$ are satisfied. Condition 4$(5)$ is satisfied,
since $\bf R$ is the center of the algebra ${\cal A}_r$. Relation
4$(2)$ is satisfied for each $h, h'\in {\sf h}_0$, since $<\beta
,h>\in \bf R$ for each $h\in \mbox{ }_0{\sf h}$. In view of 1$(3)$
it is fulfilled for all $h, h'\in \sf h$ as well. Take up to an
isomorphism $f_j\in \mbox{ }_0 \eta _-$ and $e_k\in \mbox{ }_0 \eta
_+$. Then from $(7,8)$ we infer that $(e_kf_j-f_je_k) (a) =
e_k(v_j\otimes a) - v_j \otimes e_k(a) = \delta _{k,j} \gamma _k(a)
+ v_j\otimes e_k(a) - v_j\otimes e_k(a) = \delta _{k,j} \gamma
_k(a)$, that is 4$(1)$ is satisfied. Then due to $(6,7)$  we get
\par $(hf_j-f_jh)(a) = h(v_j\otimes a) - v_j \otimes h(a) =
- <\beta _j,h> v_j\otimes a + v_j\otimes h(a) - v_j\otimes h(a)=
 - <\beta _j,h> f_j(a)$ this implies 2$(8)$. \par Mention that from
 $(6,7)$ for $s=0$ Condition 4$(3)$ follows. For $s>0$ we take
 $a=v_k\otimes b$ and $h\in \mbox{
}_0{\sf h}$, where $b\in T^{s-1}(V)$,
 so we deduce
 \par $(he_j - e_jh)(v_k\otimes b) = h(\delta _{j,k} \gamma _j(b)) +
 h(v_k \otimes e_j(b)) - e_j ( - <\beta _k,h> v_k\otimes b +
 v_k\otimes h(b)) = \delta _{j,k} \gamma _j(h(b)) - <\beta
 _k, h> v_k\otimes e_j(b) + v_k\otimes he_j(b) + <\beta _k, h> \delta
 _{j,k} \gamma _j(b) + <\beta _k,h> v_k\otimes e_j(b) - \delta
 _{j,k} \gamma _j h(b) - v_k \otimes e_jh(b) = <\beta _j,h> \delta
 _{j,k} \gamma _j(b) + v_k\otimes (he_j - e_jh) (b)$.
Apply the induction hypothesis to the second term, then we get
4$(3)$ for $\mbox{ }_0{\sf h}$, since $e_j\in \mbox{ }_0 \eta _+$.
Taking into account 2$(4)$ we get Condition 4$(3)$ for each $h\in
\sf h$.
\par In view of Conditions 4$(1-4)$ by induction prove that products
of elements from the  set $ \{ e_j, f_j: j=1,...,n; {\sf h} \} $
belong to $\eta _- + {\sf h} + \eta _+$. Consider an element $u=u_-
+ u_{\sf h} + u_+$, where $u_-\in \eta _-$, $u_+ \in \eta _+$,
$u_{\sf h} \in {\sf h}$. If $u=0$, then it acts in $T(V)$ such that
$u(1)=u_-(1) + <\lambda ,u_{\sf h}>=0$, hence $<\lambda ,u_{\sf
h}>=0$ for each $\lambda \in {\sf h}_l^*$, consequently, $u_{\sf
h}=0$.
\par Now use the mapping $f_j\mapsto v_j$, which makes the algebra
$T(V)$ as the enveloping algebra of the algebra $\eta _-$. The
tensor algebra $T(V)$ is free and is the universal enveloping
algebra $U(\eta _-)$ of $\eta _-$. This algebra is non-associative
for $r=3$, but it is associative for $r=2$. The mapping $u_- \mapsto
u_-(1)$ is the canonical embedding $\eta _-\hookrightarrow U(\eta
_-)$. Thus $u_- =0$ and the first statement is demonstrated.
\par If $A$ is an algebra over ${\cal A}_r$, then by $[A]$ we denote
the algebra obtained from $A$ by supplying it with the bracket
multiplication as above.
\par An algebra $U(L)$ over the Cayley-Dickson algebra ${\cal A}_r$
with the unit is called the universal enveloping algebra for an
algebra $L$ over ${\cal A}_r$ if there exists a homomorphism
$\epsilon : L\to [U]$ such that for each homomorphism $g: L\to [A]$
there exists a unique homomorphism $f: U\to A$ such that $g=f\circ
\epsilon $.
\par By our construction the algebra $T(V)$ has the decomposition $T(V) =
\mbox{ }_0T(V)i_0\oplus ... \oplus \mbox{ }_{2^r-1}T(V) i_{2^r-1}$,
where $\mbox{ }_0T(V),...,\mbox{ }_{2^r-1}T(V)$ are real pairwise
isomorphic algebras. Since $v_1,...,v_n\in \mbox{ }_0V$, then
$\mbox{ }_0T(V)$ is the universal enveloping algebra of $\mbox{ }_0
\eta _-$. By the Poincare-Birkhoff-Witt theorem and our choice of
$f_1,...,f_n$ above the algebra $\mbox{ }_0 \eta _-$ is free
generated by $f_1,...,f_n$. This makes Statement (2) evident.
\par Applying the involution $\omega $ we get, that $\eta _+$ is
free  generated by $e_1,...,e_n$. Using 4$(3,4)$ we get the
decomposition $\eta _{\pm } = \bigoplus _{\beta \in Q_+, \beta \ne
0} \eta_{\pm \beta }$. At the same time there is the estimate
$dim_{{\cal A}_r} \eta_{\beta } \le (2n)^{| ht ~ \beta |}$ for the
dimension of $\eta_{\beta }$ over ${\cal A}_r$. The latter implies
(4).
\par If $\tau $ is an ideal of $\eta(A)$, then it has the
decomposition $\tau = \mbox{ }_0\tau i_0 \oplus ... \oplus \mbox{
}_{2^r-1}\tau i_{2^r-1}$ with pairwise isomorphic real algebras
$\mbox{ }_0\tau ,...,\mbox{ }_{2^r-1}\tau $. Over $\bf R$ it appears
that $\mbox{ }_0\tau $ is the ideal of $\mbox{ }_0\eta(A)$. For real
algebras (5) is known. The ideal $\tau $ has the decomposition $\tau
= \bigoplus_{\beta } (\eta_{\beta }\cap \tau )$, hence $\tau =
\bigoplus_{\beta } [(\tau \cap \eta _- \cap \eta _{\beta }) \oplus
(\tau \cap \eta _+ \cap \eta _{\beta })] = (\tau \cap \eta _-)
\oplus (\tau \cap \eta _+)$ is the direct sum of ${\cal A}_r$ vector
spaces. Then $[f_j, \tau \cap \eta _+]\subset \eta _+$ and $[e_j,
\tau \cap \eta _+]\subset \eta _+$, consequently, $[\eta(A), \tau
\cap \eta _+] \subset \tau \cap \eta _+$, also $[\eta(A), \tau \cap
\eta _-] \subset \tau \cap \eta _-$. Thus the sum in (5) is the
direct sum of ideals.
\par {\bf 6. Note.} Consider the quotient algebra ${\sf g}(A) = \eta
(A)/ \tau $. It exists, since the pairwise isomorphic real quotient
algebras $(\eta (A))_j/\tau _j$ exist for each $j=0, 1,...,2^r-1$.
Each element of ${\sf g}(A)$ is of the form $b + \tau $, where $b\in
\eta (A)$. The matrix $A$ is called the Cartan matrix of the (Lie
super-) algebra ${\sf g}(A)$, while $n$ is called the rank of ${\sf
g}(A)$. The collection $({\sf g}(A), {\sf h}, \Upsilon , \Upsilon
^{\vee } )$ we call the quadruplet associated with the matrix $A$.
\par Two quadruplets $({\sf g}(A)^k, {\sf h}^k, \Upsilon ^k, \Upsilon
^{\vee ,k} )$, $k=1, 2$, are called isomorphic, if there exists an
isomorphism of algebras $\phi : {\sf g}(A)^1\to {\sf g}(A)^2$ such
that $\phi ({\sf h}^1) = {\sf h}^2$, $\phi ^*(\Upsilon ^1)= \Upsilon
^2$, $\phi (\Upsilon ^{\vee ,1}) = \Upsilon ^{\vee ,2}.$
\par We keep the notation $e_j,$ $f_j$, $\sf h$ for their images in ${\sf g}(A)$
The subalgebra $\sf h$ is called the Cartan subalgebra and its
elements $e_j, f_j$ are called Chevalley generators. They generate
the derivative sub-algebra which is by our definition ${\sf g}'(A)$
so that ${\sf g}'(A) = \mbox{ }_0{\sf g}'(A) i_0 \oplus ... \oplus
\mbox{ }_{2^r-1} {\sf g}'(A) i_{2^r-1}$ with pairwise isomorphic
real algebras $\mbox{ }_0{\sf g}'(A) \oplus ... \oplus \mbox{
}_{2^r-1} {\sf g}'(A)$, where $\mbox{ }_0 {\sf g}'(A) = [\mbox{
}_0{\sf g}(A), \mbox{ }_0 {\sf g}(A)]$ so that ${\sf g}(A) = {\sf
g}'(A) + {\sf h}$.
\par Mention that ${\sf g}(A) = {\sf g}'(A)$ if and only if $rank
(A)=n$ is maximal.
\par Then we put ${\sf h}' = \sum_{j=1}^n {\cal A}_r \gamma _j$,
hence ${\sf g}'(A) \cap {\sf h} = {\sf h}'$, ${\sf g}'(A) \cap {\sf
g}_{\beta } = {\sf g}_{\beta }$, if $\beta \ne 0$.
\par In view of Theorem 5 there exists the decomposition into root
spaces relative to $\sf h$ as
\par ${\sf g}(A) = \bigoplus_{\beta \in Q} {\sf g}_{\beta }$, \\
where ${\sf g}_{\beta } = \{ x\in {\sf g}(A): ~ [h,x] = \beta (h) x$
$\forall $ $h\in {\sf h} \} $ is the root sub-space corresponding to
$\beta $. Particularly, ${\sf g}_0 = {\sf h}$. The number $mult ~
\beta := dim_{{\cal A}_r} {\sf g}_{\beta }$ is called the dimension
of the element $\beta $. In view of \S 5 there is the estimate
\par $mult ~ \beta \le (2n)^{|ht ~ \beta |}$. \par An element $\beta \in Q$
is called a root, if $\beta \ne 0$ and $mult ~ \beta \ne 0$. If
either $\beta <0$ or $\beta >0$, then the root is called negative or
positive respectively. As usually denote by $\Delta =\Delta (A)$,
$\Delta _-$ and $\Delta _+$ the family of all elements, of all
positive elements, of all negative elements, then $\Delta $ is the
disjoint union of $\Delta _-$ and $\Delta _+$. \par So we get that
${\sf g}_{\beta } \subset \eta _+$ if $\beta >0$, ${\sf g}_{\beta
}\subset \eta _-$ if $\beta <0$. Moreover, ${\sf g}_{\beta }$ is
${\cal A}_r$ vector space spanned either on vectors $[...[[e_{j_1},
e_{j_2}], e_{j_3}]...e_{j_s}]$ with $\beta _1+...+\beta _s=\beta $
for $\beta >0$, or on vectors $[...[[f_{j_1}, f_{j_2}],
f_{j_3}]...f_{j_s}]$ with $\beta _1+...+\beta _s = - \beta $ for
$\beta >0$, since $e_j, f_j\in \mbox{ }_0{\sf g}$. Thus $g_{\beta
_j} = {\cal A}_r e_j$ for $\beta _j>0$, $g_{\beta _j} = {\cal A}_r
f_j$ for $\beta _j<0$, $g_{s\beta _j} =0$ for $|s|>1$. Each root is
either positive or negative and this implies the following.
\par {\bf 7. Lemma.} {\it If $\beta \in \Delta _+ \setminus \{
\beta _j \} $, then $(\beta + {\bf Z} \beta _j) \cap \Delta \subset
\Delta _+$.}
\par {\bf 8. Remark.} In accordance with Theorem 5$(3,5)$ the ideal
$\tau $ in $\eta (A)$ is $\omega $ invariant. Thus it induces the
automorphism of the algebra ${\sf g}(A)$ called its Chevalley
involution.  This means that $\omega (e_j) = - f_j$, $\omega (f_j) =
- e_j$, $\omega (h) = -h$ for each $h\in \sf h$. Moreover, $\omega
({\sf g}_{\beta }) = g_{-\beta }$, hence $mult ~ \beta = mult ~
(-\beta )$ and $\Delta _- = - \Delta _+$.
\par {\bf 9. Definitions.} Suppose that $G$ is a subset in the group ring
$\mbox{ }_0G i_0 \oplus_{\bf R} ... \oplus_{\bf R} \mbox{ }_{2^r-1}
G i_{2^r-1}$, where $\mbox{ }_jG$ are pairwise isomorphic
commutative groups and $ab = \sum_{j,k} \mbox{ }_ja \mbox{ }_kb$ for
all $a, b \in G$, $a = \mbox{ }_0a i_0 +...+\mbox{ }_{2^r-1}a
i_{2^r-1}$, $\mbox{ }_ja \in \mbox{ }_jG$ for every $j$. Let the
multiplication in $G$ be such that
\par $(1)$ there exists the unit element $e= \mbox{ }_0e i_0$, where
$\mbox{ }_0e$ is the unit element in $\mbox{ }_0G$;
\par $(2)$ for each $a\in G$ there exists $a^{-1}\in G$ with
$aa^{-1} = a^{-1}a=e$,
\par $(3)$ if $0\le r\le 2$, then the multiplication is associative,
if $r=3$, then the multiplication is alternative (weakly
associative) $(aa)b= a(ab)$, $b(aa) = (ba)a$ for all $a, b \in G$.
Then we call $G$ for short the quasi-commutative group instead of
algebraic associative quasi-commutative group for $r=2$ or
alternative quasi-commutative group for $r=3$.

 Let $G$ be a group. The decomposition $V =
\bigoplus_{g\in G} V_g$ of a vector space $V$ into the direct sum of
subspaces $V_g$ over ${\cal A}_r$ is called a $G$ gradation.
Elements of $V_g$ are called homogeneous of degree $g$. \par An
algebra $A$ having the decomposition $A = \bigoplus_{g\in G} A_g$
into subspaces $A_g$ such that $A_gA_h\subset A_{gh}$ for each $g,
h\in G$ is called $G$ graded. It is not supposed that $A$ is
associative or a Lie algebra. \par An algebra $A$ over ${\cal A}_r$
we call quasi-commutative, if its decomposition $A = \mbox{ }_0A
i_0\oplus ... \oplus \mbox{ }_{2^r-1}A i_{2^r-1}$ with pairwise
isomorphic real algebras $\mbox{ }_0A$,...,$\mbox{ }_{2^r-1}A$ is
such that each $\mbox{ }_jA$ is commutative.
\par {\bf 10. Proposition.} {\it  Suppose that $\sf h$ is a
quasi-commutative Lie algebra over ${\cal A}_r$, while $V$ is its
diagonalizable module, which means \par $(1)$ $V = \bigoplus_{b\in
\Lambda } V_b$, $\Lambda \subset {\sf h}^*_l$, where \par $V_b := \{
v\in V: ~ f(\mbox{ }_kv)=b(f)\mbox{ }_kv ~ \forall f\in {\sf h}, ~
\forall k=0,...,2^r-1 \} $, $\Lambda $ is a multiplicative subgroup
in ${\sf h}^*_l$.
\par Then each submodule $U$ in $V$ is graded relative to Gradation
$(1)$.}
\par {\bf Proof.} Consider an element $v\in V$ and write it in the form
$v =\sum_{j=1}^n v_j$ with $v_j \in V_{b_j}$ and $v_j = \mbox{
}_0v_j i_0+...+ \mbox{ }_{2^r-1}v_j i_{2^r-1}$ with $\mbox{
}_kv_j\in \mbox{ }_kV_j$. There exists $f\in {\sf h}$ so that the
numbers $b_j(f)$ are pairwise distinct. Then for $v\in U$ we get
\par $f^p(\mbox{ }_kv) = \sum_{j=1}^n b_j(f)^p \mbox{ }_kv_j \in U$
for each $p=0,1,...,n-1$ and $k=0,...,2^r-1$. We get the system of
linear algebraic equations over octonions with the non-degenerate
matrix. Using the gaussian algorithm we get, that all elements
$\mbox{ }_kv_j$ belong to $\mbox{ }_kU$ and hence $v_j\in U$ for
each $j$. The product of left linear functionals is also left
linear. Mention that $(\mbox{ }_kv i_k) (\mbox{ }_jv i_j) \in \mbox{
}_sV i_s$ for pure states $\mbox{ }_kv i_k\in \mbox{ }_kV i_k$,
where $i_s = i_k i_j$. Thus $U = \bigoplus_{b\in \Lambda } (U\cap
V_b)$ and this means that the submodule $U$ is $\sf h$ graded also.
\par {\bf 11. Note.} If $V$ is a graded vector space, then it can be
supplied with the formal (or in another words direct sum) topology
having a fundamental system of neighborhoods of zero consisting of
all open subsets in $V^F$, where $V^F := \bigoplus_{g\in G\setminus
F} V_g$, $F$ is a finite subset in $G$, when each $V_g$ is a
topological vector space. The completion of $V$ relative to such
topology is $\prod_{g\in G} V_g$ and the latter space is called its
formal (or direct product) completion.
\par A gradation can be introduced by counterpoising to generators $a_j$
elements $b_j\in G$ and putting by the definition that their degree
is $deg (a_j)=b_j$. This defines a unique $G$ graded Lie algebra
$\sf g$ if and only if the ideal of relations between $a_j$ is also
$G$ graded. This is the case for a free system of generators $a_j$,
$j=1,...,s$. \par Particularly let $s_1,...,s_n$ be integers and put
$deg ~ e_j = - deg ~ f_j =s_j$, $deg ~ {\sf h} =0$, then we get the
$\bf Z$ graded algebra ${\sf g}(A) = \bigoplus_{j\in \bf Z} {\sf
g}_j(s)$, where $s=(s_1,...,s_n)$, ${\sf g}_j(s) = \bigoplus_{\beta
= \sum_m k_m \beta _m \in Q; \sum_m k_m s_m = j} {\sf g}_{\beta }$.
Clearly if $s_j>0$ for each $j$, then ${\sf g}_0(s) = \sf h$ and
$dim_{{\cal A}_r} {\sf g}_j(s)<\infty $.
\par If $s=(1,...,1) = {\bf 1}$, then ${\sf g}_j({\bf 1}) =
\bigoplus_{\beta : ht (\beta ) =j } {\sf g}_{\beta }$ and ${\sf
g}_0({\bf 1}) = {\sf h}$, ${\sf g}_{-1}({\bf 1}) = \sum_m {\cal A}_r
f_m$, ${\sf g}_1({\bf 1}) = \sum_m {\cal A}_r e_m$ so that $\eta
_{\pm } = \bigoplus_{j\ge 1} {\sf g}_{\pm j} ({\bf 1})$.
\par {\bf 12. Lemma.} {\it Let $a\in \eta _+$ and either $[a,f_j]=0$ for
each $j=1,...,n$ or $[a,e_j]=0$ for all $j=1,...,n$, then $a=0$.}
\par {\bf Proof.} Suppose that $a\in \eta _+$ such that $[a, {\sf
g}_{-1}({\bf 1})]=0$, then $\sum_{k,m\ge 0} (ad ~ {\sf g}_1({\bf
1}))^k (ad ~ {\sf h})^m a$ is the subspace in $\eta _+$ over ${\cal
A}_r$ invariant relative to $ad ~ {\sf g}_1({\bf 1})$, $ad ~ {\sf
h}$ and $ad ~ {\sf g}_{-1}({\bf 1})$. Thus if $a\ne 0$ we will get
the ideal in ${\sf g}(A)$ which has with $\sf h$ the trivial
intersection. This contradicts the definition of ${\sf g}(A)$.
\par {\bf 13. Proposition.} {\it The center of the algebra ${\sf g}(A)$ or ${\sf g}'(A)$
considered over the real field $\bf R$ is $Z = \{ h\in \mbox{
}_0{\sf h}: <\beta _j, h>=0 ~ \forall ~ j=1,...,n \} $, moreover,
$dim_{\bf R} Z = n-l$.}
\par {\bf Proof.} The center of the Cayley-Dickson algebra ${\cal A}_r$ is
the field $\bf R$ of real numbers. We use the relations 1$(1-3)$.
Suppose that $p\in Z$ and $p=\sum_j p_j$ its decomposition relative
to the main gradation. Then $[p,{\sf g}_{-1}({\bf 1})]=0$ implies
$[p_j,{\sf g}_{-1}({\bf 1})]=0$ for each $j>0$. By Lemma 12 we get
that $p_j=0$ for all $j>0$. From $[p,{\sf g}_1({\bf 1})]=0$ we infer
that $p_j=0$ for each $j<0$. Thus $p\in \sf h$ and we have $[p,e_j]=
<\beta _j, p> e_j=0$, hence $<\beta _j, p> =0$ for each $j=1,...,n$.
\par Vice versa if $p\in \sf h$ and $<\beta _j, p> =0$ for each
$j=1,...,n$, then $p$ commutes with all Chevalley generators and
hence belongs to the center of the algebra. Mention also that
$Z\subset \mbox{ }_0\sf h$, since in the contrary case $dim_{\bf R}
Z >n-l$ and $\Upsilon $ will not be the ${\cal A}_r$ vector
independent set.
\par {\bf 14. Proposition.} {\it If a matrix $A$ is in $Mat_n({\cal A}_r)$,
then the algebra ${\sf g}_{{\cal A}_{r+1}}(A)$
over ${\cal A}_{r+1}$ is the smashed product of two copies of the
algebra ${\sf g}_{{\cal A}_r}(A)$ over ${\cal A}_r$ for $1\le r\le
2$.}
\par {\bf Proof.} Since $a_{j,k}\in {\cal A}_r$ for each $j, k$, then
using constant multipliers we can choose elements $\beta _j$ and
$\gamma _j$ belonging to $({\sf h}_{{\cal A}_r})^*_l$ and ${\sf
h}_{{\cal A}_r}$ respectively. The Cayley-Dickson algebra ${\cal
A}_{r+1}$ is the smashed product of two copies of ${\cal A}_r$ with
the help of the doubling procedure \cite{baez}. Therefore, each
element $v\in \eta _{{\cal A}_{r+1}}(A)$ can be decomposed in the
form $\mbox{ }_0v + \mbox{ }_1vi_{2^r}$ with $\mbox{ }_0v$ and
$\mbox{ }_1v\in \eta _{{\cal A}_r}(A)$, where $i_{2^r}$ is the
doubling generator of ${\cal A}_{r+1}$ from ${\cal A}_r$. This gives
the decomposition of $\eta _{{\cal A}_{r+1}}(A)$ as the vector space
over ${\cal A}_r$ into the direct sum $\eta _{{\cal A}_r}(A)\oplus
\eta _{{\cal A}_r}(A)i_{2^r}$, where $i_0,...,i_{2^{r+1}-1}$ are the
standard generators of ${\cal A}_{r+1}$. \par On the other hand each
pure state of $\eta _{{\cal A}_{r+1}}(A)$ is either a pure state of
$\eta _{{\cal A}_r}(A)$ or of $\eta _{{\cal A}_r}(A)i_{2^r}$ so that
all multiplication rules are the same in $\eta _{{\cal A}_{r+1}}(A)$
and in $\eta _{{\cal A}_r}(A)\oplus \eta _{{\cal A}_r}(A)i_{2^r}$.
Using the $\bf R$ bilinearity of the multiplication and the
decomposition of each element into a sum of pure states $\mbox{
}_jvi_j$ with $\mbox{ }_jv$ is belonging to the real algebra $\mbox{
}_j\eta _{{\cal A}_{r+1}}(A)$, $j=0,1,...,2^{r+1}-1$, we get that
the algebra $\eta _{{\cal A}_{r+1}}(A)$ over ${\cal A}_{r+1}$ is the
smashed product $\eta _{{\cal A}_r}(A)\otimes ^s \eta _{{\cal
A}_r}(A)$ of two copies of $\eta _{{\cal A}_r}(A)$ over ${\cal A}_r$
with the help of the doubling generator $i_{2^r}$.
\par The maximal ideal $\tau $ in $\eta _{{\cal A}_{r+1}}(A)$ having
with ${\sf h}_{{\cal A}_{r+1}}$
the trivial intersection possesses the same decomposition ${\tau
}_{{\cal A}_{r+1}} = {\tau }_{{\cal A}_r}\otimes ^s {\tau }_{{\cal
A}_r}$ so that real algebras $\mbox{ }_j {\tau }_{{\cal A}_{r+1}}$
and $\mbox{ }_k {\tau }_{{\cal A}_r}$ are pairwise isomorphic for
each $j, k$. This implies that $\mbox{ }_j{\sf g}_{{\cal
A}_{r+1}}(A)= \mbox{ }_j\eta _{{\cal A}_{r+1}}(A)/ \mbox{ }_j{\tau
}_{{\cal A}_{r+1}}$ for each $j$. Therefore, ${\sf g}_{{\cal
A}_{r+1}}(A) = \eta _{{\cal A}_{r+1}}(A)/\tau _{{\cal A}_{r+1}} =
(\eta _{{\cal A}_r}(A)/\tau _{{\cal A}_r}) \otimes ^s (\eta _{{\cal
A}_r}(A)/\tau _{{\cal A}_r}) = {\sf g}_{{\cal A}_r}(A) \otimes ^s
{\sf g}_{{\cal A}_r}(A)$.
\par {\bf 15. Lemma.} {\it Let $J_1$ and $J_2$ be two
nonintersecting subsets in $\{ 1,..., n \} $ so that
$a_{j,k}=a_{k,j}=0$ when $j\in J_1$ and $k\in J_2$. Suppose that
$\delta _s= \sum_{j\in J_s} k_{j,s} \beta _j$, $s=1, 2$ and $\beta =
\delta _1 + \delta _2$ is a root of algebra ${\sf g}(A)$ over ${\cal
A}_r$, $2\le r\le 3$. Then either $\delta _1$ or $\delta _2$ is
zero.}
\par {\bf Proof.} Let $j\in J_1$ and $k\in J_2$. Then $[\gamma _j,
e_k]=0$, $[\gamma _k, e_j]=0$, $[e_j, f_k]=0$, $[e_k, f_j]=0$,
consequently, $[e_j, e_k]=0$ and $[f_j, f_k]=0$ in accordance with
Lemma 12. Thus algebras ${\sf g}^1(A)$ and ${\sf g}^2(A)$ commute,
where ${\sf g}^s(A)$ denotes the algebra generated by $\{ e_j, f_j:
j\in J_s \} $. The algebra ${\sf g}_{\beta }(A)$ is contained in the
subalgebra generated by ${\sf g}^1(A)$ and ${\sf g}^2(A)$, hence
${\sf g}_{\beta }(A)$ is contained either in ${\sf g}^1(A)$ or in
${\sf g}^2(A)$.
\par {\bf 16. Proposition.} {\it The algebra ${\sf g}(A)$ over ${\cal A}_r$ with
$2\le r\le 3$ is simple if and only if \par $(1)$ $rank _{{\cal
A}_r} A=n$ and \par $(2)$ for each $1\le j, k \le n$ there exist
indices $j_1,...,j_s$ so that $(...(a_{j,j_1}
a_{j_1,j_2})...)a_{j_s,k}\ne 0$, where $n$ is the order of a matrix
$A$.}
\par {\bf Proof.} If either $(1)$ or $(2)$ is not satisfied, then
by Lemma 15 ${\sf g}(A)$ will contain a non trivial proper ideal
$\xi $, $\xi \ne 0$, $\xi \ne {\sf g}(A)$ and ${\sf g}(A)$ will not
be simple.
\par If Conditions $(1,2)$ are satisfied and $\delta $ is a non zero ideal in
${\sf g}(A)$, then $\delta $ contains a non zero element $h\in \sf
h$. Since $rank_{{\cal A}_r} =n$, then by Proposition 13 $Z=0$,
hence $[h,e_j]=ae_j\ne 0$ for some $j$ and inevitably $e_j\in \delta
$ and $\gamma _j = [e_j,f_j]\in \delta $. From Condition $(2)$ and
from the fact that the norm in the octonion algebra is
multiplicative we get that $e_j, f_j, \gamma _j\in \delta $ for all
$j$. On the other hand, Condition $(1)$ implies that $\sf h$ is the
${\cal A}_r$ vector span of elements $\gamma _j$. Thus we get
$\delta = {\sf g}(A)$.

\par {\bf 17. Definition.} A $n\times n$ matrix with entries in ${\cal A}_r$ is
called symmetrizable, if there exists a diagonal matrix $D=
(d_1,...,d_n)$, $d_j\ne 0$ for each $j$ and a symmetric matrix $B$,
that is $B^T=B$, so that \par $(1)$ $A=DB$, \\ where $B^T$ denotes
the transposed matrix $B$. In this case $B$ is called the
symmetrization of $A$ and an algebra ${\sf g}(A)$ is called
symmetrizable.
\par For a symmetrizable matrix $A$ with a given decomposition $(1)$
and its realization $({\sf h}, \Upsilon , \Upsilon ^{\vee })$ we fix
a complemented ${\cal A}_r$ vector space ${\sf h}_2$ to ${\sf h}_1 =
span_{{\cal A}_r} \{ \gamma _j : j \} $ in $\sf h$ a symmetric
${\cal A}_r$ valued form $(*|*)$ on $\sf h$ such that
\par $(2)$ it is left ${\cal A}_r$ linear by the left argument and right ${\cal A}_r$
linear by the right argument such that
\par $(3)$ $(\gamma _j|h) = d_j<\beta _j, h>$ for each $h\in \sf h$,
\par $(4)$ $(h_1|h_2) =0$ for all $h_1\in {\sf h}_1$ and $h_2\in {\sf h}_2$,
\par $(5)$ $(\mbox{ }_0x| \mbox{  }_0y) \in \bf R$ for each $\mbox{ }_0x,
\mbox{ }_0y\in \mbox{ }_0{\sf g}$, \par $(6)$ $(\gamma _j|\gamma _k)
=d_jb_{j,k}d_k$ with $b_{j,k}\in \bf R$ for each $j$ and $k$.

\par {\bf 18. Proposition.} {\it 1. The kernel of the restriction of
the form $(*|*)$ on ${\sf h}_1$ coincides with $Z$.
\par 2. The form $(*|*)$ is non degenerate on $\sf h$.}
\par {\bf Proof.} The first statement follows from Proposition 13.
To prove the second statement consider the condition $0= (\sum_j c_j
\gamma _j|h)$ for all $h\in \sf h$. Since $(\sum_j c_j \gamma _j|h)
= <\sum_j c_jd_j\beta _j,h>$ for each $h\in \mbox{ }_0{\sf h}$, then
$\sum_j c_jd_j\beta _j=0$, consequently, $c_j=0$ for each
$j=1,...,n$.
\par {\bf 19. Note.} The form $(*|*)$ is non degenerate and there is
the isomorphism $\nu : \mbox{ }_0h\to \mbox{ }_0{\sf h}^*_l$ having
the natural extension up to the isomorphism $\nu : {\sf h}\to {\sf
h}^*_l$ and $\nu _r: {\sf h}\to {\sf h}^*_r$ so that
\par $(1)$ $<\nu (h),p> = (h|p)$ for all $h, p\in \sf h$ and with
the  induced form $(*|*)$ on ${\sf h}^*_l$ so that \par $(2)$
$(\beta _j|\beta _k) = b_{j,k}=d_k^{-1}a_{k,j}$ for each $j, k$ and
\par $(3)$ $<p,\nu _r(h)> = (p|h)$ for each $h, p\in \sf h$ and with
the  induced form $(*|*)$ on ${\sf h}^*_r$. Clearly
\par $(4)$ $\nu (\gamma _j) = d_j \beta _j$ and
\par $(5)$ $\nu _r(\gamma _j) = \beta _jd_j$.

\par {\bf 20. Theorem.} {\it Let ${\sf g}(A)$ be a symmetrizable Lie
algebra and let 17$(1)$ be its prescribed decomposition. Then there
exists its non degenerate symmetric ${\cal A}_r$ valued form $(*|*)$
on ${\sf g}(A)$ satisfying Conditions 17$(2-5)$ and
\par $(1)$ this form $(*|*)$ is invariant on $\mbox{ }_0{\sf g}(A)$,
that is $([x,y]| z) = (x| [y,z])$ for all $x$, $y$ and $z\in \mbox{
}_0{\sf g}(A)$,
\par $(2)$ $({\sf g}_{\beta },{\sf g}_{\delta })=0$ if $\beta +
\delta \ne 0$,
\par $(3)$ the restriction $(*|*)|_{{\sf g}_{\beta } + {\sf g}_{ - \beta }}$ is
non degenerate for $\beta \ne 0$,
\par $(4)$ $[x,y] = (x|y) \nu ^{-1} (\beta )$ for each $x\in {\sf
g}_{\beta }$ and $y\in {\sf g}_{- \beta }$, $\beta \in \Delta $.}
\par {\bf Proof.} Take the principal $\bf Z$ gradation ${\sf g}(A) =
\bigoplus_{j\in \bf Z} {\sf g}_j$ and put ${\sf g}(m) =
\bigoplus_{j=-m}^m {\sf g}_j$ for $m=0, 1, 2,...$. Define the form
$(*|*)$ on ${\sf g}(0)=\sf h$ with the help of 17$(2-6)$. We extend
it on ${\sf g}(1)$ as
\par $(5)$  $(e_j|f_k) = \delta _{j,k} d_j,$  $j, k
=1,...,n$,
\par $(6)$ $({\sf g}_0|{\sf g}_{\pm 1 } ) =0$, $({\sf g}_1|{\sf
g}_1)=0$, $({\sf g}_{-1}|{\sf g}_{-1})=0$.
\par In view of 17$(2)$ $([e_j|f_k]|h)= (e_j|[f_k,h])$ for each
$h\in \sf h$, or equivalently $\delta _{j,k} (\gamma _j|h) = \delta
_{j,k} d_j <\beta _j,h>$, where we can take $e_j, f_k\in \mbox{
}_0{\sf g}$. Then the form $(*|*)$ satisfies Condition $(1)$, when
both elements $[x,y]$ and $[y,z]$ belong to ${\sf g}(1)$. Now we can
extend the form $(*|*)$ due to Rule 17$(2)$ on ${\sf g}(m)$ with the
help of induction by $m\ge 1$ so that $({\sf g}_j|{\sf g}_k)=0$, if
$|j|, |k|\le m$ and $j+k\ne 0$. Thus Condition $(1)$ is satisfied
when $[x,y]$ and $[y,z]$ belong to ${\sf g}(m)$.
\par Suppose that this extension is done for ${\sf g}(m-1)$, then we
should define $(x|y)$ for $x\in {\sf g}_{\pm m}$ and $y\in {\sf
g}_{\mp m}$ only. Write $y$ in the form $y=\sum_j [u_j,v_j]$, where
$u_j$ and $v_j$ are homogeneous elements of non zero degree,
belonging to ${\sf g}(m-1)$. Then $[x,u_j] \in {\sf g}(m-1)$ and for
$u_j, v_j\in \mbox{ }_0{\sf g}(m-1)$ we put
\par $(6)$ $(x|y) = \sum_j ([x,u_j]|v_j)$.
Then by Rule 17$(2)$ we extend it on ${\sf g}(m)$. It remains to
verify that the definition of Formula $(6)$ is correct on $\mbox{
}_0{\sf g}$. Suppose that $x_j\in \mbox{ }_0{\sf g}_j$ for $j\in \bf
Z$, take $k, j, s, t\in \bf Z$ with $|k+j| = |s+t|=m$ and $k+j+s+t
=0$; $|k|$, $|j|$, $|s|$ and $|t|<m$. Then
\par $(7)$  $([[x_k,x_j],x_s]|x_t) = (x_k|[x_j,[x_s,x_t]])$, \\
since $([[x_k,x_j],x_s]|x_t) = (([[x_k,x_s],x_j]|x_t) -
(([[x_j,x_s],x_k]|x_t)$ \\ $ = ([x_k,x_s]|[x_j,x_t]) +
(x_k|[[x_j,x_s],x_t])= (x_k|[x_s,[x_j,x_t]] + [[x_j,x_s],x_t]) = (x_k|[x_j,[x_s,x_t]])$. \\
If now $x= \sum_k [p_k,w_k]$, then from $(6,7)$ we get \par $(x|y) =
\sum_k
([x,u_k]|v_k) = \sum_k (p_k|[w_k,y])$. \\
Hence this value does not depend on a choice of expressions for $x$
and $y$.
\par From definitions it follows that Condition $(1)$ is fulfilled
as soon as $[x,y]$ and $[y,z]$ belong to ${\sf g}(m)$. Thus we have
constructed the form $(*|*)$ on $\sf g$ satisfying Conditions
$(1,2)$. Its restriction on $\sf h$ is non degenerate due to
Proposition 18.
\par The form $(*|*)$ also satisfies $(3)$, since $h\in \sf h$,
$x\in {\sf g}_{\alpha }$ and $y\in {\sf g}_{\beta }$. By the
invariance property and 17$(2)$ we infer \par $0 = ([h,x]|y) +
(x|[h,y]) = (<\beta ,h> + <\gamma ,h>)(x|y)$.
\par For $x\in {\sf g}_{\alpha }$ and $y\in {\sf g}_{ - \alpha }$
with $\alpha \in \Delta $ and $h\in \sf h$ we get
\par $([x,y] - (x|y)\nu ^{-1}(\alpha )|h) = (x|[y,h]) - (x|y)
<\alpha , h> =0$. Thus $(5)$ follows from $(2)$. From $(2,3,5)$ it
follows that the bilinear form $(*|*)$ is symmetric.
\par If $(4)$ is not satisfied, then by $(3)$ the form $(*|*)$ is
degenerate. Put ${\sf t} = Ker (*|*)$. This is ideal and by $(2)$ we
have ${\sf t} \cap {\sf h} =0$. But this contradicts to the
definition of ${\sf g}(A)$.
\par {\bf 21. Note.} Suppose that $A=(a_{k,j})$ is the symmetrizable
generalized Cartan matrix. Equation 17$(1)$ is equivalent with the
system of homogeneous linear equations and inequalities $d_j\ne 0$
so that $diag (d_1^{-1},...,d_n^{-1})A = (b_{j,k})$ with
$b_{j,k}=b_{k,j}$ for all $j, k$, since the octonion algebra ${\bf
O} = {\cal A}_3$ is alternative. But in the generalized Cartan
matrix all entries are integer, hence solutions can be chosen in the
field $\bf Q$ of rational numbers. Thus we can choose its
decomposition 17$(1)$ with $d_j>0$ so that $d_j$ and $b_{j,k}$ are
rational numbers.
\par We can suppose that $A$ is indecomposable. In view of
Proposition 16 for each $1<j\le n$ there exists a sequence
$1=j_1<j_2<...<j_{k-1}<j_k=j$ such that $a_{j_k,j_{k+1}}<0$.
Therefore, $a_{j_s,j_{s+1}} d_{j_{s+1}} = a_{j_{s+1},j_s}d_{j_s}$
for each $s=1,...,k-1$. Hence $d_jd_1>0$ for each $j$. Therefore, we
can choose $d_j>0$ for each $j$. If $A$ is indecomposable, then the
matrix $diag (d_1,...,d_n)$ is defined by 17$(1)$ uniquely up to the
multiplication on a constant. Now we fix a symmetric form satisfying
Conditions 17$(2-6)$ related with the decomposition 17$(1)$. In
accordance with Lemma 18 we get
\par $(1)$ $(\beta _j|\beta _j)>0$ for each $j=1,...,n$;
\par $(2)$ $(\beta _j|\beta _k)\le 0$ for each $j\ne k$;
\par $(3)$ $\gamma _j = 2\nu ^{-1}(\beta _j)/ (\beta _j|\beta _j)$.
Thus $a_{j,k} = 2(\beta _j|\beta _k)/ (\beta _j|\beta _j)$ for each
$j, k$. Then we take the extension of the form $(*|*)$ from $\sf h$
onto ${\sf g}(A)$ in accordance with Theorem 20. This form will be
called the standard invariant form. \par If choose for a root $\beta
$ dual bases $\{ e_{\beta }^j \} $ and $ \{ e_{-\beta }^j \} $ in
$\mbox{ }_0{\sf g}_{\beta }$ and in $\mbox{ }_0{\sf g}_{-\beta }$ so
that $(e_{\beta }^j|e_{-\beta }^k) =\delta _{j,k}$ for each $j, k
=1,...,mult ~ \beta $, then for $x\in {\sf g}_{\beta }$ and $y\in
{\sf g}_{-\beta }$ the identity
\par $(4)$ $(x|y) = \sum_j (x|e_{-\beta }^j) (y|e_{\beta }^j)$ is
satisfied.
\par {\bf 22. Lemma.} {\it If $\alpha , \beta \in \Delta $ and $z\in
{\sf g}_{\beta -\alpha }$, then in ${\sf g}(A)\otimes {\sf g}(A)$
the identity
\par $(1)$  $\sum_s e_{-\alpha }^s \otimes [z,e_{\alpha } ^s] = \sum_s
[e_{-\beta }^s,z]\otimes e_{\beta }^s$ is satisfied.}
\par {\bf Proof.} We define the form $(*|*)$
satisfying 17$(2-6)$ by the formula $(\mbox{ }_kxi_k\otimes \mbox{
}_jyi_j|\mbox{ }_swi_s\otimes \mbox{ }_tzi_t) = (-1)^{\zeta
(k,j,s,t)}(\mbox{ }_kxi_k|\mbox{ }_swi_s) (\mbox{ }_jyi_j|\mbox{
}_tzi_t)$ for pure states and extend it by $\bf R$ bi-linearity on
${\sf g}(A)\otimes {\sf g}(A)$, where $\zeta (k,j,s,t) \in \{ 0,1 \}
$ is such that $(i_ki_j)(i_si_t) = (-1)^{\zeta (k,j,s,t)} (i_ki_s)
(i_ji_t)$. Then we take $e\in \mbox{ }_0{\sf g}_{\alpha }$ and $f\in
\mbox{ }_0{\sf g}_{-\beta }$. Therefore,
\par $\sum _s (e_{-\alpha }^s \otimes [z,e_{\alpha }^s]|e\otimes f)
= \sum_s (e_{-\alpha }^s|e) ([z,e_{\alpha }^s]|f) = \sum_s
(e_{-\alpha }^s|e) (e_{\alpha
}^s|[f,z]) = (e|[f,z])$ also \\
\par $\sum_s ([e_{-\beta }^s,z]\otimes e_{\beta }^s|e\otimes f)
=\sum_s (e_{-\beta }^s|[z,e]) (e_{\beta }^s|f) = ([z,e]|f)$ \\
by Theorem 20 and Formula 21$(4)$. In view of 17$(2)$ the last two
formulas imply $(1)$, since the octonion algebra is alternative.
\par {\bf 23. Corollary.} {\it Let conditions of Lemma 22 be
satisfied, then
\par $(1)$ $\sum_s [e_{-\alpha }^s, [z, e_{\alpha }^s]] =  \sum_s
[[e_{-\beta }^s, z], e_{\beta }^s]$ in ${\sf g}(A)$;
\par $(2)$ $\sum_s e_{-\alpha }^s [z,e_{\alpha }^s] =  \sum_s
[e_{ -\beta }^s, z]e_{\beta }^s$ in $U({\sf g}(A))$.}
\par {\bf Proof.} This follows from the application of mappings
${\sf g}(A)\otimes {\sf g}(A)\ni x\otimes y \mapsto [x,y]\in {\sf
g}(A)$ and ${\sf g}(A)\otimes {\sf g}(A)\ni x\otimes y \mapsto xy\in
U({\sf g}(A))$ to Formula 22$(1)$.
\par {\bf 24. Remark.} Let ${\sf g}(A)$ be a Lie algebra
corresponding to a matrix $A$ and let $\sf h$ be its Cartan
subalgebra, ${\sf g} = \bigoplus_{\beta } {\sf g}_{\beta }$ be its
decomposition into root subspaces relative to $\sf h$. A ${\sf
g}(A)$ module (or ${\sf g}'(A)$ module) $V$ is called bounded if for
each $v\in V$ we have ${\sf g}_{\beta }(v)=0$ for all roots besides
a finite number of positive roots $\beta $.
\par Introduce the functional $\rho \in {\sf h}^*_r$ by the formula:
\par $(1)$ $<\rho ,\gamma _j> = a_{jj}/2$, $j=1,...,n$, where $\gamma _j\in
\Upsilon ^{\vee }$. If $rank ~ A <n$, then this does not define
$\rho $ uniquely, so we can take any functional satisfying these
relations. In accordance with Formulas 19$(2,4)$ we get
\par $(2)$ $(\rho |\beta _j) = (\beta _j|\beta _j)/2$ for each
$j=1,....,n$.
\par For each positive root $\beta $ we choose a basis $ \{ e_{\beta
}^j \} $ in the space ${\sf g}_{\beta }$ and take the dual basis $
\{ e_{-\beta }^j \} $ in ${\sf g}_{ - \beta }$ so that they belong
to $\mbox{ }_0{\sf g}_{\beta }$ and $\mbox{ }_0{\sf g}_{ - \beta }$
correspondingly. Then we put
\par $(3)$  $ \Omega _0 = 2 \sum_{\beta \in \Delta _+} \sum_j e_{-\beta }^j
e_{\beta }^j$ and this operator does not depend on a choice of the
dual basis. For each $v\in V$ only finite number of additives
$e_{-\beta }^j e_{\beta }^j(v)$ is non zero, hence $\Omega _0$ is
defined correctly on $V$. Let $u_1, u_2,...$ and $u^1, u^2,...$ be
dual bases of the subalgebra $\sf h$ belonging to $\mbox{ }_0\sf h$.
We define the generalized Casimir operator as
\par $(4)$ $\Omega = 2\nu ^{-1}_r(\rho ) + \sum_j u^ju_j + \Omega
_0$. Choose $\rho \in \mbox{ }_0{\sf h}^*_r$.
\par We consider now the decomposition into root subspaces for
$U({\sf g}(A))$ relative to $\sf h$:
\par $U({\sf g}(A)) = \bigoplus_{\beta \in Q} U_{\beta }$, where
\par $U_{\beta } = \{ x\in U({\sf g}(A)): ~ [h,x] = <\beta ,h> x ~
\forall h\in {\sf h} \} $ \\ and put ${U'}_{\beta } = U({\sf
g}'(A))\cap U_{\beta }$ so that $U({\sf g}'(A)) = \bigoplus_{\beta }
{U'}_{\beta }$. \par {\bf 25. Theorem.} {\it Let ${\sf g}(A)$ be a
symmetrizable algebra over ${\cal A}_r$, $2\le r \le 3$. If $V$ is a
bounded ${\sf g}'(A)$ module and $u \in {U'}_{ \beta }$, then
\par $(1)$  $[\Omega _0, u] = - (2(\rho |\beta ) + (\beta |\beta )
+ 2\nu _r^{-1}(\beta )) u$. \par If $V$ is a bounded ${\sf g}'(A)$
module, then
\par $(2)$ $\Omega $ commutes with the action of ${\sf g}(A)$ on
$V$.}
\par {\bf Proof.} Since $\lambda = \sum_j <\lambda ,u^j>\nu _r(u_j)
= \sum_j <\lambda ,u_j> \nu (u^j)$ for each $\lambda \in {\sf
g}(A)$, then
\par $(3)$ $\sum_j <\lambda ,u^j> <u_j,\mu > = (\lambda |\mu )$.
Moreover,
\par $(4)$ $[\sum_j u^ju_j,x] = ((\beta |\beta ) + 2\nu
_r^{-1}(\beta ))x$ for each $x\in {\sf g}_{\beta }$, since
\par $[\sum_j u^ju_j,x] = \sum_j <\beta ,u^j>xu_j + \sum_j u^j
<u_j,\beta >x$ \\
$=\sum_j <\beta , u^j>  <u_j,\beta > x + (\sum_j u^j<u_j, \beta > +
<\beta ,u^j>u_j)x$.
\par Thus Statement $(2)$ follows from $(1)$ and Formula $(4)$.
\par Take now elements $e_{\beta _j}$ and $e_{-\beta _j}$ with
$j=1,...,n$, which generate the algebra ${\sf g}'(A)$. If either
$u=e_{\beta _j}$ or $u=e_{-\beta _j}$, then due to Lemmas 7 and 23
we infer:
\par $[\Omega _0,e_{\beta _j}] = 2 \sum_{\beta \in \Delta _+} \sum_s
([e_{-\beta }^s,e_{\beta _j}]e_{\beta }^s + e_{-\beta }^s [e_{\beta
}^s,e_{\beta _j}]$
\par $ = 2 [e_{-\beta _j},e_{\beta _j}]e_{\beta _j} + 2 \sum_{\beta
\in \Delta _+\setminus \{ \beta _j \} } (\sum_s [e_{-\beta
}^s,e_{\beta _j}] e_{\beta }^s + \sum_s e_{-\beta +\beta _j}^s
[e^s_{\beta - \beta _j}, e_{\beta _j}]) $
\par $ = - 2 \nu _r^{-1}(\beta _j) e_{\beta _j} = - 2(\beta _j|\beta
_j) e_{\beta _j} - 2 e_{\beta _j} \nu _r^{-1} (\beta _j)$.
\par Analogously we have $[\Omega _0, e_{-\beta _j}] = 2e_{-\beta
_j}[e_{\beta _j},e_{-\beta _j}] = 2e_{-\beta _j} \nu _r^{-1}(\beta
_j)$. Thus we have got $(1)$ for $u=e_{\beta _j}$ and $u=e_{-\beta
_j}$. If $u\in {U'}_{\alpha }$ and $v\in {U'}_{\beta }$, then
\par $[\Omega _0,uv] = [\Omega _0,u]v + u[\Omega _0,v]$\par $ = - \{ (2(\rho
|\alpha ) + (\alpha |\alpha ) + 2\nu ^{-1}_r(\alpha ))u \} v - u\{
(2(\rho |\beta ) + (\beta |\beta ) + 2\nu _r^{-1}(\beta ))v \} $
\par $= - ((2(\rho |\alpha ) + (\alpha |\alpha ) + 2\nu
_r^{-1}(\alpha ) + 2(\alpha |\beta ) + 2(\rho |\beta ) + (\beta
|\beta ) + 2\nu _r^{-1}(\beta )) uv$
\par $= (2(\rho |\alpha + \beta ) + (\alpha + \beta |\alpha + \beta
) + 2\nu _r^{-1} (\alpha + \beta )) uv$, \\
since $\alpha , \beta \in \mbox{ }_0{\sf h}^*$, $e_{\beta } \in
\mbox{ }_0{\sf g}_{\beta }$,  $\rho \in \mbox{ }_0{\sf h}^*_r$ and
hence $\nu _r^{-1}(\rho ) \in \mbox{ }_0{\sf h}$, while $\bf R$ is
the center of the Cayley-Dickson algebra ${\cal A}_r$.
\par Thus $(1)$ is proved in general also.
\par {\bf 26. Corollary.} {\it If suppositions of Theorem 25$(2)$
are satisfied and there exists a vector $v\in V$ so that $e_j(v)=0$
for all $j=1,...,n$ and $h(v)=<b,h>v$ for some $b\in {\sf h}^*_r$
and all $h\in \sf h$, then
\par $(1)$ $\Omega (v) = (b+2\rho |b)v$. \par Moreover, if $U({\sf
g}(A))v=V$, then
\par $(2)$ $\Omega = (b+2\rho |b) I_V$.}
\par {\bf Proof.} Formula $(1)$ follows from the definition of
$\Omega $ and Formula 25$(3)$. Then Formula $(2)$ follows from
Formula $(1)$ and Theorem 25.
\par {\bf 27. Proposition.} {\it Let $\sf g$ be an algebra over
${\cal A}_r$, $2\le r\le 3$, with an invariant form $(*|*)$
satisfying Conditions 17$(2-6)$, let also $ \{ x_j: j \} $ and $ \{
y_j: j \} $ be dual bases belonging to $\mbox{ }_0{\sf g}$, that is
$(x_j|y_k) = \delta _{j,k}$ for each $j, k$. Suppose that $V$ is an
$\sf g$ module so that for each pair of elements $u, v\in V$ either
$x_j(u)=0$ or $y_j(v)=0$ for all $j$ besides a finite number of $j$.
Then the operator
\par $\Omega _2 := \sum_j x_j\otimes y_j$ \\
is defined on $V\otimes V$ and commutes with actions of all elements
of $\sf g$ on $V$.}
\par {\bf Proof.} Consider commutators $[z,x_j] = \sum_k c_{j,k}
x_k$, $[z,y_j] = \sum_k p_{k,j} y_k$, where $c_{j,k}, p_{k,j}\in
{\cal A}_r$. Taking the scalar products we deduce that $c_{j,k} =
([z,x_j]|y_k)$ and $p_{k,j} = ([z,y_j]|x_k)$. From the invariance of
$(*|*)$ it follows that $c_{j,k} = (z|[x_j,y_k])$ and $p_{k,j} =
(z|[y_j,x_k])$, consequently, $c_{j,k} = - p_{k,j}$ for all $j, k$
and inevitably $\sum_j ([z,x_j]\otimes y_j + x_j\otimes [z,y_j])=0$,
since $x_j, y_j \in \mbox{ }_0{\sf g}$ for each $j$.
\par {\bf 28. Example.} Consider the $n\times n$ zero matrix $A=0$ with either $n\in
\bf N$ or $n=\infty $. This means that $[e_j,e_k]=0$, $[f_j,f_k]=0$,
$[e_j,f_k]=\delta _{j,k}\gamma _j$ for all $j, k=1,...,n$, hence
${\sf g}(0) = {\sf h} \oplus \sum_j {\cal A}_r e_j \oplus \sum_j
{\cal A}_r f_j$. The center of ${\sf g}(0)$ is $Z = \sum_j {\bf R}
\gamma _j$. Moreover, $dim_{{\cal A}_r} {\sf h} =2n$ and we can
choose elements $d_1,...,d_n\in {\sf h}$ so that
\par ${\sf h} = {\cal A}_r\otimes _{\bf R} Z + \sum_j {\cal A}_r d_j$ and $[d_j,e_k] = \delta
_{j,k} e_k$, $[d_j,f_k] = - \delta _{j,k} f_k$ for each $j,
k=1,...,n$. A non degenerate symmetric invariant form satisfying
Conditions 17$(2-6)$ is defined as $(e_j|f_j) = 1$, $(\gamma _j|d_j)
= 1$ for each $j$ and all others scalar products are zero. In the
considered algebra $\rho =0$ and the Casimir operator takes the form
\par  $ \Omega = 2\sum_j \gamma _jd_j + 2 \sum_j f_je_j$.
\par We put $q = \sum {\cal A}_r (\gamma _j - \gamma _k)\subset {\cal A}_r\otimes_{\bf R} Z$ and the
algebra $H := {\sf g}'(0)/q$ we call the Heisenberg algebra over the
Cayley-Dickson algebra ${\cal A}_r$ so that $\mbox{ }_0H = \mbox{
}_0{\sf g}'(0)/\mbox{ }_0q$, $H = \mbox{ }_0Hi_0\oplus ... \oplus
\mbox{ }_{2^r-1}H i_{2^r-1}$ with pairwise isomorphic real algebras
$\mbox{ }_0H,...,\mbox{ }_{2^r-1} H$, where ${\sf g}' = \mbox{
}_0{\sf g}' i_0 \oplus ... \oplus \mbox{ }_{2^r-1}{\sf g}'
i_{2^r-1}$.

\section{Residues of octonion meromorphic functions}
\par For subsequent proceedings we need residues of (super)
differentiable functions of Cayley-Dickson variables. In this
section they are studied in more generality and with new details as
in previous papers.
\par {\bf 1. Definitions.}
\par Let $f: V\to {\cal A}_r$ be a function, where $V$ is a neighborhood
of $z\in {\hat {\cal A}}_r$, where ${\hat {\cal A}}_r$ is the
one-point compactification of ${\cal A}_r$ with the help of the
infinity point or as the non commutative analog of the $2^r$
dimensional Riemann sphere \cite{luoyst,luoyst2}. Then $f$ is said
to be meromorphic at $z$ if $f$ has an isolated singularity at $z$
and $f$ is ${\cal A}_r$-holomorphic in $V\setminus \{ z \} $. If $U$
is an open subset in ${\hat {\cal A}}_r$, then $f$ is called
meromorphic in $U$ if $f$ is meromorphic at each point $z\in U$. If
$U$ is a domain of $f$ and $f$ is meromorphic in $U$, then $f$ is
called meromorphic on $U$.
\par This definition has the natural generalization. Let
$W$ be a closed connected subset in ${\hat {\cal A}}_r$ and its
codimension $codim (W)= 2^r - dim (W)\ge 2$ and $W\cap ( z +
{\widehat{ {\bf R}\oplus M{\bf R}}})$ is a set consisting of
isolated points for each purely imaginary Cayley-Dickson number
$M\in {\cal A}_r$ with $|M|=1$ and every $z\in W\cup \{ 0 \} $,
where $dim (W)$ is the topological covering dimension of $W$ (see \S
7.1 \cite{eng}). If there exists an open neighborhood $V$ of $W$ so
that a function $f$ is ${\cal A}_r$-holomorphic in $V\setminus W$
and may have singularities at points of $W$, then $f$ is called
meromorphic at $W$. If $U$ is an open subset in ${\hat {\cal A}}_r$
and different subsets $W_k$ in $U$ are isolated from each other,
that is $\inf \{ |\zeta -\eta |: j=1,2,k-1,k+1,...; \zeta \in W_k,
\eta \in W_j \}
>0$ for each marked $k$ and $f$ is meromorphic at each $W_k$ and
${\cal A}_r$-holomorphic in $U\setminus [\bigcup_k W_k]$ and
$\bigcup_k W_k$ is closed in ${\cal A}_r$, then $f$ is called
meromorphic in $U$.
\par Mainly we shall consider meromorphic functions with point isolated
singularities if another will not be specified and denote by ${\bf
M}(U)$ the set of all meromorphic functions on $U$ with singleton
singularities $W= \{ z \} $. Let $f$ be meromorphic on a region $U$
in the set ${\hat {\cal A}}_r$. A point
$$c \in \bigcap_{V\subset U, V \mbox{ is closed and bounded }}
cl (f(U\setminus V)) $$ is called a cluster value of $f$.

\par Let $V$ be an open subset in ${\cal A}_r$. Define the residue of a
meromorphic function $f$ at $W$ with a singularity at a point $a\in
W\subset {\cal A}_r$ as $$(i)\quad Res (a,f).M := (2\pi
)^{-1}\lim_{y \to 0}(\int_{\gamma _y}f(z)dz)$$ whenever this limit
exists,
$$\mbox{where }~ \gamma _y(t)=
a+ y \rho \exp (2\pi tM)\subset V\setminus W,$$ $\rho >0$, $|M|=1$,
$M\in {\cal I}_r$, $t\in [0,1]$, $0<y\le 1$, $\gamma := \gamma _1$
encompasses only one singular point $a$ of $f$ in the complex plane
$(a+ {\bf R}\oplus M{\bf R})$, ${\cal I}_r := \{ z: ~ z\in {\cal
A}_r, Re (z)=0 \} $, $Re (z) := (z+{\tilde z})/2$. Here as usually
we suppose that a function \par $(R1)$ $f$ is ${\cal
A}_r$-holomorphic on $V\setminus W$, \par $(R2)$ $W$ is a closed
connected subset in ${\hat {{\cal A}}}_r$ of codimension not less
than $2$,
\par $(R3)$ the intersection $W\cap (z+\widehat{{\bf R}\oplus N{\bf R}})$
consists of isolated points for each purely imaginary Cayley-Dickson
number $N\in {\cal I}_r$, $|N|=1$, and every $z\in W\cup \{ 0 \} $.

\par Extend $Res (a,f).M$ by Formula $(i)$ on ${\cal I}_r$ as
$$(ii)\quad Res (a,f).M := [Res (a,f).(M/|M|)] |M|,$$
$\forall M\ne 0$; $Res (a,f).0 := 0$,  when $Res (a,f).M$ is finite
for each $M\in {\cal I}_r$, $|M|=1$.

\par {\bf 2. Definition.} For a metric space $X$ with a metric
$\rho $ let $dist (x,A) := \inf \{ \rho (x,y): y\in A \} $ denotes a
distance from a point $x$ to a subset $A$ in $X$. Let $z$ be a
marked point in the Cayley-Dickson algebra and $\gamma (t) := z + R
\exp (2\pi Mt)$ be a circle with center at $z$ of radius $0<R<\infty
$ in the plane $z+ ({\bf R}\oplus M{\bf R})$ embedded into ${\cal
A}_r$, where $t\in [0,1]$. Then $\gamma $ encompasses $z$ in the
usual sense. \par We say that a loop $\psi (t)$ in the
Cayley-Dickson algebra ${\cal A}_r$ encompasses the point $z$
relative to $\gamma $, if there exists a continuous monotonously
increasing function $\phi : [0,1]\to [0,1]$ being piecewise
continuously differentiable so that $|\psi (t) - \gamma (\phi (t))|
< \min ( |\psi (t)-z|,~ R)$ for each $t\in [0,1]$.
\par If ${\cal A}_m$ is the subalgebra of ${\cal A}_r$, $m<r$, then
there exists the projection $P_m: {\cal A}_r\to {\cal A}_m$ as for
real linear spaces.
\par {\bf 3. Theorem.} {\it Let $f$ be an ${\cal A}_r$-holomorphic
function on an open domain $U$ in ${\cal A}_r$, $\infty \ge r\ge 2$.
If $(\gamma +z_0)$ and $\psi $ are presented as piecewise unions of
paths $\gamma _j+z_0$ and $\psi _j$ with respect to parameter
$\theta \in [a_j,b_j]$ and $\theta \in [c_j,d_j]$ respectively with
$a_j<b_j$ and $c_j<d_j$ for each $j=1,...,n$ and
$\bigcup_j[a_j,b_j]=\bigcup_j[c_j,d_j]=[0,1]$ homotopic relative to
$U_j\setminus \{ z_0 \} $, where $U_j\setminus \{ z_0 \} $ is a
$(2^r-1)$-connected open domain in the Cayley-Dickson algebra ${\cal
A}_r$ such that $\pi _{s,p,t}(U_j\setminus \{ z_0 \})$ is simply
connected in $\bf C$ for each $s=i_{2k}$, $p=i_{2k+1}$,
$k=0,1,...,2^{r-1}-1$ ($\forall 0\le k\in {\bf Z}$ and
$P_m(U_j\setminus \{ z_0 \} )$ is $(2^m-1)$-connected for each $4\le
m\in \bf N$ if $r=\infty $), each $t\in {\cal A}_{r,s,p}$ and $u\in
{\bf C}_{s,p}$ for which there exists $z=t+u\in {\cal A}_r$. If
$(\gamma +z_0)$ and $\psi $ are closed rectifiable paths (loops) in
$U$ such that $\gamma (\theta )=\rho \exp (2\pi \theta M)$ with
$\theta \in [0,1]$ and a marked $M\in {\cal I}_r$, $|M|=1$ and
$0<\rho <\infty $, also $z_0\notin \psi $ and $M = [\int_{\psi } d
Ln (\zeta -z)]/(2\pi )$ and $\psi $ encompasses $z_0$ relative to
$\gamma +z_0$. Then

$$(1)\quad (2\pi )f(z)M=\int_{\psi }f(\zeta )dLn (\zeta -z) $$

for each $z\in U$ such that $$|z-z_0|< \inf_{\zeta \in \psi ([0,1])}
|\zeta -z_0|.$$  If either ${\cal A}_r$ is alternative, that is,
$r=2, 3$, or $f(z)\in {\bf R}$ for each $z$, then

$$(3.11)\quad f(z)=(2\pi )^{-1}(\int_{\psi }f(\zeta )d Ln (\zeta -z)^{-1}
) M^*.$$}

\par {\bf Proof.} The logarithmic function is $\zeta $-differentiable
so that $dLn (\zeta -z) = [D_{\zeta } Ln (\zeta -z)].d\zeta $ for a
marked $z$ and the variable $\zeta $ with $\zeta \ne z$,
consequently, the considered integrals exist. Take a Cayley-Dickson
number $z$ satisfying the conditions of the theorem. Put $\psi
_y(\theta ) := z + y(\psi (\theta ) -z)$ and $\gamma _y(\theta )= y
\rho \exp (2\pi \theta M)$ for each $\theta \in [0,1]$ and every
$0<y\le 1$. Up to a notation we can consider $\phi (\theta )=\theta
$ for each $\theta $ denoting $\gamma \circ \phi $ for simplicity by
$\gamma $. Then $|\psi _y(\theta ) -(\gamma _y(\theta )+z)|< |\psi
_y(\theta ) -z|$ for each $0<y\le 1$.
\par Join paths $\gamma _y +z$ and $\psi _y$ by a
rectifiable path $\omega _y$ such that $z\notin \omega $, which is
going in one direction and the opposite direction, denoted $\omega
_y ^-$, such that $\omega _{y,j}\cup \psi _{y,j}\cup \gamma
_{y,j}\cup \omega _{y,j+1}$ is homotopic to a point relative to
$U_j\setminus \{ z_0 \}$ for suitable $\omega _{y,j}$ and $\omega
_{y,j+1}$, where $\omega _{y,j}$ joins $\gamma _y(a_j)+z$ with $\psi
_y(c_j)$ and $\omega _{y,j+1}$ joins $\psi _y(d_j)$ with $\gamma
_y(b_j)+z$ such that $z$ and $z_0\notin \omega _j$ for each $j$.
\par  Then the equality
$$\int_{\omega _{y,j}} f(\zeta ) d Ln (\zeta -z) = - \int_{\omega
_{y,j}^-}f(\zeta ) d Ln (\zeta -z) $$ is accomplished for each $j$.

\par Mention that $[D_{\zeta } Ln (\zeta -z)].1 = (\zeta -z)^{-1}$
for each $\zeta \ne z$ and the function $f(\zeta )$ is ${\cal
A}_r$-holomorphic in $U$. The branching of the logarithmic function
$Ln (\zeta -z)$ by the variable $\zeta $ with $\zeta \ne z$ for a
marked Cayley-Dickson number $z$ is independent from $|\zeta -z|>0$
as shows its non-commutative Riemannian surface described in
\cite{luoyst,luoyst2,norfamlud}, since $Ln (y(\zeta -z)) = \ln y +
Ln (\zeta -z)$ for each $y>0$ and each $\zeta \ne z$, so that
branches of $Ln (\zeta -z)$ are indexed by purely imaginary
Cayley-Dickson numbers $M\in {\cal I}_r$. Consider now sub-domains
$U_j\setminus \{ z_0 \} $ and loops there composed of fragments
$\psi ([c_j,d_j])$ and $\psi _y([c_j,d_j])$ and joining their
corresponding ends rectifiable paths gone in a definite direction,
as well as loops composed of fragments $\gamma ([a_j,b_j]) +z$ and
$\gamma _y([a_j,b_j]) +z$ and joining their respective ends
rectifiable paths gone in a definite direction. Therefore, due to
the homotopy Theorem 2.15 \cite{luoyst,luoyst2} and the conditions
of this theorem we infer, that
$$ \int_{\psi } f(\zeta ) d Ln (\zeta -z) =
\int_{\psi _y} f(\zeta ) d Ln (\zeta -z)\mbox{  and}$$
$$ \int_{\gamma +z} f(\zeta ) d Ln (\zeta -z) =
\int_{\gamma _y +z} f(\zeta ) d Ln (\zeta -z)$$ for each $0<y\le 1$.
\par Since $\gamma _y+z$ is
a circle around $z$ its radius $y\rho >0$ can be chosen so small,
that $f(\zeta )=f(z)+ \alpha (\zeta ,z)$, where $\alpha $ is a
continuous function on $U^2$ such that the limit $\lim_{\zeta \to
z}\alpha (\zeta ,z)=0$ exists, then
$$\int_{\gamma _y+z}f(\zeta ) d Ln (\zeta -z) =$$
$$\int_{\gamma _y+z}f(z) d Ln (\zeta -z) +\delta (y\rho )=
2\pi f(z) M + \delta (y\rho ),$$ where

$$|\delta (y\rho )|\le |\int_{\gamma _y + z}\alpha (\zeta ,z)
d Ln (\zeta -z) | \le 2\pi \sup_{\zeta \in \gamma _y} |\alpha (\zeta
,z)| C_1 \exp (C_2 (y\rho )^m),$$

where $C_1$ and $C_2$ are positive constants, $m=2+2^r$ (see
Inequality $(2.7.4)$ \cite{luoyst,luoyst2}), hence there exists
$\lim_{y \to 0, y
>0}\delta (y\rho )=0$. Analogous estimates are for $\psi _y$ instead
of $\gamma _y+z$.
\par We have that $M = [\int_{\psi _y} d Ln (\zeta -z)]/(2\pi )$ and
also $M = [\int_{\gamma _y+z} d Ln (\zeta -z)]/(2\pi )$ for each $0<
y\le 1$ due to conditions of this theorem, since $Ln (y(\zeta -z)) =
\ln (y) + Ln (\zeta -z)$ for each $\zeta \ne z$ and $y$ is the
positive parameter independent from $z$, where $\ln $ is the
standard natural logarithmic function for positive numbers.
\par On the other hand, $Ln (1+z) = z - z^2/2 + z^3/3 - z^4/4 +...+
(-1)^{n+1}z^n/n +...$ for each Cayley-Dickson number of absolute
value $|z|<1$ less than one, since each $z$ can be written in the
form $z= Re (z) + Im (z)$, where $Re (z) =(z+{\tilde z})/2$, $Im (z)
= z - Re (z)$, $Im (z)$ is a purely imaginary number so that $(Im
(z))^2 = - | Im (z)|^2$. The latter series uniformly converges in a
ball of a given radius $0<R<1$ with the center at zero. Therefore,
the winding numbers of $\psi _y$ and $\gamma +z$ around $z$ are the
same for each $0<y\le 1$ and equal to $\int_{\gamma +z}dLn (\zeta
-z)/(2\pi M)$. Taking the limit while $y >0$ tends to zero yields
the conclusion of this theorem, since
$$\lim_{y\to 0} \int_{\gamma _y+z}f(\zeta )d Ln (\zeta -z) =
\lim_{y\to 0} \int_{\psi _y}f(\zeta ) d Ln (\zeta -z) .$$

\par If either $r=2, 3$, or
$f(z)\in \bf R$ for each $z$, then $((2\pi )f(z)M)M^*=2\pi f(z)$.

\par {\bf 4. Definitions.} This definition 1 of the residue spreads also on loops
$\psi $ not necessarily in a definite plane like $\psi $ in \S \S 2
and 3. But then the purely imaginary Cayley-Dickson number
substitutes on the mean value of $M$ computed with the help of the
line integral for $d Ln (z)$ along the loop $\psi $ divided on $2\pi
n$, that is $M=\int_{\psi }d Ln (z)/(2\pi n)$, where $n\in \bf N$ is
the winding number of the loop $\psi $. Due to \S \S 3 we have that
this operator $Res$ at a marked point $a=z_0$ is the same for
$\gamma +z_0$, $\gamma _y+z_0$, $\psi $ and $\psi _y$ loops for each
$0<y\le 1$, when the loops $\gamma _y+z_0$ encompass only one
singular point $z_0$ in the complex plane $(z_0+{\bf R}\oplus M{\bf
R})$ for $\gamma _y+z_0$, while $\psi _y$ encompasses $z_0$ relative
to $\gamma _y+z_0$ in the sense of \S 2.

\par For the
infinite point $a=\infty $ the definitions of the index and the
residue change so that take instead of a circle $\gamma $ or a loop
$\psi $ it with the opposite orientation $\gamma ^-$ or $\psi ^-$
respectively. Thus taking the mapping $z\mapsto 1/z$ reduce the pole
at the infinity into the pole at zero. In addition in the residue's
definition we take a sufficiently large radius $\rho
>0$ such that in the complex plane ${\bf R}\oplus M{\bf R}$ for a
marked purely imaginary Cayley-Dickson number the circle $\gamma ^-$
with the center at zero encompasses only one singular point $\infty
$ from $W$ in the non-commutative analog ${\hat {\cal A}}_r$ of the
Riemann sphere. We also consider circles $\gamma ^-_y$ with $y\ge 1$
and take the limit
$$ (i)\quad Res (\infty ,f).M := (2\pi )^{-1} \lim_{y\to \infty }\int_{\gamma
^-_y}f(z)dz$$  whenever it exists.

\par Mainly we shall consider residues at isolated singular
points, when another will not be specified.
\par If $f$ has an isolated singularity at $a\in {\hat {\cal A}}_r$,
then coefficients $b_k$ of its Laurent series (see \S 3
\cite{luoyst,luoyst2}) are independent from $\rho >0$. The common
series is called the $a$-Laurent series. If $a=\infty $, then
$g(z):=f(z^{-1})$ has a $0$-Laurent series $c_k$ such that
$c_{-k}=b_k$. Let $$\beta := \sup_{b_k\ne 0} \eta (k),$$ where $\eta
(k)=k_1+...+k_m$, $m=m(k)$ for $a=\infty $;
$$\beta =\inf_{b_k\ne 0} \eta (k)$$ for $a\ne \infty $. We say that
$f$ has a removable singularity, pole, essential singularity at
$\infty $ according as $\beta \le 0$, $0<\beta <\infty $, $\beta
=+\infty $. In the second case $\beta $ is called the order of the
pole at $\infty $. For a finite $a$ the corresponding cases are:
$\beta \ge 0$, $-\infty <\beta <0$, $\beta =-\infty $. If $f$ has a
pole at $a$, then $|\beta |$ is called the order of the pole at $a$.
\par A value of a function
$$\partial _f(a) := \inf \{ \eta (k): b_k\ne 0 \} $$
is called a divisor of $f$ at $a\ne \infty $,
$$\partial _f(a) := \inf \{ - \eta (k): b_k\ne 0 \} $$ for $a=\infty $,
where $b_k\ne 0$ means that $b_{k,1}\ne 0$,....,$b_{k,m(k)}\ne 0$.
Then $$\partial _{f+g}(a)\ge \min \{ \partial _f(a), \partial _g(a)
\} $$ for each $a\in dom (f)\cap dom (g)$ and $$\partial
_{fg}(a)=\partial _f(a)+\partial _g(a).$$  For a function $f$
meromorphic on an open subset $U$ in ${\hat {\cal A}}_r$ the
function $\partial _f(p)$ by the variable $p\in U$ is called the
divisor of $f$.

\par {\bf 5. Example and Remark.} Consider another example
of a multi-valued locally analytic function. Let $\xi $ be a
Cayley-Dickson number and consider the power function $z = \xi ^n$,
where $2\le n\in \bf N$, $n$ is a natural number. Each
Cayley-Dickson number has the polar form $z = \rho \exp (2\pi \theta
M)$, where $\theta \in \bf R$, $M$ is a purely imaginary
Cayley-Dickson number of the unit absolute value $|M|=1$, $\rho :=
|z|$. For definiteness we can consider $M$ of the form $M = M_1i_1
+...+ M_ki_k$ with $1\le k\le 2^r-1$, $M_1,...,M_k\in \bf R$,
$M_k>0$, that to exclude repeating caused by the identity $(-\theta
)(-M)= \theta M$. For each marked $M$ as above the function $\exp
(2\pi \theta M)$ is periodic so that $\exp (2\pi (\theta +n))M)=\exp
(2\pi \theta M)$ for each integer number $n$, $n\in \bf Z$. \par
Therefore, the inverse function $\xi = z^{1/n}$ on ${\cal
A}_r\setminus \{ 0 \} $ is multi-valued with $n$ branches. Each
branch is $\xi _j = \rho ^{1/n} \exp (2\pi (\theta +(j-1))M/n)$,
where $M=M(z)$, $b^{1/n}>0$ denotes the positive value of the $n$-th
root for each $b>0$, $j=1,...,n$. When a loop, for example, a circle
with the center at zero in the plane ${\bf R}\oplus M{\bf R}$ is
gone around $0$ on $2\pi $, then there is the transition from the
$j$-th branch into $j+1$ for each $1\le j<n$, while the $n$-th
branch transits into the first branch. That is, the point $z=0$ is
the branching point of the function $z^{1/n}$. Branches of this
function are indexed by $j=1,...,n$ and depend on purely imaginary
Cayley-Dickson numbers $M\in {\cal I}_r$. Since $z^{1/n} = \exp [(Ln
(z))/n]$, then $z^{1/n}$ is ${\cal A}_r$-holomorphic on ${\cal
A}_r\setminus \{ 0 \} $. This also follows from the inverse function
theorem: if $f$ is ${\cal A}_r$-holomorphic on an open subset $V$
and it has an inverse function $f^{-1}$, then its inverse function
$f^{-1}$ is also ${\cal A}_r$-holomorphic \cite{luoyst,luoyst2,
norfamlud}.
\par The non-commutative Riemannian surface for it is constructed
from $n$ copies of the Cayley-Dickson algebra ${\cal A}_r$ cut by
the subset $Q$ as in \S 3.7 \cite{luoyst,luoyst2} such that they are
embedded into ${\cal A}_r^2$. Each cut copy $Y_j$ of ${\cal A}_r$ is
slightly bend and each edge $\mbox{ }_2Q_j$ is glued with $\mbox{
}_1Q_{j+1}$ for each $j=1,...,n-1$ and $\mbox{ }_2Q_n$ with $\mbox{
}_1Q_1$ by straight rays with initially parallel directing vectors.
Then for such non-commutative Riemannian $2^r$-dimensional surface
${\cal R} = {\cal R}(z^{1/n})$ the function $z^{1/n}: ( {\cal
A}_r\setminus \{ 0 \} ) \to {\cal R}$ is already univalent.

\par Generally suppose that $f: (U\setminus W)\to {\cal A}_r$ is a
${\cal A}_r$-holomorphic function on an open set $(U\setminus W)$ in
the one-point compactification ${\hat {\cal A}}_r$ of the
Cayley-Dickson algebra ${\cal A}_r$, so that the compactification is
relative to the weak topology for $card (r)\ge \aleph _0$. We have
${\hat {\cal A}}_r\setminus {\cal A}_r = \{ \infty \} $. Here we
consider $W$ a closed connected subset of codimension at least $2$,
$codim (W)\ge 2$, in the Cayley-Dickson algebra ${\cal A}_r$ such
that $W\cap (z + {\widehat{{\bf R}\oplus M{\bf R}}})$ is a set
consisting of isolated points for each purely imaginary $M\in {\cal
I}_r$ and each $z\in W\cup \{ 0 \} $, $U$ is open in ${\cal A}_r$.
We also suppose that $f$ is multi-valued and each its branch is
defined on $(U\setminus W)$. \par If $z\in W$ and along each loop in
$(z + {\widehat {{\bf R}\oplus M{\bf R}}})$ encompassing $z$ and
neither encompassing nor containing another points from $W$ the
complete circuit around the loop on $2\pi $ leads to a transition
from the $j$-th branch of $f$ into a definite branch $k(j)$ for each
$j$, then we say that $z$ is a branching point of the function $f$.
\par A non-commutative surface ${\cal R} = {\cal R}(f)$
composed from copies of $Y_k$ properly embedded into a suitable
space ${\cal A}_r^s$, $card (s)\ge 2$, and with suitable gluing of
$Y_j$ with $Y_{k(j)}$ by the corresponding edges so that $f:
(U\setminus W)\to {\cal R}$ becomes univalent and ${\cal
A}_r$-holomorphic, then ${\cal R}$ is called the non-commutative
Riemannian surface of $f$.
\par If $f$ is such function, then for suitable $U$ satisfying
conditions of Theorem 3 with $z_0$ replaced on $W$, due to the
monodromy Theorem 2.41 \cite{norfamlud} and the homotopy Theorem
2.15 \cite{luoyst,luoyst2} the branching of $f$ depends only on
$M\in {\cal I}_r$ and is independent from $|\zeta -z|$ for $z\in W$.
This follows from the consideration of loops $\gamma $ in planes $z+
(\widehat{{\bf R}\oplus M{\bf R}})$ with $z\in W$ and $\gamma
\subset ((U_j\setminus W)\cap [z+ (\widehat{{\bf R}\oplus M{\bf
R}})])$ for each $j$.

\par Mention that generally conditions of Theorem 2.15 \cite{luoyst,luoyst2} can not be
replaced on local homotopies within each $U_j\setminus W$. This is
caused by several reasons. At first, consider two circles $\gamma
_1$ and $\gamma _2$ of radii $0<R_1<R_2<\infty $ with the center at
a point $z_0$ in planes $z_0 + ({\bf R}\oplus M{\bf R})$ and $z_0 +
({\bf R}\oplus N{\bf R})$, where $N = aM + bN_1$, $a, b\in {\bf
R}\setminus \{ 0 \} $, $Re (MN_1) =0$, $|M|=|N|=|N_1|=1$, $M$ and
$N$ and $N_1$ are purely imaginary marked Cayley-Dickson numbers.
Let $U$ and $z_0$ satisfy conditions of Theorem 3. Consider $\gamma
_1 (t) = z_0 + R_1\exp (2\pi tM)$ and $\gamma ^-_2(t) = z_0 + R_2
\exp (- 2\pi tN)$, $t\in [0,1]$. Divide them on arcs respective to
the partition $[(j-1)/n,j/n]$ with $4\le n\in \bf N$ and
$j=1,...,n$. Take the corresponding arcs of these circles and join
their respective ends $\gamma _1(j/n)$ with $\gamma _2(j/n)$ for
each $j$ by segments $w_j$ of straight lines. We get loops $\eta _j$
from such segments and arcs of $\gamma _1$ and $\gamma ^-_2$ circuit
in accordance with $\gamma _1$ and $\gamma ^-_2$. \par Now consider
two circles $\gamma _1$ and $\gamma _2$ embedded into $z_0 + ({\bf
R}\oplus M{\bf R}\oplus N_1{\bf R})$. Take a new system of
coordinates in the latter three dimensional manifold with the origin
at $z_0$. The axis $e_1$ take parallel to $i_0$, the axis $e_2$
parallel to $M$ and $e_3$ parallel to $N_1$. Therefore, $\gamma _2$
is in a half-space above $\gamma _1$ and in the other half-space
below $\gamma _1$. Then join $\gamma _1(t)$ with $\gamma _2(t)$ by
the segment of the straight line for each $t\in [0,1)$. So we get a
set $L$ which is homeomorphic to the M\"obius band, which is a
non-oriented two dimensional surface. Suppose that $L\subset
U\setminus W$.
\par There exists $k$ for which $w_k$ is gone twice in the same
direction, because $L$ is non-oriented. Then the sum of the
integrals by all $j$-th loops $\eta _j$ generally does not vanish
for an ${\cal A}_r$-holomorphic function $f$ in $(U\setminus \{ z_0
\} )$. Generally it has the rest of the type $2 c \int_{w_k}f(z)dz$,
where $c=1$ or $c= -1$.
\par The surface $Q$ cutting ${\cal A}_r$ has the unit codimension,
$codim (Q)=1$, consequently, the loops $\gamma $ and $\psi - z_0 $
from  \S 3 both intersect $Q$, $Q\cap \gamma ([0,1])\ne \emptyset $
and $Q\cap \psi ([0,1])\ne \emptyset $. So the integral by the loop
$\eta _j$ may be different from zero for some $j$, when $\gamma $
and $\psi $ are not in one plane ${\bf R}\oplus M{\bf R}$ and $z_0$
is a branching point of $f$.
\par If $f$, $U$, $z_0$, $\psi $ and $\gamma +z_0$ are as in Theorem 3
such that $\psi ([0,1])$ does not intersect $\gamma ([0,1])+ z_0$,
$\psi ([0,1])$ is not contained in the plane $z_0 + ({\bf R}\oplus
M{\bf R})$, then from $\psi $ and $\gamma $ construct a loop $\eta $
taking their arcs and joining their ends by rectifiable paths in $U$
so that $\eta \subset U\setminus W$. Then there exists a purely
imaginary Cayley-Dickson number $N\in {\cal I}_r$ such that it is
not contained in $M{\bf R}$ for which a projection of $\eta $ on the
plane $z_0 + ({\bf R}\oplus N{\bf R})$ encompasses the point $z_0$.
If $f$ has a singularity at $z_0$, which may be a branching point,
then generally $\int_{\gamma +z_0}f(z)dz$ is different from
$\int_{\psi }f(z)dz$, since a branching of $f$ generally depends on
all purely imaginary Cayley-Dickson numbers $N\in {\cal I}_r$,
$|N|=1$. This is easily seen on the examples of the functions $Ln
(\zeta -z_0)$ and $(\zeta - z_0)^{1/n}$ for $2\le n\in \bf N$.
Moreover, for the domain $U$ satisfying conditions of Theorem 3
generally conditions of the monodromy Theorem 2.41 in
\cite{norfamlud} are not satisfied.

\par In the particular case of $\gamma $ and $\psi $ contained
in the same plane $z_0 + ({\bf R}\oplus M{\bf R})$ the above
situation with $N_1\perp M$ is already excluded. On the other hand,
the branching of $f$ is independent from $|\zeta - z_0|$ (see
above). If there are not another singular points in $U$ besides
$z_0$ and $U$ satisfies conditions of Theorem 3, and when $\gamma
+z_0$ and $\psi $ are homotopic relative to $U\cap (z_0 + ({\bf
R}\oplus M{\bf R}))$, then due to the homotopy Theorem 2.15
\cite{luoyst,luoyst2} and since $M$ is fixed we can conclude as in
the complex case that $\int_{\gamma +z_0} f(z)dz = \int_{\psi }
f(z)dz$.
\par Therefore, the conditions of the homotopy Theorem 2.15
\cite{luoyst,luoyst2} are substantially stronger, than in the
complex case.

\par {\bf 6. Theorem.} {\it Let $U$ be an open region in the set
${\hat {\cal A}}_r$, $2\le r\le \infty $, with $n$ distinct marked
points $p_1,...,p_n$, and let $f$ be an ${\cal A}_r$-holomorphic
function on $U\setminus \{ p_1,...,p_n \} =:U_0$ and $\psi $ be a
rectifiable closed curve lying in $U_0$ such that $U_0$ satisfies
the conditions of Theorem 3 for each $z_0\in \{ p_1,...,p_n \} $.
Then
$$\int_{\psi }f(z)dz= 2 \pi \sum_{j=1}^n Res (p_j,f).
{\hat I}n (p_j,\psi )$$ and $Res (p_j,f).M$ is the $\bf
R$-homogeneous ${\cal I}_r$-additive (of the variable $M$ in ${\cal
I}_r$) ${\cal A}_r$-valued functional for each $j$.}
\par {\bf Proof.} In the considered case each singular point $p_j$ is
isolated, hence in the Definitions 1, 4 of the residue operator $Res
(p_j,f)$ the limit by $y$ can be omitted, since all integrals with
different $0<y\le 1$ are equal, when $\gamma _y$ are in the same
complex plane (see \S 3 above and \S 3.9.3 \cite{luoyst,luoyst2}).
\par For each $p_j$ consider the principal part $T_j$ of a Laurent
series for $f$ in a neighborhood of $p_j$, that is,
$$T_j(z)=\sum_{k, \eta (k)<0} \{ (b_k,(z-p_j)^k) \} _{q(m(k)+ \eta (k))},$$
where $\eta (k)=k_1+...+k_n$ for $k=(k_1,...,k_n)$ (see Theorem 3.21
\cite{luoyst,luoyst2}). Therefore, $$h(z):=f(z)-\sum_jT_j(z)$$ is a
function having an ${\cal A}_r$-holomorphic extension on $U$. In
view of Theorem 3 for an ${\cal A}_r$-holomorphic function $g$ in a
neighborhood $V$ of a point $p$ and a rectifiable closed curve
$\zeta $ we have
$$g(p) {\hat I}n (p,\zeta )=(2\pi )^{-1}(\int_{\zeta }g(z)(z-p)^{-1}dz)$$
(see \S 3.22 \cite{luoyst,luoyst2}). Let loops $\gamma $ and each
$\zeta _j$ be lying in planes each of which is parallel to the
complex plane ${\bf R}\oplus M{\bf R}$ with a marked purely
imaginary $M$ so that $M$ and $\psi $ satisfy conditions of Theorem
3.
\par We may consider small loops $\zeta _j$ around each $p_j$ with
${\hat I}n(p_j,\zeta _j)={\hat I}n (p_j,\gamma )$ for each
$j=1,...,n$. Then we infer that
$$\int_{\zeta _j}f(z)dz=\int_{\zeta _j}T_j(z)dz$$ for each $j$.
Representing $U_0$ as a finite union of open regions $U_j$ and
joining $\zeta _j$ with $\gamma $ by paths $\omega _j$ going in one
and the opposite direction as in Theorem 3 we get
$$\int_{\gamma }f(z)dz+\sum_j\int_{\zeta _j^-}f(z)dz=0,$$
consequently,
$$\int_{\gamma }f(z)dz=\sum_j\int_{\zeta _j}f(z)dz=
\sum_j 2\pi Res (p_j,f) {\hat I}n (p_j,\gamma ) ,$$ where ${\hat I}n
(p_j,\gamma )$ and $Res (p_j,f)$ are invariant relative to
homotopies satisfying conditions of Theorem 3 within a given complex
plane. \par Since the integral $\int_{\zeta _j} g(z)d Ln (z-p_j)$ is
$\bf R$-homogeneous and ${\cal I}_r$-additive relative to a
directing vector $M\in {\cal I}_r$ of a loop $\zeta _j$, then $Res
(p_j,f)M$ defined by Formulas 1$(i,ii)$, 4$(i)$ is $\bf
R$-homogeneous ${\cal I}_r$-additive of the argument $M$ in ${\cal
I}_r$.

\par {\bf 7. Theorem.} {\it Let $f$ be a meromorphic function
in $V$ (see \S 1), so that $f(\zeta )=g(\zeta )[1/v(\zeta )]$ and
$1/f(\zeta )= v(\zeta )[1/g(\zeta )]$ for each $\zeta \in V$, where
$g$ and $v$ are two ${\cal A}_r$-holomorphic functions on an open
set $V$, $g$ has not zeros in $V$, $v^{-1}(0)=W$, $W\subset V$,
$2\le r<\infty $, $codim (W)\ge 2$, then $Res (z,f)$ is an operator
continuously depending on $z\in W$.}
\par {\bf Proof.} The set $W$ is closed and connected in ${\hat
{\cal A}}_r$ and the intersection $W\cap (z + \widehat{{\bf R}\oplus
M{\bf R}})$ consists of isolated points for each $z\in W\cup \{ 0 \}
$ and purely imaginary $M\in {\cal I}_r$. Therefore, ${\hat {\cal
A}}_r\setminus W$ is open in ${\hat {\cal A}}_r$. Since $v$ is
continuous and $W=v^{-1}(0)$ and $g$ has not zeros in $V$, then each
point $z$ in $W$ is singular for $f$. Making a shift $z\mapsto
z+z_0$ we can consider that $0$ and $\infty $ do not belong
simultaneously to $W$. To treat a singularity at $\infty $ we can
use the transformation $z\mapsto 1/z$ under which we can consider
instead a singularity at zero.  So we can suppose without loss of
generality that each $z\in W$ is a finite singular point.
\par We can take a neighborhood $U$ of $z$ satisfying conditions of
Theorem 3 and choose the family $\gamma _y+z$ of loops in $U$
encompassing only one singular point $z$ of $f$ in the complex plane
$z+{\bf R}\oplus M{\bf R}$ for each $M\in {\cal I}_r$ with $|M|=1$
and each $0<y\le 1$ and every $z\in W\cap U$. The integral
$$\int_{\gamma _y +z} f(\zeta )d\zeta =\int_{\gamma _y} f(\zeta
+z)d\zeta $$ is the continuous functional relative to $f$ and $z$,
since $f$ is holomorphic in $U\setminus W$, $\gamma
_y([0,1])+z\subset U\setminus W$ for each $z\in W\cap U$ and every
$0<y\le 1$. \par If $z_1\in W$ is another either singular point
different from $z$ or any other point $z_1\in U\setminus \{ z \} $,
then $|z-z_1|>0$. Choose $0<y<|z-z_1|/\rho $, where $0<\rho <\infty
$ is a radius of the circle $\gamma $, then $\gamma _y+z$ does not
encompass $z_1$ in the entire $U$ as well. Since $2 \pi ~ Res
(z,f).M = \lim_{y\to 0} \int_{\gamma _y +z} f(z)dz$, when $M$ is
fixed, hence $Rez (z,f).M$ may depend on $z$ and the residue $Res
(z,f).M$ is independent from each $z_1$ not equal to $z$, $z_1\ne
z$.

\par The functions $g$ and $v$ being ${\cal A}_r$-holomorphic are
locally $\zeta $-analytic in $V$ and for each $z\in W$ there exists
a ball with center at $z$ of some radius $0<\delta <\infty $ such
that $g$ and $v$ have non-commutative power series expansions in it.
Since $g(\zeta )\ne 0$ for each $\zeta \in V$, then $g(z)\ne 0$ and
$v(z)=0$ for each $z\in W$, consequently, $f(z)=\infty $ and
$1/f(z)=0$, and inevitably $1/f(\zeta )=v(\zeta )[1/g(\zeta )]$ is
${\cal A}_r$ holomorphic in $V$ and hence in a neighborhood $V_z$ of
$z\in W$, $V_z\subset V$. Thus the function $1/f(\zeta )$ has the
power series expansion in a ball with the center at $z\in W$ of some
radius $0<\epsilon <\infty $. Let its expansion coefficients be
$A_m=(a_{m,m_1},...,a_{m,m_k})$, $a_{m,j}\in {\cal A}_r$ for each
$m, j$, $m=(m_1,...,m_k)$, $m_j\ge 0$ for each $j=1,...,k$, $k\in
\bf N$ so that $$(1)\quad 1/f(\zeta ) = \sum_m \{ A_m, (\zeta -z)^m
\} _{q(m)},$$ where $ \{ A_m, z^m \} _{q(m)} := \{
a_{m,m_1}z^{m_1}...a_{m,m_k}z^{m_k} \} _{q(m)}$, $q(m)$ is the
vector indicating on an order of the multiplication. Let us seek the
function $f$ in the form:
$$(2)\quad f(\zeta ) = \sum_p \{ B_p, (\zeta -z)^p \}
_{q(p)},$$ where $p=(p_1,...,p_k)$, $k=k(p)\in \bf N$, $p_j\in \bf
Z$ for each $j$, for each $p$ either all $p_j\ge 0$ or all $p_j\le
0$ simultaneously. Put $\eta (m) := m_1+...+m_k$, where $k=k(m)$,
denote $$Q_n(\xi ) := \sum_{\eta (p)=n} \{ B_p, \xi ^p \}
_{q(p)}\mbox{ and}$$  $$P_n(\xi ) := \sum_{\eta (m)=n} \{ A_m, \xi
^m \} _{q(m)}$$ homogeneous terms so that $Q_n(t\xi )=t^nQ_n(\xi )$
and $P_n(t\xi )=t^nP_n(\xi )$ for each $0\ne t\in \bf R$ and every
$\xi \in {\cal A}_r\setminus \{ 0 \} $. When $n>0$, then we have
$P_n(0)=0$ and $Q_n(0)=0$. On the other hand, for a negative number
$n$ the term $Q_n(1/\xi )$ is defined for each finite Cayley-Dickson
number $\xi $. \par For the function $1/f$ we have $\eta (m)=n\ge
1$, since $1/f$ is holomorphic and $1/f(z)=0$. If each $P_n(\xi )$
would be an identically zero polynomial, then $1/f$ would be
identically zero, that is not the case. Thus there exists the
maximal $\alpha
>0$ for which $P_{\alpha }(\xi )$ is the nontrivial polynomial,
while $P_n(\xi )$ is identically zero for each $n<\alpha $. That is
$\alpha $ is the order of zero $z$ of the function $1/f$.
\par We have the equation $$(3)\quad f(\zeta )[1/f(\zeta )]=1$$ identically
in the set $B({\cal A}_r,z,\epsilon )\setminus W$ and by the
continuity this equation extends on the entire ball. Since
$$(4)\quad Q_{n_1}(t\xi )P_{n_2}(t\xi ) = t^{n_1+n_2} Q_{n_1}(\xi )P_{n_2}(\xi
)$$ for each $0\ne t\in \bf R$ and every $\xi \in {\cal
A}_r\setminus \{ 0 \} $, then we infer the inequality $n_1\ge
-\alpha $. In view of $(4)$ Equation $(3)$ gives the system of
equations
$$(5)\quad \sum_{n_1 + n_2=l} Q_{n_1}(\xi )P_{n_2}(\xi )=\delta
_{l,0}$$ for each $0<|\xi |<\epsilon $, where $0\le l\in \bf Z$,
$\delta _{i,j} =0$ for each $i\ne j\in \bf Z$, $\delta _{j,j}=1$ for
each $j\in \bf Z$, consequently, $\alpha \le n_2\le l+\alpha $ in
each $l$-th equation $(5)$.
\par The Cayley-Dickson algebra ${\cal A}_r$ has the finite
dimension $2^r$ over the real field $\bf R$, hence for each $n$ the
number of different $\bf R$-linearly independent terms $ \{ B_p,\xi
^p \} _{q(p)}$ with $\eta (p)=n$ is finite, as well as a number of
different $\bf R$-linearly independent terms $ \{ A_m, \xi ^m \}
_{q(m)}$ is finite for $\eta (m)=n$. Thus each term $Q_n$ and $P_n$
is a finite sum and a number of expansion coefficients in them is
finite.
\par Now using the multiplication table in the Cayley-Dickson algebra
${\cal A}_r$, where $2\le r<\infty $, and exploiting the
decomposition $w=\mbox{ }_0w i_0+...+\mbox{ }_{2^r-1}w i_{2^r-1}$
for each $w\in {\cal A}_r$, where $\mbox{ }_0w,...,\mbox{
}_{2^r-1}w\in \bf R$, it is possible by induction on $l=0,1,...$
resolve the system $(5)$ relative to $B_p$ through $A_m$. The
procedure is the following. For $l=0$ we get the power relative to
$\mbox{ }_0\xi ,...,\mbox{ }_{2^r-1}\xi $ Equation $(5)$ expressing
vectors $B_p$ with $\eta (p)= -\alpha $ through vectors $A_m$ with
$\eta (m)=\alpha $. For $l=1$ Equation $(5)$ together with the
previous equation expresses vectors $B_p$ with $\eta (p) = - \alpha
+1 $ through vectors $A_m$ with $\eta (m)\in \{ \alpha , \alpha +1
\} $. By induction equations $(5)$ with $l=0,...,\beta $ express
vectors $B_p$ with $\eta (p)=\beta $ through vectors $A_m$ with
$\eta (m)\in \{ \alpha ,...,\alpha +\beta \} $. Thus each $B_p$ is a
continuous function of a finite number of $A_m=A_m(z)$. At the same
time each vector $A_m(\xi )$ is a continuous function of $\xi $ in
$B({\cal A}_r,z,\epsilon /2)$, since the function $1/f(\zeta -z)$ is
${\cal A}_r$-holomorphic by $\zeta , z\in {\cal A}_r$ with $\zeta -z
\in V$, where $V$ is open in the Cayley-Dickson algebra ${\cal A}_r$
by the conditions of this Theorem. Consequently, each $B_p(\xi )$ is
a continuous functions of $\xi $ in the ball $B({\cal
A}_r,z,\epsilon /2)$.
\par Evidently we have
$$\int_{\gamma _y}\{ B_p,(\zeta -z)^p \} _{q(p)} d\zeta =0$$
for each $\eta (p)\ne -1$ and each $0<y\le 1$, so that $z$ here is
equal to $a$ in the Definition of the residue in \S \S 1, 4.
Therefore, the limit
$$2\pi ~ Res (z,f).M :=
\lim_{y\to 0} \int_{\gamma _y}f(\zeta ) d\zeta $$
$$ = \lim_{y\to 0} \sum_{p, \eta (p)= -1}
\int_{\gamma _y} \{ B_p,(\zeta -z)^p \} _{q(p)} d\zeta $$
$$ = \sum_{p, \eta (p)= -1}
\int_{\gamma } \{ B_p,(\zeta -z)^p \} _{q(p)} d\zeta $$ exists,
since $\int_{\gamma _y}$ is the additive functional by integrands
and it is continuous on the space of bounded continuous functions
(see Theorem 2.7 \cite{luoyst,luoyst2}) and due to the homotopy of
$\gamma _y$ with $\gamma $ in the same complex plane
$z+\widehat{{\bf R}\oplus M{\bf R}}$ relative to $(V\setminus W)\cap
(z+\widehat{{\bf R}\oplus M{\bf R}})$ (see \S 3 above and 3.9.3
\cite{luoyst,luoyst2}). \par This implies that $Res (z,f)$ is
completely defined by the power series expansion of $f(z)$, namely
by vector coefficients $B_p$ with $\eta (p) = -1$ only.  The family
of coefficients with $\eta (p)= -1$ is finite for a finite $2\le
r<\infty $. Therefore, we get that $Res (z,f).M$ is the continuous
operator-valued function by the variable $z\in W$.

\par {\bf 8. Corollary.} {\it Let $U$ be an open region
in the set ${\hat {\cal A}}_r$, $2\le r\le \infty $, with $n$
distinct points $p_1,...,p_n$, let also $f$ be an ${\cal
A}_r$-holomorphic function on $U\setminus \{ p_1,...,p_n \} =:U_0$,
$p_n=\infty $, and $U_0$ satisfies conditions of Theorem 3 with at
least one $\psi $, $\gamma $ and each $z_0\in \{ p_1,...,p_n \} $.
Then $$\sum_{p_j\in U} Res (p_j,f)M =0.$$}
\par {\bf Proof.} If $\gamma $ is a closed curve encompassing
points $p_1$,...,$p_{n-1}$, then $\gamma ^-(t):=\gamma (1-t)$, where
$t\in [0,1]$, encompasses $p_n=\infty $ with positive going by
$\gamma ^-$ relative to $p_n$. Since $$\int_{\gamma
}f(z)dz+\int_{\gamma ^-}f(z)dz=0,$$  we get from Theorem 6, that
$$\sum_{p_j\in U} Res (p_j,f)M=0$$ for each $M\in {\cal I}_r$, hence
$$\sum_{p_j\in U} Res (p_j,f)M =0$$ is the zero $\bf R$-homogeneous
${\cal I}_r$-additive ${\cal A}_r$-valued functional on ${\cal
I}_r$.
\par {\bf 9. Remark.} If a rectifiable loop $\psi $ is as in \S 1,
then there exists a sequence $\psi _p$ of rectifiable loops composed
of arcs of circles centered at $z_0$ and segments of straight lines
such that straight lines contain $z_0$, so that $\psi _p$ converges
uniformly to $\psi $ with $p$ tending to the infinity. Therefore, if
$f$ is the ${\cal A}_r$ differentiable function in an open set $U$
so that $\psi \subset U$ and $\psi _p\subset U$ for each $p$, then
\par $\lim_{p\to \infty } \int_{\psi _p} f(z)dz = \int_{\psi }
f(z)dz$.
\par There are useful Moufang identities in the octonion
algebra: \par $(M1)$  $(xyx)z = x(y(xz))$,
\par $(M2)$ $z(xyx) = ((zx)y)x$,
\par $(M3)$ $(xy)(zx) = x(yz)x$ for each $x, y, z\in \bf O$ \\
(see page 120 in \cite{harvey}). For calculations it is also worth
that $a{\tilde a} = {\tilde a} a = |a|^2$ and the real valued scalar
product has the symmetry properties $Re (a{\tilde b}) = Re ({\tilde
b}a) = Re ({\tilde a} b) = Re (b{\tilde a })$, where $a^* = {\tilde
a}$ denotes the conjugated Cayley-Dickson number.
\par Mention also that $(de^z/dz).h = \sum_{n=1}^{\infty }
\sum_{j=0}^{n-1} z^jhz^{n-j-1}/n!$ for all $z, h\in {\cal A}_3={\bf
O}$, hence generally $(de^z/dz).h$ is not simply $e^zh$ besides the
case of $Im (h)\in Im (z){\bf R}$, where $Im (z) = z - Re (z)$ and
the order of the multiplication in each additive of the series
corresponds to the order of multiplication from right to left (right
order of brackets). Due to the Moufang identities an order of
multiplications either left or right or $((z^j)h)(z^{n-j-1})$ or
$(z^j)(h(z^{n-j-1}))$ gives the same result.  \par The inverse
operator $(de^z/dz)^{-1}$ to $de^z/dz$ provides $d Ln (y)/dy$ for
$y=e^z$. The operator $de^z/dz$ is ${\cal A}_r$ additive and $\bf R$
homogeneous, hence its inverse operator is such also.  Particularly,
if $Im (h)\in Im (y)\bf R$, where $y\ne 0$, then $(dLn (y)/dy).h
=y^{-1}h$. Therefore, $\int_{\gamma _k} dLn (y) = \int_{\gamma _k}
y^{-1}dy$ for each $\gamma _k$ either arc of a circle with the
center at zero or a segment of a straight line such that the
straight line contains $0$, consequently, $\int_{\psi _p-z_0} dLn
(y) = \int_{\psi _p-z_0} y^{-1}dy$ for each $p$ (see also \S 12).
\par Thus due to Definition 1 we deduce that Theorems 3 and
7, Corollary 8 spread on more general rectifiable loops $\psi $, but
then $M$ in formulas there substitutes on the mean value of $M$
computed as \par $(1)$ $M = \lim_{p\to \infty } \int_{\psi _p-z_0} d
Ln (z)/(2\pi j) = \int_{\psi -z_0}d Ln (z)/(2\pi j)$, where $j\in
\bf N$ is the winding number of the loop $\psi $ (see also
Definitions 1, 2 and Theorem 3). This formula for $M$ is essential,
since generally the logarithmic function over the Cayley-Dickson
algebra has more complicated non-commutative Riemannian surface,
than in the commutative complex case. If $a, b, c$ and $e$ are
constants in ${\cal A}_r$ so that $ \{ ab M ce \} _{q_1(5)} = \{
abce {\tilde M} \} _{q_2(5)}$, then
\par $(2)$ $\int_{\psi }d \{ ab Ln (z) ce \} _{q_1(5)} /(2\pi n) =
\{ abce {\tilde M} \} _{q_2(5)}$, \\ where vectors $q_1(5)$ and
$q_2(5)$ prescribe an order of the multiplication in brackets,
${\tilde M} = - M$ for purely imaginary Cayley-Dickson number $M$.
On the other hand, if $f$ is a (super) differentiable function in a
neighborhood of $\gamma $, then $\int_{\gamma } \{ abf(z)ce
\}_{q_1(5)} dz = \{ ab (\int_{\gamma } f(z)dz) ce \} _{q_1(5)}$ due
to the definition of the line integral \cite{luoyst,luoyst2}. If $N$
and $M$ are two orthogonal purely imaginary octonions, $Re (M{\tilde
N})=0$, also $Im (z) \in M{\bf R}$, then $N(zN^*) = (Nz)N^* =
{\tilde z}$.
\par Therefore, if for constants $a_1,...,b_{k+1}\in {\cal A}_r$ the identity
\par $(3)$ $ \{ a_1z^{n_1}...a_k z^{n_k}a_{k+1} \} _{q_1(2k+1)} = \{
b_1...b_{k+1} (z^s {\tilde z}^m) \} _{q_2(k+2)}$ is satisfied in an
open neighborhood of a rectifiable loop $\psi $ for $z_0=0$ and with
the winding number $j$, where $k\in \bf N$, $n_1,...,n_k, s, m \in
\bf Z$, then
\par $(4)$ $(2\pi j)^{-1} \int_{\psi } \{ a_1z^{n_1}...a_k z^{n_k}a_{k+1} \}
_{q_1(2k+1)} dz = \{ b_1...b_{k+1} M \} _{q_2(k+2)}$  for $m-s=1$
with $n= -1$;
\par $(5)$ $(2\pi j)^{-1} \int_{\psi } \{
a_1z^{n_1}...a_k z^{n_k}a_{k+1} \} _{q_1(2k+1)} dz = \{
b_1...b_{k+1} {\tilde M} \} _{q_2(k+2)}$  for $s-m=1$ with $n= -1$;
\par $(6)$ $(2\pi j)^{-1} \int_{\psi } \{
a_1z^{n_1}...a_k z^{n_k}a_{k+1} \} _{q_1(2k+1)} dz = 0$ for either
$|m-s|\ne 1$ or $n\ne -1$, where $n=n_1+...+n_k$. Since $s+m=n$,
then Case $(4)$ implies $m=0$ and $s= -1$, Case $(5)$ implies $s=0$
and $m= -1$.
\par {\bf 10. Definitions.} Let $f$ be an ${\cal A}_r$-holomorphic
function, $2\le r\le \infty $, on a neighborhood $V$ of a point
$z\in {\cal A}_r$. Then the infimum: $$\eta (z;f):=\inf \{ k: k\in
{\bf N}, f^{(k)}(z)\ne 0 \} $$ is called a multiplicity of $f$ at
$z$. Let $f$ be an ${\cal A}_r$-holomorphic function on an open
subset $U$ in the set ${\hat {\cal A}}_r$, $2\le r\le \infty $.
Suppose $w\in {\hat {\cal A}}_r$, then the valence $\nu _f(w)$ of
the function $f$ at $w$ is by the definition $$\nu _f(w):=\infty ,
\mbox{ when the set } \{ z: f(z)=w \} \mbox{ is infinite,}$$ and
otherwise
$$\nu _f(w):=\sum_{z, f(z)=w}\eta (z;f).$$

\par {\bf 11. Theorem.} {\it Let $\gamma $ and $\psi $ be two
rectifiable paths in ${\cal A}_r$ contained in open sets $U_{\gamma
}$ and $U_{\psi }$ respectively and let a diffeomorphism $\xi :
U_{\psi }\to U_{\gamma }$ exists with (super) differentiable $\xi $
and $\xi ^{-1}$ such that $\xi (\psi ([0,1]))=\gamma ([0,1])$. If
$f$ is a continuous ${\cal A}_r$ valued function on $U_{\gamma }$,
then
\par $(1)$ $\int_{\gamma }f(z)dz = \int_{\psi } f(\xi (y)).(\xi '(y).dy)$.}
\par {\bf Proof.} For $f$ in a  canonical closed bounded neighborhood $V$
of $\gamma $ so that $V\subset U$ take a sequence of (super)
differentiable functions $f_n$ converging uniformly on $V$ to $f$
with the corresponding phrases $\eta _n$ fixing a $z$-representation
of $f_n$ and of $f$ (see \S 2 in \cite{luoyst,luoyst2}). Then the
limit $\lim_{n\to \infty } \int_{\gamma } f_n(z)dz = \int_{\gamma }
f(z)dz$ exists (see Theorem 2.7 \cite{luoyst,luoyst2}). Therefore,
it is sufficient to prove this theorem for $f_n$ or in the case,
when $f$ is (super) differentiable. \par In accordance with
Proposition 2.6 \cite{luoyst,luoyst2} take a function $g(z)$ and its
phrase $\nu $ given by either left or right algorithm so that
$g'(z).1=f(z)$ and $\nu '(z).1 = \eta (z)$ for each $z\in V$, where
$g'(z) = dg(z)/dz$, $~ \eta (z)$ is the phrase of $f$. This
specifies the branch of the non commutative line integral
\par $(2)$ $\int_{\gamma } f(z)dz = \lim_{\delta (P)\to 0} \sum_{j=1}^m
(dg(z)/dz)|_{z=z_j}.\Delta _jz$, \\
where $P$ denotes a partition of $[0,1]$ with points $t_0=0<
t_1<...<t_m=1$ and $\tau _j\in [t_{j-1},t_j]$ for each $j=1,...,m$,
$\Delta _jz := \gamma (t_j) - \gamma (t_{j-1})$, $z_j := \gamma
(\tau _j)$, $\delta (P) := \max_{j=1,...,m}  (t_j-t_{j-1})$.
\par Since by the chain rule $(dg(\xi (y))/dy).h = (dg(z)/dz)|_{z=\xi
(y)}. ((d\xi (y)/dy).h)$ for each $y\in U_{\psi }$ and $h\in {\cal
A}_r$ (see also \S 2 \cite{luoyst,luoyst2}), then
\par $(dg(z)/dz)|_{z=\xi (y_j)}.((d\xi (y)/dy).\Delta _jy) =
(dg(z)/dz)|_{z=z_j}.\Delta _jz + o(\Delta _jz)$, \\ where $y_j =
\psi (\tau _j)$, $\Delta _jy = \psi (t_j) - \psi (t_{j-1})$,
$z_j=\xi (y_j)$, $\Delta _jz = \xi (\psi (t_j)) - \xi (\psi
(t_{j-1}))$. Therefore, \\ $(3)$ $\lim_{\delta (P)\to 0}
\sum_{j=1}^m (dg(z)/dz)|_{z=z_j}.\Delta _jz = \lim_{\delta
(P)\to 0} \sum_{j=1}^m (dg(z)/dz)|_{z=\xi (y_j)}.(\xi '(y_j).\Delta _jy)$, \\
since paths $\gamma $ and $\psi $ are rectifiable. Thus Formula
$(1)$ follows from $(2,3)$.
\par {\bf 12. Remark.} Mention that $\int_{\gamma } f(z)dz$ is
independent from parametrization of $\gamma $ in the following
sense. Two paths $\gamma _1: [a_1,b_1]\to {\cal A}_r$ and $\gamma
_2: [a_2,b_2] \to {\cal A}_r$ are called equivalent, $\gamma _1\sim
\gamma _2$, if there exists a continuous monotonously increasing
function $\phi : [a_1,b_1] \to [a_2,b_2]$ so that $\gamma _1(t) =
\gamma _2(\phi (t))$ for each $t\in [a_1,b_1]$, where $a_1<b_1$,
$a_2<b_2$. This relation is reflexive ($\gamma \sim \gamma $),
symmetric (if $\gamma _1\sim \gamma _2$, then $\gamma _2\sim \gamma
_1$), transitive (if $\gamma _1\sim \gamma _2$ and $\gamma _2\sim
\gamma _3$, then $\gamma _1\sim \gamma _3$). A class of equivalent
paths is called a curve.
\par If $\gamma _1\sim \gamma _2$ are two rectifiable paths, then
the mapping $\phi $ gives the bijective correspondence between
partitions of $\gamma _1$ and $\gamma _2$ and between the
corresponding integral sums. From the definition of the non
commutative line integral it follows, that if $\gamma _1\sim \gamma
_2$, then $\int_{\gamma _1} f(z) dz = \int_{\gamma _2} f(z)dz$ for a
continuous function $f$. Thus the non commutative line integral
depends on curves.

\par {\bf 13. Theorem.} {\it  If $\gamma $ is a rectifiable path $\gamma : [a,b]\to {\cal A}_r$
in the Cayley-Dickson algebra ${\cal A}_r$ and $f: U\to {\cal A}_r$
is a continuous function, where $U$ is an open subset in ${\cal
A}_r$ and $\gamma \subset U$, $-\infty < a<b<\infty $, then the non
commutative line integral reduces to the non commutative
Lebesgue-Stieltjes integral:
\par $(1)$ $\int_{\gamma } f(z) dz = (L-S) \int_a^b {\hat f}(\gamma
(t)).d\gamma (t)$.
\par Moreover, if $\gamma (t)$ is absolutely continuous, then the
non commutative line integral reduces to the non commutative
Lebesgue integral: \par $(2)$ $\int_{\gamma } f(z) dz = (L) \int_a^b
{\hat f}(\gamma (t)).\gamma '(t)dt$.}
\par {\bf Proof.} Recall that a function $f: [a,b]\to {\cal A}_r$ is with finite variation,
if there exists $C>0$ so that \par $V_a^b := \sup \sum_{j=1}^m
|f(t_j) - f(t_{j-1}| \le C$, \\
where the supremum  is taken by all finite partitions of the segment
$[a,b]$, $m\in \bf N$, $V_a^b(f)$ is called the variation of $f$ on
$[a,b]$. This implies that each $f_j$ is also with the finite
variation, where $f(t) =\mbox{ }_0f(t)i_0+...+\mbox{
}_{2^r-1}f(t)i_{2^r-1}$, $\mbox{ }_jf(t)\in \bf R$ for each
$j=0,...,2^r-1$ and every $t\in [a,b]$.
\par A path $\gamma : [a,b]\to {\cal A}_r$ is called rectifiable,
if it is continuous and of finite variation. \par A function $g:
[a,b]\to {\cal A}_r$ is called absolutely continuous, if for each
$\delta >0$ there exists $\epsilon >0$ so that for each system of
pairwise disjoint intervals $(a_j,b_j)$ in $[a,b]$, $j=1,...,m$,
$m\in \bf
N$, the inequality \par  $\sum_{j=1}^m |f(b_j)-f(a_j)|<\epsilon $ \\
is satisfied. A continuous function is called singular if its
derivative is almost everywhere relative to the Lebesque measure
equal to zero. \par In view of \S VI.4 \cite{kolmfom} each function
$f$ of finite variation can be presented as $f=H+\psi +\chi $, where
$H$ is the jump function, $\psi $ is the absolutely continuous
function and $\chi $ is the singular function. For a rectifiable
path $\gamma $ we get $H=0$. By Theorem VI.4.2 \cite{kolmfom} the
function $F(x) = (L) \int_a^xp(t)dt$ of a Lebesgue integrable
function $p(t)$ on the segment $[a,b]$ is absolutely continuous,
where $(L) \int_a^b p(t)dt$ is the Lebesgue integral. In view of the
Lebesgue Theorem IV.4.3 \cite{kolmfom} the derivative $p'=P$ of an
absolutely continuous function $p$ on the segment $[a,b]$ is
Lebesgue integrable on it and
\par $(L)\int_a^x p(t)dt = P(x) - P(a)$ \\ for each $a\le x\le b$.
Thus $\gamma $ is almost everywhere on $[a,b]$ differentiable.
\par The operator $\hat f$ is defined for each (super) differentiable function
with the help of either the left or the right algorithm, where
${\hat f}.h := g'(z).h$ for each $z\in V$, $h\in {\cal A}_r$.
Moreover, $\hat f$ is (super) differentiable, if $f$ is such. Then
we use the continuous extension of the continuous functional
$\int_{\gamma }$ on the space of continuous functions on $V$
specifying the branch of the integral and get the operator $\hat f$
for a chosen sequence of (super) differentiable functions $f_n$ and
their phrases $\eta _n$ converging uniformly on $V$ to $f$ (see
\cite{luoyst,luoyst2} and \S 11 above).
\par Each ${\cal A}_r$ additive $\bf R$ homogeneous operator $A$ on
${\cal A}_r$ can be written in the form: \par $A.h =
\sum_{k=0}^{2^r-1}A.[\mbox{ }_khi_k] = \sum_{j,k=0}^{2^r-1} \mbox{ }_{j,k}A.[\mbox{ }_kh]i_j$, \\
where each $\mbox{ }_{j,k}A$ is a linear functional on $\bf R$ for
all $j, k$. The analogous decomposition is for an operator ${\cal
A}_r$ additive $\bf R$ homogeneous valued function $A(t).h$, $t\in
[a,b]$. If each $\mbox{ }_{j,k}A(t)$ function is Lebesgue-Stieltjes
integrable on $[a,b]$ relative to each function $\mbox{ }_kP$ of
bounded variation, then there is defined the non commutative
Lebesgue-Stieltjes integral
\par $(L-S) \int_a^b A(t).dP = \sum_{j,k} [\int_a^b \mbox{ }_{j,k} A(t).d
\mbox{}_kP(t)] i_j$. On the the other hand, if $v(t) = \mbox{ }_0v
i_0 +...+ \mbox{ }_{2^r-1}v i_{2^r-1}$ is a function so that each
$\mbox{ }_jv(t)$ is Lebesgue integrable on $[a,b]$, then there is
defined the non commutative Lebesgue integral
\par $(L)\int_a^b v(t)dt = \sum_{j=0}^{2^r-1} \int_a^b \mbox{ }_jv(t)dt$.
\par Consequently, \par
$\int_{\gamma } f(z)dz = (L-S) \int_a^b {\hat f}(\gamma (t)).d\gamma
(t)$, \\ where ${\hat f}$ is an ${\cal A}_r$ additive $\bf R$
homogeneous operator corresponding to $f$, $(L-S) \int_a^b {\hat
f}(\gamma (t)).d\gamma (t)$ denotes the non commutative
Lebesgue-Stieltjes integral. If $\gamma (t)$ is absolutely
continuous, then the non commutative Lebesgue-Stieltjes integral
reduces to the non commutative Lebesgue integral
\par $(L-S) \int_a^b {\hat f}(\gamma (t)).d\gamma (t)= (L) \int_a^b
{\hat f}(\gamma (t)).\gamma '(t)dt$, since \\
$(L-S) \int_a^b \mbox{ }_{j,k}A(t)d\mbox{}_k\gamma (t) = (L)
\int_a^b \mbox{ }_{j,k}A(t)\mbox{ }_k\gamma '(t)dt$ \\
for each $j, k$ (see also \S VI.6.2 \cite{kolmfom}).
\par {\bf 14. Particular cases of integrals and
residues.}
\par For convenience we can
choose the parametrization of the curve so that its path $\gamma :
[\alpha ,\beta ]\to {\cal A}_r$ has $[\alpha ,\beta ]=[0,1]$ if
another is not specified, where $\gamma (0)=a$, $\gamma (1)=b$ are
Cayley-Dickson numbers.
\par While calculation of line integrals of ${\cal A}_r$
differentiable functions $f: U\to {\cal A}_r$, where $U$ is a domain
in ${\cal A}_r$ satisfying conditions of the homotopy Theorem 2.15
\cite{luoyst,luoyst2}, for a rectifiable path $\gamma $ in $U$
contained in a complex plane $\widehat{{\bf R}\oplus M{\bf R}}$ with
a marked purely imaginary Cayley-Dickson number $M$, $|M|=1$, it is
possible to make simplifications in algorithms. For this choose
purely imaginary Cayley-Dickson numbers $N_1$,...,$N_{2^r-1}$ being
generators in ${\cal A}_r$, $2\le r\in \bf N$, such that $N_1=M$,
$|N_1|=1$,...,$|N_{2^r-1}|=1$, $N_j\perp N_k$ for each $1\le j\ne
k\le 2^r-1$, that is $Re (N_jN_k)=0$, and hence certainly $N_jN_k= -
N_kN_j$ for all $1\le j\ne k\le 2^r-1$,...,$N_0=1$, $N_0N_j=N_jN_0$
for all $j$. This can be done standardly by induction using doubling
generators $N_2$, $N_4$,...,$N_{2^{r-1}}$, so that $N_3=N_1N_2$,
...,$N_{2^{r-1}+p} = N_pN_{2^{r-1}}$ for each $p=1,...,2^{r-1}-1$.
With the doubling procedure of ${\cal A}_{s+1}$ from ${\cal A}_s$,
$1\le s$, the multiplication rule is given by the formula:
\par $(1)$ $(u+vl)(w+xl) = (uw - {\tilde x} v) + (xu + v {\tilde w} )l$, \\
for each $u, v, w, x \in {\cal A}_s$, where $l=N_{2^s}$, $u+vl$ and
$w+xl \in {\cal A}_{s+1}$ (see \cite{baez,kansol,kurosh}).
\par With this new basis of generators write $f$ in the form
\par $(2)$ $f(z) = \mbox{ }_0g(z)N_0 + \mbox{ }_1g(z)N_2 + ... + \mbox{
}_{2^{r-1}-1}g(z)N_{2^r-2}$, \\
where $\mbox{ }_pg(z)\in {\bf C}_M$ for each $z\in U$ and all
$p=0,...,2^{r-1}-1$, while ${\bf C}_M := {\bf R}\oplus M{\bf R}$
denotes the complex plane embedded into ${\cal A}_r$. If $z$ is an
arbitrary Cayley-Dickson number $z\in {\cal A}_r$, then it can be
written as \par $(3)$ $z = x + y N_z$, \\ where $x=x_z\in {\bf
C}_M$, $y=y_z\in \bf R$, $N_z$ is a purely imaginary Cayley-Dickson
number may be dependent on $z$, $|N_z|=1$, $N_z\perp M$, which
follows also from $(2)$. Since $N_z\perp M$ and they are purely
imaginary, then $N_zM = -MN_z$.  Thus $z^k = x_k + y_k N_{z^k}$ for
any integer $k\in \bf Z$ with $x_k\in {\bf C}_M$ and $y_k\in \bf R$,
with purely imaginary $N_{z^k}$, $N_{z^k}\perp M$.
\par Each $\mbox{ }_pg$ is the ${\bf C}_M$ valued function, consequently, up to the
isomorphism it is the complex locally analytic function. Therefore,
its restriction on ${\bf C}_M$ can be written in the form
\par $(4)$ $\mbox{ }_pg(x)|_{{\bf C}_M} = \sum_k c_{k,p}(x-x_0)^k$, where $x_0\in {\bf
C}_M$ is a marked point and $c_{k,p}\in {\bf C}_M$ for each $k$.
Suppose that $\nu $ is some phrase of $f$ in the $z$-representation
of $f$ on $U$ prescribed by Equations $(2,4)$ on the entire domain
$U$ with a common $x_0\in {\bf C}_M$ independent from $p$ so that
the series converge uniformly on $U$. Using the translation
$z\mapsto z-z_0$ we can consider for simplicity, that $x_0=0$. Then
\par $(5)$ $\int_{\gamma } \sum_{k,p} c_{k,p} x^kN_p dx = \sum_{k,p}
[c_{k,p}/(k+1)] (b^{k+1}- a^{k+1})N_p$ when $f$ does not contain any
singularities in $U$, that is ${\cal A}_r$ differentiable in $U$,
where $a=\gamma (0)$, $b= \gamma (1)$, $\gamma (t)$ with $t\in
[0,1]$ is the rectifiable path in ${\bf C}_M$. \par If $z_0\in {\bf
C}_M\cap U$ is an isolated pole of $f$ encompassed by $\gamma $,
then \par $(6)$ $Res (z_0,\nu ).M = \sum_p (c_{-1,p} M)N_p$.
\par Indeed, each Cayley-Dickson number can be written in the polar form also
\par $(7)$ $z=|z|\exp (Arg (z))$, where $Arg (z)$ is a purely imaginary
Cayley-Dickson number, so that $Arg (z) = \alpha M + \beta N_z$ with
$\alpha , \beta \in \bf R$, consequently, $z^k = |z|^k \exp (k Arg
(z))$ for every $k\in \bf Z$. Therefore, $z^k = |z|^k \exp (kM\phi
)$ for each $z\in {\bf C}_M$, where $\phi =\phi (z) \in \bf R$. The
line integral is additive, hence $$\int_{\gamma } \sum_p \mbox{
}_pg(z)N_p dz = \sum_p (\int_{\gamma } \mbox{ }_pg(z)dz)N_p.$$ Then
each term of the form $\eta (z) = \{ az^kb \} _{q(3)}$ on $U$ with
integer $k\ne -1$ has the function $v(z)$ given by any either left
or the right algorithm as $v(z) = \{ az^{k+1}b \} _{q(3)}$ with
$(dv(z)/dz).1 = \{ az^kb \} _{q(3)}$,  since $dv(z)/dz =
\sum_{j=0}^k \{ a((z^j{\bf 1})z^{k-j})b \} _{q(3)}$, where $\bf 1$
denotes the unit operator on ${\cal A}_r$, while $q(3)$ indicates on
the order of multiplication for $r\ge 3$, for $r=2$ the quaternion
algebra is associative. In each ${\bf C}_M$ for $\eta $ with $k= -1$
the function $v$ has the restriction $v|_{{\bf C}_M} = \{ a ~ Ln(z)
~ b \} _{q(3)}$ for $z\ne 0$, since $(de^z/dz).h = e^zh$ and $(d Ln
(z)/dz ).h = z^{-1}h$ for each $0\ne z\in {\bf C}_M$, $h\in {\bf
C}_M$ and the logarithmic function $Ln (z)$ is the inverse function
of the exponential function (see also
\cite{luoyst,luoyst2,ludfejm}). Thus $Res ( z_0, \{ a(z-z_0)^{-1}b
\} _{q(3)}).M = \{ a M b \} _{q(3)}$ for any purely imaginary
Cayley-Dickson number $M$.  Another valuable identities for
calculating residues follow from Formulas 9$(3-6)$.
\par Take the word $\eta (z) =
\{ a_1z^{n_1}...a_k z^{n_k}a_{k+1} \} _{q_1(2k+1)}$ the restriction
of which on ${\bf C}_M$ has the form $\{ b_1...b_{k+1} (z^s {\tilde
z}^m) \} _{q_2(k+2)}$, where constants $a_1,...,b_{k+1}$ belong to
${\cal A}_r$, $2\le r\le 3$, $k\in \bf N$, $n_1,...,n_k, s, m \in
\bf Z$.  Then for the rectifiable loop $\gamma $ encompassing zero
in ${\bf C}_M$  with the winding number $j$ we infer, that
\par $(8)$ $(2\pi j)^{-1} \int_{\gamma } \{ a_1z^{n_1}...a_k z^{n_k}a_{k+1} \}
_{q_1(2k+1)} dz = \{ b_1...b_{k+1} M \} _{q_2(k+2)}$  for $m-s=1$
with $n= -1$;
\par $(9)$ $(2\pi j)^{-1} \int_{\gamma } \{ a_1z^{n_1}...a_k
z^{n_k}a_{k+1} \} _{q_1(2k+1)} dz = \{ b_1...b_{k+1} {\tilde M} \}
_{q_2(k+2)}$  for $s-m=1$ with $n =-1$;
\par $(10)$ $(2\pi j)^{-1}
\int_{\psi } \{ a_1z^{n_1}...a_k z^{n_k}a_{k+1} \} _{q_1(2k+1)} dz =
0$ for either $|m-s|\ne 1$ or $n\ne -1$, where $n=n_1+...+n_k$, $s+m=n$, \\
since $\exp (-2\pi Mt) M \exp (2\pi Mt) = M$ in ${\cal A}_r$, also
$(de^z/dz).h = e^zh$ for $Im (h)\in Im (z){\bf R}$.
\par If $z_1=a + bM$ is some Cayley-Dickson number with $a, b\in \bf
R$ and a purely imaginary number $M$, then any other $z_2\in {\cal
A}_r$ can be written in the form $z_2 = \alpha + \beta M + \phi N$,
where $\alpha , \beta , \phi \in \bf R$, $N$ is a purely imaginary
number orthogonal to $M$, $N\perp M$. Therefore, \par $(11)$ $z_1z_2
= (\alpha + \beta M) z_1 + \phi N {\tilde z}_1$, since $NM = - MN$.
\\ This identity and the formulas given above together with the additivity
and $\bf R$ homogeneousity of the non commutative line integral over
${\cal A}_r$ can be used for calculations of integrals along paths
in planes such as ${\bf C}_M$ and the corresponding residues $Res
(z_0,\mu ).M$.
\par Generally if there is given $f$ and its phrase $\mu $ is
specified, then such transformations to the form $(2,4)$ may change
the phrase, so $\mu $ may be not equal to $\nu $. If $\int_{\gamma
}\nu (z) dz$ or $Res (z_0, \nu ).M$ is calculated, then using
transition formulas from $\mu $ to $\nu $ and vise versa one may
calculate these quantities for $\mu $ if these transition formulas
do not change such integrals and residues, that generally may be not
a case. Also mention that $K\exp (Mt) = K\cos (t) +KM \sin (t) =:
N(t)$ for two purely imaginary octonion numbers $K$ and $M$ with
$|M|=1$ and the real variable $t$. If $K$ and $M$ are perpendicular,
$K\perp M$, and $|K|=1$, then $K\exp (Mt) = \exp (\pi N(t)/2)$ and
hence $Ln (K\exp (Mt)) = \pi N(t)/2 + 2\pi N(t)k$, $k\in \bf Z$
depending on the branch of $Ln$, since $|N(t)|=1$.
\par There is also an interesting particular case of terms \par $ \eta (z) = \{
c_{1,n_1}z^{n_1}...c_{k,n_k}z^{n_k}c_{k+1,n_{k+1}} \} _{q(2k+1)}$,
\\ when there exist $K$ and $M$ purely imaginary numbers and constants
$c_{j,n_j}$ in the Cayley-Dickson algebra ${\cal A}_r$, $2\le r \le
3$, so that $Im (c_{j,n_j}) \perp K{\bf C}_M$ for each $j$, where
$|K|=|M|=1$, $K\perp M$, $Im (z) := (z-{\tilde z})/2$, $n_j\in {\bf
N} = \{ 1, 2, 3,... \}$ for each $j$, $n = n_1+...+n_k>0$. The real
field $\bf R$ is the center of ${\cal A}_r$ and each purely
imaginary number $S$ orthogonal to $K{\bf C}_M$ anti-commutes with
each $z\in K{\bf C}_M$. Then using multiplications of generators and
distributivity of the multiplication $(a+b)z=az+bz$ in ${\cal A}_r$
one finds the restriction $v(z)|_{K{\bf C}_M}$ of a function $v(z)$
which can be reduced to the form $\alpha z^{n+1}/(n+1)$  and it has
the extension on ${\cal A}_r$ such that
\par $(dv(z)/dz).h = \sum_{j=0}^n \alpha ((z^jh)z^{n-j})/(n+1)$
for each $z$ and $h\in {\cal A}_r$, \\
where a constant $\alpha $ is in ${\cal A}_r$, $Im (\alpha )\perp
K{\bf C}_M$, since the octonion algebra is alternative and the
quaternion skew field is associative. Therefore, $(dv(z)/dz).1 =
\eta (z)$ for each $z \in K{\bf C}_M$. On the other hand, each
function given by either the left or right algorithm of integration
is of total degree by $z$ equal to $n+1$ and can be reduced to the
form $\alpha z^{n+1}/(n+1)$ on $K{\bf C}_M$. Evidently $(v(b)+
w(x))-(v(a)+w(x))=v(b)-v(a)$ for each function $w(x)$ with values in
${\cal A}_r$ and the argument $x\in {\cal A}_r\ominus K{\bf
C}_M\ominus {\bf R}$, $Im (x)\perp K{\bf C}_M$. This implies that
\par $(12)$ $\int_{\gamma } \eta (z)dz = \alpha (b^{n+1}-a^{n+1})/(n+1)$,
when $n\in \bf N$, $a$ and $b\in K{\bf C}_M$ and either conditions
of the homotopy Theorem 2.15 \cite{luoyst,luoyst2} are satisfied or
$\gamma ([0,1])\subset K{\bf C}_M$.
\par Another example is of a function $f$ which can be written in
the form $f(z)=(a(z)((b(z)1/(z-y))c(z)))e(z)$ in a neighborhood of
$y\in {\cal A}_r$, where $a(z)$, $b(z)$, $c(z)$ and $e(z)$ are
${\cal A}_r$-holomorphic and $a(y)\ne 0$, $b(y)\ne 0$, $c(y)\ne 0$
and $e(y)\ne 0$, $2\le r\le 3$. Then the residue operator is:
$$(13)\quad Res (y,f).N = (2\pi )^{-1} \lim_{0<\beta \to 0} \int_{\gamma
_{\beta } }(a(z)((b(z)(1/(z-y))c(z)))e(z) dz$$  $$ =
(a(y)((b(y)N)c(y))e(y),$$ where to $a$ in Formula 1$(i)$ here
corresponds $y$, while to $y$ in 1$(i)$ here corresponds $\beta $.
\par In relation with words and phrases there are also some useful
identities. Let $q(n+2)$ be a vector indicating an order of the
multiplication of $n+2$ multipliers, then
\par $(14)$ $(d \{ az^nb \}_{q(n+2)} /dz).h = \sum_{j=0}^{n-1} \{ az^jhz^{n-j-1}b
\}_{q(n+2)}$ \\
for each $h\in {\cal A}_r$, where $z^n$ is treated as the product of
$n$ multipliers $z$. Since the Cayley-Dickson algebra ${\cal A}_r$,
$r\ge 2$, is power associative, then
\par $(15)$ $(d \{ z^n \}_{q(n)} /dz).h = \sum_{j=0}^{n-1} \{ z^jhz^{n-j-1}
\}_{q(n)}= \sum_{j=0}^{n-1} \{ z^jhz^{n-j-1} \} _{q_l(n)} $ \\ for
each $q(n)$, but $q(n)$ is the same for all additives in the sum
independently from $j$, where $q_l(n)$ corresponds to the left order
of brackets, $ \{ a_1...a_n \}_{q(n)} =
((...(a_1a_2)...)a_{n-1})a_n$. In the octonion algebra an order of
multiplications in the term $\{ z^jhz^{n-j-1} \}_{q(n)}$ is not so
important (see \S 9) due to the Moufang identities. \par The
complete differential is $(D\eta (z,{\tilde z})).h = (\partial \eta
(z,{\tilde z})/\partial z).h + (\partial \eta (z,{\tilde
z})/\partial {\tilde z}).h$ for a $(z,{\tilde z})$ (super)
differentiable phrase $\eta (z,{\tilde z})$. Particularly for $\eta
(z,{\tilde z}) = z^n{\tilde z}^m$ the complete differential is the
same for ${\tilde z}^mz^n$, since $z$ and $\tilde z$ commute. This
implies, that
\par $(16)$ $((dz^n/dz).h) {\tilde z}^m) + z^n (d({\tilde z}^m)/d{\tilde z}).h) = {\tilde z}^m ((dz^n/dz).h) +
((d{\tilde z}^m/d{\tilde z}).h)z^n $ \\
for all $z, h\in {\cal A}_r$. \par Thus Formulas $(1-16)$ of this
section and Formulas 9$(M1-M3,1-6)$ can serve for practical
calculations of residues.
\par In the case of Cayley-Dickson algebras calculations of line
integrals are more complicated because of non commutativity, non
associativity and possible behavior of functions around branching
points $z_0$ branches of which generally depend on $M/|M|$ for $M=
Im (z-z_0)\ne 0$.
\par  Evaluate some useful integrals. In the notation of \cite{luoyst,luoyst2} ${\hat f}$
is the operator such that ${\hat f}(z).dz = (dg(z)/dz).dz$ for a
(super) differentiable function $f(z)$ with $(dg(z)/dz).1=f(z)$ in
an open subset $U$ and $g$ is calculated with the help of either
left or right algorithm for each power series. Then \par
$(z^n)^{\hat .}.\Delta z = \Delta (z^{n+1}/(n+1)) + o(\Delta z)$
with $n\in \bf N$ and
\par $(e^z)^{\hat .}.\Delta z = \Delta (e^z) + o(\Delta z)$ for each
$z\in {\cal A}_r$, \par $(1/(1+y))^{\hat .}.\Delta y = \Delta (Ln
(1+y)) +o (\Delta y)$ for $|y|<1$, $y\in {\cal A}_r$, since the
power series \par $Ln (1+y) = \sum_{n=1}^{\infty } (-1)^{n+1}y^n/n$
\\ has real expansion coefficients and uniformly converges in each
ball of radius $0<\rho <1$ with the center at zero, where $\Delta y$
is a sufficiently small increment of $y$. Thus
\par $(17)$ $\int_{\gamma } z^ndz = (b^{n+1}-a^{n+1})/(n+1)$,
\par $(18)$ $\int_{\gamma } e^zdz = e^b - e^a$ for each rectifiable path $\gamma $
such that $\gamma (0)=a$ and $\gamma (1)=b$;
\par $(19)$ $\int_{\gamma } dLn (1+y) = \int_{\gamma } (1+y)^{-1} dy = Ln (1+b)
- Ln (1+a)$ \\ for $\gamma $ such that $|\gamma (t)|<1$ for each
$t\in [0,1]$. On the other hand, if $Im (z)$ and $Im (y)\in M\bf R$
for some marked purely imaginary Cayley-Dickson number $M$, then $z$
and $y$ commute, consequently, $Ln (z(1+y)) = Ln (z) + Ln (1+y)$ for
such $y$ and $z$. Therefore, formula $(19)$ is valid also for each
rectifiable $\gamma $ in the $({\bf C}_M\setminus \{ z\in {\cal
A}_r: ~ z\le 0 \} ) -1$, where ${\bf C}_M = {\bf R}\oplus M{\bf R}$
is the complex plane embedded into the Cayley-Dickson algebra ${\cal
A}_r$. \par Then for any given rectifiable path $\gamma $ in ${\cal
A}_r\setminus Q$, where $Q$ is the slit (cut) sub-manifold used for
construction of the non commutative analog of the Riemann surface
over ${\cal A}_r$, choose a sequence $\gamma _n$ of rectifiable
paths in ${\cal A}_r\setminus Q$ so that $\gamma _n$ is the
combination of paths $\gamma _{n,j}: [b_{j-1},b_j]\to {\cal A}_r$ in
${\bf C}_{M_{n,j}}$, $\gamma _n(t)=\gamma _{n,j}(t)$ for each $t\in
[b_{j-1},b_j]$, $j=1,...,m(n)\in \bf N$, $n\in \bf N$,
$b_0=0<b_1...<b_{m(n)}=1$.  Then for each $\gamma _{n,j}$ and hence
for each $\gamma _n$ Formula $(19)$ is also valid in $({\cal
A}_r\setminus Q) -1$. The limit by $n$ then gives
\par $(20)$ $\int_{\gamma } dLn (z) = \int_{\gamma } z^{-1}dz = Ln (b) - Ln (a)$ \\
for each rectifiable path $\gamma $ in ${\cal A}_r\setminus Q$ with
a chosen branch of $Ln$. This is natural, since $(d Ln (z)/dz).1=
1/z$ for each $z\ne 0$ in ${\cal A}_r$.

\par If $f(z)$ is a univalent function in an open domain $U$
in ${\cal A}_r$ having a power series decomposition $f(z) =
\sum_{n=0}^{\infty }a_nz^n$ uniformly converging in $U$ with real
expansion coefficients $a_n$, then due to $(17)$ we deduce, that
\par $(21)$ $\int_{\gamma } f(z)dz = \sum_{n=0}^{\infty }
a_n(b^{n+1}-a^{n+1})/(n+1)$, \\
where $\gamma : [0,1]\to U$ is a rectifiable path in $U$ with
$\gamma (0)=a$ and $\gamma (1)=b$. In particular this gives:
\par $(22)$ $\int_{\gamma } \sin(z) dz = \cos (a) - \cos (b)$,
\par $(23)$ $\int_{\gamma } \cos (z) dz = \sin (b) - \sin (a)$, \\
where $\sin (z) = \sum_{n=0}^{\infty } (-1)^nz^{2n+1}/(2n+1)!$,
$\quad \cos (z) = \sum_{n=0}^{\infty } (-1)^nz^{2n}/(2n)!$.
\par If $f$ is a (super) differentiable function of the octonion
variable not equal to zero on a open set $U$, then $f(z)(1/f(z)) =1$
on $U$ and $((df(z)/dz).h) (1/f(z)) + f(z)((d[1/f(z)]/dz).h)$ for
each $z\in U$ and $h\in {\bf O}$, consequently, $(d[1/f(z)]/dz).h =
- [1/f(z)] ( (df(z)/dz).h) [1/f(z)]$ on the corresponding domain
(see Proposition 3.8.2 \cite{luoyst2}). The Cayley-Dickson algebra
is power associative, so repeating this for $z^n$ with $n\in \bf N$
instead of $f$ we get, that $(dz^{-n}/dz).h = -
z^{-n}((dz^n/dz).h)z^{-n}$ for each $z\ne 0$ and $h$ in ${\cal A}_r$
with $Im (z)$ and $Im (h) \in M{\bf R}$ for some marked purely
imaginary number $M$. In particular, $(dz^n/dz).1 = n z^{n-1}$ for
each integer number $n$ with $z\ne 0$, when $n<0$, $dz^0/dz=0$. Thus
\par $(24)$ $\int_{\gamma } z^ndz = (b^{n+1}-a^{n+1})/(n+1)$ for each integer $n\ne -1$\\
an for every rectifiable path $\gamma $ with $a=\gamma (0)$,
$b=\gamma (1)$, so that $\gamma ([0,1])$ does not contain $0$, when
$n<0$. \par If $\alpha \in {\cal A}_r$ then we define $z^{\alpha }
:= \exp (\alpha Ln (z))$ for each $z\ne 0$. So we calculate its
derivative \par $(dz^{\alpha }/dz). h = (de^y/dy)|_{y=\alpha Ln (z)}
.(\alpha [(dLn (z)/dz).h)])$ for each $h\in {\cal A}_r$.
Particularly,
\par $(dz^{\alpha }/dz).1 = (de^y/dy)|_{y = \alpha Ln (z)}. (\alpha
z^{-1})$. Take $\alpha \in \bf R$, then $(dz^{\alpha }/dz).1= \alpha
z^{\alpha -1}$, hence \par $(25)$ $\int_{\gamma } z^{\alpha }dz =
(b^{\alpha +1} - a^{\alpha +1})/ (\alpha +1)$ for each real $\alpha
\ne -1$ and the rectifiable path in ${\cal A}_r\setminus Q$.
\par If $\alpha \in \bf R$, then we infer the power series
decomposition \par $(1+z)^{\alpha } =\sum_{n=0}^{\infty } {\alpha
\choose n}z^n$ for each $|z|<1$, where ${\alpha \choose n} = \alpha
(\alpha -1)...(\alpha -n+1)/n!$. Particularly, for $\alpha = -1/2$
and $z= - y^2$ we deduce, that
\par $(26)$ $\int_{\gamma } (1/\sqrt{1-y^2})dy = \arcsin (b) -
\arcsin (a)$ \\ for $\gamma ([0,1])\subset [{\cal A}_r\setminus \{
z\in {\cal A}_r: Re (z) =0 \} ]$, where the square root branch is
taken $\sqrt{x}>0$ for each $x>0$, $|\gamma (t)|<1$ for each $t\in
[0,1]$.
\par For the tangent function its power series has the form:
\par $\tan (z) = \sum_{n=1}^{\infty } 2^{2n} (2^{2n} -1) B_n z^{2n-1}
/ (2n)!$ and \\
the cotangent function is:
\par $z \cot (z) = 1 - \sum_{n=1}^{\infty } 2^{2n} B_n z^{2n}
/ (2n)!$, \\
where $B_n$ are the Bernoulli numbers and these power series
absolutely converge in the balls $|z| <\pi /2$ and $|z|<\pi $
respectively. The Bernoulli numbers appear from the generating
function
\par $x/ (e^x-1) = 1 + \sum_{n=1}^{\infty } \beta _nx^n/n!$ for each
$x\in (-\delta ,\delta ) \subset \bf R$ with sufficiently small
$\delta
>0$ so that $\beta _{2n} = (-1)^{n-1}B_n$, $~ \beta _{2n+1}=0$ for
each $n\ge 1$ (see \S XII.4 (445, 449) in \cite{fihteng}). Due to
Formula $(25)$ and the periodicity of the trigonometric functions we
get:
\par $(27)$ $\int_{\gamma } \cos ^{-2} (z) dz = \tan (b) - \tan
(a)$, when $|\gamma (t) - \pi m|<\pi /2$ for each $t\in [0,1]$ for
some marked $m\in \bf Z$,
\par $(28)$ $\int_{\gamma } \sin ^{-2} (z) dz = \cot (a) - \cot
(b)$, when $|\gamma (t) - \pi m - \pi /2|<\pi /2$ for each $t\in
[0,1]$ for some marked $m\in \bf Z$. On the other, hand $\tan (z)$
and $\cot (z)$ are meromorphic functions so that $(d \tan (z)/dz).1
= \cos ^{-2}(z)$ for each $z \in V := \{ y\in {\cal A}_r: Re (y)\ne
\pi (m+1/2), m\in {\bf Z} \} $, also $(d \cot (z)/dz).1 = - \sin
^{-2}(z)$ for each $z \in U := \{ y\in {\cal A}_r: Re (y)\ne \pi m,
m\in {\bf Z} \} $. Therefore, Formulas $(27,28)$ extend on any
rectifiable path $\gamma $ in $V$ or $U$ correspondingly.
\par Mention also that Formula $(18)$ implies:
\par $(29)$ $\int_{\gamma } \cosh (z) dz = \sinh (b) - \sinh (a)$,
\par $(30)$ $\int_{\gamma } \sinh (z) dz = \cosh (b) - \cosh (a)$,
where $\cosh (z) =(e^z + e^{-z})/2$, $\sinh (z) = (e^z - e^{-z})/2$.
\par For $\coth (z) = \cosh (z)/\sinh (z)$ the power series
\par $z \coth(z) = 1 + \sum_{n=1}^{\infty } (-1)^{n-1} 2^{2n} B_n z^{2n}
/ (2n)!$ \\
absolutely converges in the ball $|z|<\pi $, hence
\par $(31)$ $\int_{\gamma } \sinh ^{-2} (z)dz = \coth (a) - \coth
(b) $ \\
for $0<|\gamma (t)|<\pi $ for each $t\in [0,1]$. We can mention that
if $y$ is written in the polar form $y= \rho e^S$, where $\rho
=|y|\ge 0$ and $S$ is a purely imaginary Cayley-Dickson number, then
$y^2=-1$ if and only if $\rho =1$ and $|S| = \pi (m+1/2)$ with $m=
0, 1, 2,3,...$.
\par Since $\tanh (z)$ and $\coth (z)$ are locally analytic and $(d \tanh
(z)/dz).1 = \cosh ^{-2}(z)$ for each $z\in V:= {\cal A}_r\setminus
\{ y: Re (y)=0, |Im(y)| = \pi (m+1/2), ~ m=0, 1, 2, 3,... \}$ and
$(d\coth (z)/dz).1 = - \sinh ^{-2}(z)$ for each $z\in U := {\cal
A}_r\setminus \{ y: Re (y)=0, |Im (y)| = \pi m, ~ m = 0, 1, 2, 3,...
\} $, then $(31)$ extends on $\gamma \subset U$ and
\par $(32)$ $\int_{\gamma } \cosh ^{-2}(z) dz = \tanh (b) - \tanh (a) $
for each $\gamma \subset V$.
\par The series $1/(1+z^2) = \sum_{n=0}^{\infty } {{-1/2}\choose n} z^{2n}/n!$
absolutely converges for $|z|<1$ and its integral gives $\arctan $
in this ball. The function $\arctan (z)$ is locally analytic and
$(d\arctan (z)/dz).1 = 1/(1+z^2)$ for $z^2\ne -1$, since $(d \tan
(z)/dz).1 = 1/\cos ^2(z)$ and $(df^{-1}(z)/dz)|_{z=f(y)} =
(df(y)/dy)^{-1}$ on the corresponding domains (see Proposition 3.8.1
in \cite{luoyst2}). Thus
\par $(33)$ $\int_{\gamma } [1/(1+z^2)] dz = \arctan (b) - \arctan
(a)$ for $\gamma \subset \{ y\in {\cal A}_r: y^2\ne -1 \} $.
\par In the octonion algebra due to its alternativity
$(d Ln (z+ \sqrt{z^2+\alpha })/dz).1 = 1/\sqrt{z^2+\alpha }$ for a
non zero real number $\alpha $ and each $z^2\ne -\alpha $ in ${\bf
O}$, where $Ln (z+ \sqrt{z^2+\alpha })$ and $\sqrt{z^2+ \alpha }$
are locally analytic. Therefore, with the help of Formula $(20)$ we
infer:
\par $(34)$ $\int_{\gamma } [1/\sqrt{z^2+\alpha }]dz = \int_{\gamma } d Ln (z+ \sqrt{z^2+\alpha
}) $ \par $ = Ln (b+ \sqrt{b^2+\alpha }) - Ln (a+ \sqrt{a^2+\alpha
})$ \\ for each path $\gamma $ in the set $U := \{ y\in {\bf O}:
y^2\ne - \alpha , Re (y)\ne 0 \} $, since in it branches of
functions $1/\sqrt{z^2+\alpha }$ and $Ln (z+ \sqrt{z^2+\alpha })$
are specified. One mentions that $y^2= - \alpha $ if and only if
$\rho ^2 = |\alpha |$ and $|S| = \pi (m + \kappa (\alpha )/2)$ for
$m\in \bf Z$, where $y=\rho e^S$, $\rho =|y|$, $Re (S)=0$, $\kappa
(\alpha ) =1$ for $\alpha
>0$, $\kappa (\alpha ) = 0$ for $\alpha \le 0$.
\par Consider now the general algorithm in more details.
Let $f_1$ and $f_2$ be two analytic functions of $z$ on $U$. Denote
$f^0 =f$, $f^{-n} = f^{(n)}$, $f^{(n)}(z) :=
(d^nf(z)/dz^n).1^{\otimes n}$, $f^n$ is such that
$(df^n(z)/dz).1=f^{n-1}(z)$ for each $n\in \bf N$. Particularly,
$(\{ az^kb \}_{q(3)})^n = \{ a z^{k+n} b \}_{q(3)}
[(k+1)...(k+n)]^{-1}$ for constants $a, b\in {\cal A}_r$, $k=0, 1,
2,...$, where the symbol $z^0$ is also integrated for convenience so
that each phrase of a $z$ (super) differentiable function is a
series of addends having obligatory components $z^n$ or $e= {\bf
1}(1)$ (see also \S 2.14 \cite{norfamlud}). In view of this
$(df^1(z)/dz).h = {\hat f}.h$ and $f\mapsto f^l$ is the
anti-derivation operation of order $l$, $l\in \bf N$. Then the left
algorithm is
\par $(35)$ $(f_1f_2)^1 = \sum_{s=0}^{n_2} (-1)^s
f_1^{1+s}f_2^{-s}$, where $n_2$ is the least natural number such
that $f_2^{-n_2-1}=0$. The right algorithm is symmetric to the left.
The work with phrases and power series reduces to products of
polynomials and their sums, so consider this algorithm for $(\{
f_1...f_k \} )^l$, where $l\in \bf N$, $f_2$,...,$f_k$ are
polynomials by $z$ of degrees $n_2,...,n_k$ respectively. Then \par
$(36)$ $(f_1f_2)^2 = \sum_{s_1=0}^{n_2} (-1)^{s_1}
(f_1^{1+s_1}f_2^{-s_1})^1 =$ \par $ \sum_{s_1=0}^{n_2}
\sum_{s_2=0}^{n_2-s_1} (-1)^{s_1+s_2} f_1^{2+s_1+s_2} f_2^{-s_1-s_2}
=\sum_{s=0}^{n_2}(-1)^s (s+1) f_1^{2+s}f_2^{-s}$. \\
Continue these calculations by induction using that
\par $\sum_{s_1=0}^s \sum_{s_2=0}^{s-s_1} ...
\sum_{s_p=0}^{s-s_1-...-s_{p-1}} 1 = {{s+p-1}\choose s}$, \\
where ${n\choose  m} = n!/(m!(n-m)!)$ denotes the binomial
coefficient. Then
\par $(37)$ $(f_1f_2)^l = \sum_{s=0}^{n_2} (-1)^s {{s+l-1}\choose s}
f_1^{s+l} f_2^{-s}$. \\
Applying Formula $(37)$ to product of $k$ terms we deduce that \\
$(38)$ $( \{ f_1f_2...f_k \}_{q(k)})^l = \sum_{s_k=0}^{n_k}
(-1)^{s_k} {{s_k+l-1}\choose {s_k}} \{ [f_1...f_{k-1}]^{s_k+l}
f_k^{-s_k} \}_{q(k)} = \sum_{s_k=0}^{n_k}\sum_{s_{k-1}=0}^{n_{k-1}}$ \\
$ ... \sum_{s_2=0}^{n_2} (-1)^{s_2+...+s_k} {{s_k+l-1}\choose {s_k}}
{{s_{k-1} + s_k+l-1}\choose {s_{k-1}}}... {{s_2+...+s_k+ l
-1}\choose {s_2}} \{ f_1^{l+s_2+...+s_k}
f_2^{-s_2}...f_k^{-s_k} \} _{q(k)}$, \\
where the notation $[f_1...f_{k-1}]^{s_k+l}$ does not mean any order
of multiplications but only that the anti-derivation operator is of
order $s_k+l$ for a collection of these terms in square brackets.
For the skew field of quaternions curled brackets can be omitted due
to the associativity. Thus either left or right algorithm specifies
the branch of the non commutative line integral for $f$ (see also
\cite{luoyst,luoyst2,norfamlud}) and for this Formulas $(35-37)$ are
helpful.
\par If $f(z)$ is an analytic function on a ball $B$ with the center at zero
and the uniformly converging power series of $f$ has real expansion
coefficients, also $\alpha \in \bf O$, $\alpha \ne 0$, then the
mapping $z\mapsto \alpha z\alpha ^{-1}$ in $\bf O$ induces the
mapping $f(z)\mapsto \alpha f(z)\alpha ^{-1}$ so that $\alpha
f(z)\alpha ^{-1} = f(\alpha z\alpha ^{-1})$. But the line integral
changes as \par $(39)$ $\alpha (\int_{\gamma } f(z).dz)\alpha ^{-1}
= \int_{\gamma } f(\alpha z\alpha ^{-1}).(\alpha dz \alpha ^{-1}) =
\int_{\psi } f(y).dy$
\\ for such function $f$, since the octonion algebra is alternative and $\bf R$ is its center, where
$\psi (t) = \alpha \gamma (t) \alpha ^{-1}$ for each $t\in [0,1]$.
If $f$ is a multi-valued function, for example, $Ln (z)$ or
$\sqrt[n]{z}$, $n\ge 2$, then the mapping $z\mapsto \alpha z$ or
$z\mapsto z\alpha $ may cause a transition from one branch of $f$
into another.
\par {\bf 15. Proposition.} {\it Let $f_1$ and $f_2$ be two (super)
differentiable univalent functions on an open set $U$ in the
Cayley-Dickson algebra ${\cal A}_r$ and $\gamma $ be a rectifiable
path in $U$. Then \par
\par $\int_{\gamma } f_1(z)[(df_2(z)/dz).1]dz = f_1(z)f_2(z)|_a^b -
\int_{\gamma } [(df_1(z)/dz).1]f_2dz$.}
\par {\bf Proof.} Consider the (super) derivative
\par $(df_1(z)f_2(z)/dz).h = [(df_1(z)/dz).h] f_2(z) +
f_1(z)[(df_2(z)/dz).h]$, \\ where $z\in U$ and $h\in {\cal A}_r$.
Then we take the integral sums $\sum_{j=1}^n$ corresponding to
partitions $P$ of $\gamma $. Therefore, in the notation of \S 11 we
get
\par $\int_{\gamma } [(df_1(z)f_2(z)/dz).1]dz = \lim_{\delta (P)\to 0} \sum_{j=1}^n
(df_1(z)f_2(z)/dz)|_{z=z_j}.\Delta _jz$ \par $= f_1(z)f_2(z)|_a^b =
\lim_{\delta (P)\to 0} \sum_{j=1}^n
\widehat{([(df_1(z)/dz).1]f_2(z))}.\Delta _jz$ \par $ + \lim_{\delta
(P)\to 0} \sum_{j=1}^n \widehat{(f_1(z)[(df_2(z)/dz).1])}.\Delta _j
z$ \par $ = \int_{\gamma }[(df_1(z)/dz).1]f_2(z)dz + \int_{\gamma
}f_1(z)[(df_2(z)/dz).1] dz$, consequently,
\par $\int_{\gamma } f_1(z)[(df_2(z)/dz).1]dz = f_1(z)f_2(z)|_a^b -
\int_{\gamma } [(df_1(z)/dz).1]f_2dz$.

\section{Affine algebras over octonions}
\par {\bf 1. Definition.} If a generalized Cartan matrix $A$
has positive all proper main minors and $det (A)=0$, then one says
that $A$ is of the affine type. The algebra ${\sf g}(A)$ (see \S
2.6) corresponding to the generalized Cartan matrix $A$ of the
affine type is called the affine algebra over ${\cal A}_r$.
\par {\bf 2. Remarks.} Let $f$ and $g$ be two meromorphic
functions on an open set $V$ in the Cayley-Dickson algebra with
singularities at a set $W$ satisfying conditions 3.1$(R1-R3)$. Then
\par $(df(z)g(z)/dz).h = [(df(z)/dz).h] g(z) + f(z)[(dg(z)/dz).h]$
\\ for each $h\in {\cal A}_r$ and every $z\in V\setminus W=: U$. In particular,
take $h=1$. On the other hand, $\int_{\gamma } [(df(z)g(z)/dz).1]dz
= 0$ for each rectifiable loop $\gamma $ in $U$ and we get \par
$(1)$ $Res (z_0, (f'(z).1)g(z)).M = - Res (z_0, f(z)(g'(z).1)).M$ \\
for each $z_0\in V$ and every purely imaginary Cayley-Dickson number
$M$ (see also \S \S 3.1, 3.6 and 3.15). In the notation $\psi
_{z_0}(f',g) := Res (z_0, (f'(z).1)g(z))$ Formula $(1)$ takes the
form
\par $(2)$ $\psi _{z_0}(f',g) = - \psi _{z_0}(f,g')$ for each $z_0\in
U$, where $\psi _{z_0}(f',g).M$ is ${\cal A}_r$ additive by $f$ and
by $g$ also, by purely imaginary Cayley-Dickson number $M$,
moreover, it is $\bf R$ homogeneous by $f$ and by $g$ and by $M$,
which follows from the properties of the residue operator (see
Theorem 3.6 also).
\par Consider the set theoretical composition of mappings $f_1\circ
f_2(y) = f_1(f_2(y))$. In view of Proposition 2.2.1 and Theorems
2.11, 2.15, 2.16 and 3.10 \cite{luoyst2} we can reformulate the
result.
\par {\bf 2.1. Proposition.} {\it Let $g: U\to {\cal A}_r^m$, $r\ge 2$,
and $f: W\to {\cal A}_r^n$ be two differentiable functions on $U$
and $W$ respectively such that $g(U)\subset W$, $U$ is open in
${\cal A}_r^k$, $W$ is open in ${\cal A}_r^m$, $k, n, m\in \bf N$,
where $f$ and $g$ are simultaneously either $(z,{\tilde z})$, or
$z$, or $\tilde z$-differentiable. Then the composite function
$f\circ g(z) := f(g(z))$ is differentiable on $U$ and
\par $(3)$ $(Df\circ g(z)).h = (Df(g)).((Dg(z)).h)$ \\
for each $z\in U$ and each $h\in {\cal A}_r^k$, and hence $f\circ g$
is of the same type of differentiability as $f$ and $g$.}
\par {\bf 2.2. Multiplications of operators.} The set theoretic composition serves as the forgetful
functor on an order of the derivatives: \par $(4)$ $(d(f_3\circ
f_2)\circ f_1(x)/dx).h = (df_3\circ
f_2(y)/dy)|_{y=f_1(x)}.[(df_1(x)/dx).h] = $ \par $ (df_3\circ
f_2(y)/dy)|_{y=f_1(x)}.v|_{v=(df_1(x)/dx).h}$ \par $ =
[(df_3(z)/dz)|_{z=f_2(y)}.[(df_2(y)/dy)|_{y=f_1(x)}.v]|_{v=(df_1(x)/dx).h}$
\\ $=(df_3(z)/dz)|_{z=f_2(y)}.[(df_2(y)/dy)|_{y=f_1(x)}.[(df_1(x)/dx).h]]=
(df_3(f_2(f_1(x)))/dx).h$
\\ for every $h\in
{\cal A}_r^{n_1}$, for open domains $U_j$ in ${\cal A}_r^{n_j}$ for
$j=1, 2, 3$ so that $f_j: U_j\to U_{j+1}$ for $j=1, 2$, $f_3: U_3\to
{\cal A}_r^n$, $n_1, n_2, n_3, n\in \bf N$, $f_j$ is $z$
differentiable on $U_j$ for each $j=1, 2,3$. If consider ${\cal
A}_r$ additive $\bf R$ homogeneous operators $A_j$ on vector spaces
$X_j$ over ${\cal A}_r$, $A_j: X_j\to X_{j+1}$, with another type of
composition $A_2\diamond A_1$ induced by that of the matrix
multiplication of matrices with entries in ${\cal A}_r$, also
Cayley-Dickson numbers and also by the multiplication of matrices on
the corresponding vectors, then it will be necessary to describe an
order of such new multiplication indicating a vector $q(m)$
prescribing such order $ \{ A_m\diamond ... \diamond A_1 \}_{q(m)}$,
where $A_j\diamond h_j:= A_j(h_j)$ for each $h_j\in X_j$. The set
theoretic composition of operators
\par $(5)$ $(A_3\circ A_2)\circ A_1(h) = (A_3\circ A_2)(v)|_{v=A_1(h)}$
\par $ = A_3(A_2(v))|_{v=A_1(h)} = A_3(A_2(A_1(h)))$  for each $h\in X_1$ \\
is always associative and they correspond as
\par $(6)$ $A_m\circ ... \circ A_1(h) = A_m\diamond (A_{m-1} \diamond ...
(A_2\diamond (A_1\diamond h))...)$ \\
for each $h\in X_1$, particularly, \par $(7)$ $(A_2\circ A_1)(h) =
A_2\diamond (A_1\diamond h)$ \\ such that the set theoretic
composition of operators can be considered as the particular case of
more general multiplication of operators with the right composition
on the right side of Formulas $(6,7)$. Mention that generally
$(A_2\diamond A_1) \diamond h $ may be not equal to $A_2\diamond
(A_1 \diamond h)$ over $\bf O$. Over the quaternion skew field $\bf
H$ brackets can be omitted, since $\bf H$ is associative.
\par {\bf 2.3.}  Take the algebra ${\cal M}(V,W)$ of all
meromorhic functions in $V$ with singularities contained in $W$,
where $V$ is open in $\hat {\cal A}_r$ and $W$ satisfy Conditions
3.1$(R1-R3)$, $W\subset V$. In particular, for a singleton $\{ 0 \}
$ instead of $W$ this gives the algebra ${\cal L}(V)$ of all Loran
series on the set $V\setminus \{ 0 \} $ open in ${\cal A}_r$:
\par $(8)$ $\eta (z) = \sum_{m, j, p} \{ A_{m,j}, z^m \} _{q(2p)}$ \\
with center at zero, where
\par $(9)$ $ \{ A_{m,j}, z^m
\}_{q(2p)} : = \{ a_{m,m_1,1}z^{m_1}...a_{m,m_p,p}z^{m_p} \}_{q(2p)}$, \\
$a_{m,k,p}\in {\cal A}_r$ and $m_k\in {\bf Z}$ for each $k\in \bf
N$, $j, p\in \bf N$, $m=(m_1,...,m_p)$ (see \cite{luoyst,luoyst2}).
To each phrase $\eta (z)$ the function $f_{\eta }(z) = ev_z(\eta )
=\eta (z)$ corresponds, where $ev_z$ denotes the valuation operation
at $z\in U$.
\par There is the natural equivalence relation of phrases
prescribed by the rules:
\par $(10)$ $b\eta _1\eta _2 = \eta _1 b \eta _2 = \eta _1 \eta _2
b$  for each real number $b$,
\par $(11)$ $(z^lz^n) = (z^{n+l})$ for all integers $l, n$,
\par $(12)$ $(\eta _1 +\eta _2) - \eta _2 = \eta _1$ for all phrases
$\eta _1$ and $\eta _2$, since ${\cal A}_r$ is the power associative
algebra with the center $\bf R$. Another equivalences are caused by
associativity of the quaternion skew field $\bf H$ and the
alternativity of the octonion algebra $\bf O$. Over $\bf H$ brackets
in phrases can be omitted.
\par If $W$ is not a singleton, then the restriction of $f\in {\cal M}(V,W)$
on  $z_0+{\bf C}_M$ has the Loran series $(9)$ by $z-z_0\in {\bf
C}_M$ instead of $z$ for each $z_0\in W$ and each purely imaginary
Cayley-Dickson number $M\in {\cal A}_r$. This follows from Formulas
$(3.22,23)$ \cite{luoyst,luoyst2} applied to rectifiable loops
$\gamma _1$ and $\gamma _2$ in $z_0+{\bf C}_M$ and Conditions
3.1$(R1-R3)$.

\par {\bf 2.3. Definitions.} For a
generalized Cartan matrix $A$ take the matrix ${\AA}$ obtained from
$A$ by deleting its zero column and zero row and put $ \check{\sf
g}= {\sf g}({\AA })$. Define the wrap algebra \par ${\cal
M}(\check{\sf g}) : = {\cal M}(V,W) \otimes_{{\cal A}_r} \check{\sf
g}$ with the multiplication:
\par $[P\otimes x, Q\otimes y]_0 := PQ \otimes [x,y]$ for all
$P, Q\in {\cal M}(V,W)$ and $x, y \in \check{\sf g}$.
\par Extend the form $(x|y)$ from $\check{\sf g}$ onto ${\cal M}(\check{\sf
g})$ by the formula:
\par $(P\otimes x|Q\otimes y)_z = PQ(x|y)$.
\par {\bf 3. Proposition.} {\it There exists an ${\cal A}_r$ valued graded $2$-cocycle
on the wrap algebra ${\cal M}(\check{\sf g})$ over ${\cal A}_r$,
$2\le r \le 3$.}
\par {\bf Proof.} The differentiation \par $(1)$ ${\sf D} P(z) := (dP(z)/dz).1$ of the algebra
${\cal M}(V,W)$ for $z\in V\setminus W$ \\ we extend to the
differentiation of ${\cal M}(\check{\sf g})$ by the formula
\par $(2)$ ${\sf D} (P\otimes x) := ({\sf D} P)\otimes x$.
\par Define the mapping
\par $(3)$ $\omega (z_0;P\otimes x, Q\otimes y)
:= Res (z_0, ({\sf D} P)Q).(x|y)$ \\
for each $x, y \in \check{\sf g}$ and all $P, Q\in {\cal M}(V,W)$
and every $z_0\in W$, where the form $(*|*)$ is given by Theorem
2.20. Its restriction $(*|*)|_{\mbox{ }_0{\sf g}}$ on $\mbox{
}_0{\sf g}$ is real valued. We have
\par $(4)$ $Res (z_0, f).b=0$ for each real number $b$ \\ and every
$f \in {\cal M}(V,W)$ with $z_0\in V$, since in this case the line
integral is calculated along a rectifiable loop $\gamma $ contained
in $\bf R$, hence
\par $(5)$ $Res (z_0,f).h = Res (z_0, f).(Im (h))$ \\ for each
Cayley-Dickson number $h\in {\cal A}_r$. For the definiteness of
residues we consider them for specified phrases of functions.
Consider the ordered product $\{ f_1f_2f_3 \}_{q(3)}$ (usual
point-wise) of $z$ differentiable functions $f_1, f_2, f_3$ on the
open set $U$ in ${\cal A}_r$. For each $z_0\in U$ we have $Res
(z_0,(\{ f_1f_2f_3 \}_{q(3)})'.1) =0$, consequently, \par $(6)$ $Res
(z_0,\{ (f_1'.1)f_2f_3 \}_{q(3)}) + Res (z_0,\{ f_1(f_2'.1)f_3
\}_{q(3)})$ \par $ + Res (z_0,\{ f_1f_2(f_3'.1) \}_{q(3)}) =0$.
\par There is the decomposition
\par $(7)$ $P = \mbox{ }_0P i_0+...+ \mbox{ }_{2^r-1}P i_{2^r-1}$, \\
where $\mbox{ }_jP$ is the $\bf R$ valued function for each $j$, so
that they are related by the formulas
\par $(8)$ $\mbox{ }_0P = (P + (2^r-2)^{-1} \{ - P +
\sum_{j=1}^{2^r- 1} i_j(Pi_j^*) \} )/2$,
\par $(9)$ $\mbox{ }_kP = (i_k (2^r-2)^{-1} \{ - P + \sum_{j=1}^{2^r- 1}
i_j(Pi_j^*) \} - Pi_k )/2$ for each $k=1,...,2^r-1$. \\ We have the
natural identity
\par $(10)$ $\mbox{ }_kP i_k \otimes \mbox{ }_kv = \mbox{ }_kP \otimes \mbox{
}_kvi_k$. \par On the other hand, \par $(11)$ $\omega (P\otimes x,
Q\otimes y) = \psi _0((DP), Q).(x|y)$, consequently, for pure states
$\mbox{ }_ka i_k = P_k\otimes \mbox{ }_kx i_k\in \mbox{ }_k{\cal
M}(\check{\sf g})i_k$ and $bi_j = Q_j\otimes \mbox{ }_jy i_j\in
\mbox{ }_j{\cal M}(\check{\sf g})i_j$ we get
\par $(12)$ $\omega (\mbox{ }_kP\otimes \mbox{ }_kx i_k, \mbox{ }_jQ\otimes \mbox{ }_jy i_j)= (-1)^{\eta
(k,j)+1} \omega (\mbox{ }_jQ\otimes \mbox{ }_jy i_j, \mbox{ }_kP\otimes \mbox{ }_kx i_k)$, \\
since the form $(x|y)$ is symmetric on $\mbox{ }_0{\sf g}(A)$ and
hence on $\mbox{ }_0\check{\sf g}$ and satisfies 17$(2)$. \par
 At the same time we have the identities:
\par $(13)$ $ Res (z_0, ((\mbox{ }_kP\mbox{ }_jQ)'.1)\mbox{ }_sR)
+ Res (z_0, ((\mbox{ }_jQ\mbox{ }_sR)'.1)\mbox{ }_kP) + Res (z_0,
((\mbox{ }_sR\mbox{ }_kP)'.1)\mbox{ }_jQ)$
\par $ = 2 (Res (z_0, ((\mbox{ }_kP)'.1)\mbox{ }_jQ\mbox{ }_sR) + Res (z_0,
\mbox{ }_kP((\mbox{ }_jQ)'.1)\mbox{ }_sR) $
\par $+ Res (z_0, \mbox{ }_kP \mbox{ }_jQ ((\mbox{ }_sR)'.1))) = 2 Res(z_0,(\mbox{ }_kP \mbox{ }_jQ \mbox{ }_sR)'.1)= 0
$ \\ and for pure states $\mbox{ }_ka i_k= \mbox{ }_kP\otimes \mbox{
}_kx i_k$, $\mbox{ }_jb i_j = \mbox{ }_jQ\otimes \mbox{ }_jy i_j$,
$\mbox{ }_sc i_s = \mbox{ }_sR\otimes \mbox{ }_sz i_s$ we deduce
taking into account Theorem 2.20 and Properties 17$(1-5)$ and
Identity $(13)$:
\par $(14)$ $\omega (z_0; [\mbox{ }_kai_k,\mbox{ }_jbi_j],\mbox{ }_sci_s) +
(-1)^{\xi (k,j,s)} \omega (z_0;[\mbox{ }_jbi_j,\mbox{
}_sci_s],\mbox{ }_kai_k)$\par $ + (-1)^{\xi (k,j,s) + \xi (j,s,k)}
\omega (z_0;[\mbox{ }_sci_s,\mbox{ }_kai_k],\mbox{ }_jbi_j) $
\par $ = Res (z_0, (D(\mbox{ }_kP \mbox{ }_jQ))\mbox{ }_sR).([\mbox{ }_kx i_k,\mbox{ }_jy i_j]|\mbox{ }_sz i_s)
+$ \par $ (-1)^{\xi (k,j,s)} Res (z_0, (D(\mbox{ }_jQ\mbox{
}_sR))\mbox{ }_kP).([\mbox{ }_jy i_j,\mbox{ }_sz i_s]|\mbox{ }_kx
i_k)$ \par $+ (-1)^{\xi (k,j,s) + \xi (j,s,k)} Res (z_0, (D(\mbox{
}_sR \mbox{ }_kP))\mbox{ }_jQ).([\mbox{ }_szi_s,\mbox{
}_kxi_k]|yi_j)$
\par $= 2 Res(z_0,D(\mbox{ }_kP \mbox{ }_jQ \mbox{ }_sR)).([\mbox{ }_kxi_k,\mbox{ }_jyi_j]|\mbox{ }_szi_s) =0$. \par Thus
Formulas $(12, 14)$ mean that $\omega $ is the graded $2$-cocycle
over the Cayley-Dickson algebra ${\cal A}_r$.
\par {\bf 4. Definitions.} Put $\check{\cal M}(\check{\sf g}) := {\cal M}(\check{\sf g})\oplus F(W,{\cal
A}_r)K$ with the multiplication
\par $(1)$ $[a \oplus \alpha (z_0) K, b \oplus \beta (z_0) K] := [a,b]_0 + \omega (z_0;a,b)K$
for each $a, b \in {\cal L}(\check{\sf g})$ and $\alpha , \beta \in
F(W,{\cal A}_r)$, $z_0\in W$ \\ where $K$ is the additional
generator, $F(W,{\cal A}_r)$ denotes the ${\cal A}_r$ vector space
of all functions $\alpha , \beta : W\to {\cal A}_r$.
\par We introduce the operator $d_lP(z) := - (dP(z)/dz).z^{l+1}$ for each $P\in {\cal M}(V,W)$
and all $z\in V\setminus W$, also $d_lK := 0$, $d_l(P\otimes x) =
d_l(P)\otimes x$, where $l\in \bf Z$.
\par Define new algebra ${\hat {\cal M}}(\check{\sf g}) = \check{\cal M}(\check{\sf g})\oplus {\cal A}_r{\sf d}$
with the differentiation ${\sf d} P(z) := (dP(z)/dz).z = - d_0P(z)$
for each $z\in U\setminus \{ 0 \} $ and every $P\in {\cal M}(U)$ and
${\sf d}K :=0$. We put as the multiplication
\par $(2)$ $[\mbox{ }_kP\otimes \mbox{ }_kxi_k\oplus \alpha i_kK \oplus ei_k{\sf d},
\mbox { }_jQ \otimes \mbox{ }_jyi_j\oplus \beta i_jK \oplus ti_j{\sf
d}] := $
\par $(\mbox{ }_kP\mbox { }_jQ\otimes [\mbox{ }_kxi_k, \mbox{ }_jyi_j] \oplus
ei_k({\sf d} \mbox{ }_jQ)\otimes \mbox{ }_jyi_j - (-1)^{\eta (k,j)}
t i_j({\sf d} \mbox{ }_kP)\otimes \mbox{ }_kxi_k) \oplus \omega
(z_0; \mbox{ }_kP\otimes \mbox{ }_kxi_k, \mbox{
}_jQ\otimes \mbox{ }_jyi_j)$ \\
for pure states $\mbox{ }_kP\otimes xi_k\oplus \alpha i_kK \oplus
ei_k{\sf d} \in \mbox{ }_k{\hat {\cal M}}(\check{\sf g})i_k$ and
$\mbox { }_jQ \otimes yi_j\oplus \beta i_jK \oplus ti_j{\sf d} \in
\mbox{ }_j{\hat {\cal M}}(\check{\sf g})i_j$, where $\alpha
(W)\subset {\bf R}$, $\beta (W)\subset {\bf R}$, $e, t\in \bf R$.
Extend this multiplication on all elements of ${\hat {\cal
M}}(\check{\sf g})$ by $\bf R$ bi-linearity (see Formulas 2.1$(3)$).
\par {\bf 5. Corollary.} {\it The set $\check{\cal M}(\check{\sf g})$ is
the Lie super-algebra over ${\cal A}_r$.}
\par {\bf Proof.} By the construction $\check{\cal M}(\check{\sf g})$ is the ${\cal A}_r$ vector space.
Since ${\cal M}(\check{\sf g})$ is the Lie super-algebra and $\omega
$ satisfies Conditions $(12,14)$, then the multiplication in
$\check{\cal M}(\check{\sf g})$ defined by Formula 4$(1)$ satisfies
Identities 2.1$(1-3)$.
\par {\bf 6. Proposition.} {\it The operator $d_l$ from Definition 4 is the differentiation
of the algebra $\check{\cal M}(\check{\sf g})$.}
\par {\bf Proof.} From the definition of $d_l$ we get
\par $(1)$ $d_l ([P\otimes x, Q\otimes y]_0) = d_l (PQ\otimes [x,y]) $
\par $=(d_l(PQ))\otimes [x,y] = -  \{ ((dP(z)/dz).z^{l+1})Q + P((dQ(z)/dz).z^{l+1}) \} \otimes
[x,y]$
\par $= ((d_lP)Q + P(d_lQ))\otimes [x,y] = [(d_lP)\otimes x, Q\otimes y]_0 + [P\otimes x, (d_lQ)\otimes
y]_0$. By the bi-additivity we get
\par $(2)$ $d_l[a,b]_0 = [d_la,b]_0 + [a,d_lb]_0$, consequently,
\par $d_l[a + \alpha K,b+\beta K] = [d_la,b]_0 + [a,d_lb]_0$
for all $a, b \in {\check{\cal M}}(\check{\sf g})$ and each $\alpha
, \beta \in F(W,{\cal A}_r)$. On the other hand,
\par $[d_la,b] = [d_la,b]_0 + \omega (d_la,b)K$.
\\ Since $ Res (z_0, (d_0d_lP)Q) + Res (z_0, (d_lP)(d_0Q)) = Res
(z_0, d_0((d_lP)Q)) =0$ and symmetrically $ Res (z_0, Pd_0d_lQ) +
Res (z_0, (d_0P)(d_lQ)) = Res (z_0, d_0(Pd_lQ)) =0$, then
\par $(Res (z_0,(d_0d_lP)Q) + Res
(z_0,(d_lP)(d_0Q))).(x|y)$\par $ = (- Res (z_0,(d_lP)(d_0Q)) + Res
(z_0,(d_lP)(d_0Q))).(x|y) =0$ and inevitably
\par $(3)$ $\omega (d_la,b) + \omega (a,d_lb) =0$ for each
$a=P\otimes x$ and $b=Q\otimes y$. Using $\bf R$ linearity of $Res
(z_0, (d_0P)Q).(x|y)$ by $P$, $Q$, $x$ and $y$ and Equations $(1-3)$
we get the statement of this proposition.
\par {\bf 7. Proposition.} {\it The commutator of differentiations
$d_k$ and $d_j$ on the space ${\cal M}(V,W)$ is
\par $(1)$ $[d_k,d_j] = (k-j)d_{j+k}$.}
\par {\bf Proof.} Since $d_kd_jPQ = d_k((d_jP)Q + Pd_jQ) =
(d_kd_jP)Q + (d_jP)(d_kQ) + (d_kP)(d_jQ) + P(d_kd_jQ)$, then
$[d_k,d_j]PQ = (d_kd_j-d_jd_k)PQ = ([d_k,d_j]P)Q + (d_jP)(d_kQ) -
(d_kP)(d_jQ) + (d_kP)(d_jQ) - (d_jP)(d_kQ) + P([d_k,d_j]Q)$
$=([d_k,d_j]P)Q + P([d_k,d_j]Q)$. \par For each term $\alpha z^l$
with a constant $\alpha \in {\cal A}_r$ we have $[d_k,d_j]\alpha z^l
= l(j+l)\alpha z^{j+l+k} - l(k+l)\alpha z^{k+l+j} = (j-k)l\alpha
z^{j+l+k} = (k-j)d_{j+k}\alpha z^l$ on $V\setminus W$ for each $l\in
\bf Z$. We also have \par $(2)$ $d_j \{ \alpha
_{m_1,j}z^{m_1}...\alpha _{m_p,j}z^{m_p} \} _{q(2p)}=$
\par $\sum_{s=1}^p \{
\alpha _{m_1,j}z^{m_1}...\alpha _{m_{s-1},j}z^{m_{s-1}}(\alpha
_{m_s,j}d_lz^{m_s})\alpha _{m_{s+1},j}z^{m_{s+1}}...\alpha
_{m_p,j}z^{m_p} \} _{q(2p)}$. \par By the non commutative analog of
the Stone-Weierstrass theorem (see \S 2 \cite{luoyst,luoyst2}) for
each compact canonical closed subset $J$ in $V\setminus W$ the set
of polynomials over ${\cal A}_r$ is dense in the set of continuous
${\cal A}_r$ valued functions on $J$. For each $z$ differentiable
$f$ function on $J$ there exists a sequence of polynomials $f_n$
converging to $f$ with ${f'}_n$ converging to ${f'}$ uniformly on
$J$ and $J\times \{ h\in {\cal A}_r: ~ |h|\le 1 \} $ respectively.
Applying identities deduced above
\par $(3)$ $[d_k,d_j]PQ = ([d_k,d_j]P)Q + P ([d_k,d_j]Q)$ and
\par $(4)$ $[d_k,d_j]\alpha z^l = (k-j)d_{j+k}\alpha z^l$
to ordered products of terms $\alpha _{m_k,j}z^{m_k}$ with
associators $q(2p)$, where $\alpha _{m_k,j}$ are Cayley-Dickson
constants, $p\in {\bf N}$, $m_1,...,m_p\in \bf Z$, we get that
\par $(4)$ $[d_k,d_j]f = (k-j)d_{k+j}f$ for each $f\in {\cal M}(V,W)$
on the domain $V\setminus W$.
\par {\bf 8. Remark.} From the definition of operators $d_j$ it follows, that
\par $(1)$ $[p \beta d_j, d_k ] = [p d_j,\beta d_k] = p\beta
[d_j,d_k ]$ for each real number $p$ and every Cayley-Dickson number
$\beta $. Certainly
\par $(2)$ $[\alpha d_s + \beta d_j, d_k] = [\alpha d_s,d_k] + [\beta
d_j,d_k]$ and $[d_k, \alpha d_s + \beta d_j] = [d_k, \alpha d_s] +
[d_k, \beta d_j]$ for each $\alpha , \beta \in {\cal A}_r$. Define
the algebra ${\bf \partial } := \bigoplus_{j\in \bf Z} {\cal A}_r
d_j$ of vector fields on ${\cal A}_r \setminus \{ 0 \} $. The non
commutative analog of the Virasoro algebra is ${\cal V}_r := {\bf
\partial } \oplus {\cal A}_rc$, where $c$ denotes the additional
generator and with the multiplication
\par $(3)$ $[d_j + s c, d_k + t c ] = (j-k)d_{j+k} + [(j^3-j)\delta _{j, -k
}/12] c$, where $s, t \in {\cal A}_r$. \par Then it is possible to
define the semidirect product ${\cal V}_r\otimes ^s \check{{\cal
M}}(\check{{\sf g}})$ with the multiplication given by Formulas
$(1-3)$, 4$(1)$ and \par $(4)$ $[c, \check{{\cal M}}(\check{{\sf
g}})] =0$.
\par This algebra may be useful in the quantum field theory over
quaternions and octonions.

\end{document}